\documentclass[paper=A4, DIV=12, USenglish, numbers=noenddot,fontsize=12pt]{scrartcl} 
\usepackage{cmap}		
\usepackage[utf8]{inputenc}	
\usepackage[T1]{fontenc}	
\usepackage[USenglish]{babel}
\usepackage[]{microtype}
\usepackage{booktabs} 

\usepackage{amsmath}
\usepackage{amssymb}
\usepackage{mathtools}		
\usepackage{xfrac}

\usepackage{aliascnt}		
\usepackage{amsthm}

\usepackage{pdflscape}
\usepackage{multirow,bigdelim}

\usepackage{tikz}		
\usetikzlibrary{matrix,arrows,patterns,intersections,calc,decorations.pathmorphing,decorations.markings}

\usepackage[pdftitle={Large-genus asymptotics of saddle connection Siegel--Veech constants},pdfauthor={Anja Randecker},pdfborder={0 0 0}]{hyperref}
\usepackage[figure, table]{hypcap}	

\newtheoremstyle{break}
{} 
{} 
{\itshape} 
{} 
{\bfseries} 
{} 
{\newline} 
{\thmname{#1}\thmnumber{ #2}\thmnote{ (#3)}} 

\newtheoremstyle{breakdef}
{} 
{} 
{} 
{} 
{\bfseries} 
{} 
{\newline} 
{\thmname{#1}\thmnumber{ #2}\thmnote{ (#3)}} 

\newtheoremstyle{remark}
{} 
{} 
{} 
{} 
{\itshape} 
{.} 
{0.5em} 
{\thmname{#1}\thmnumber{ #2}\thmnote{ {\normalfont (}#3{\normalfont )}}} 

\theoremstyle{breakdef}

\theoremstyle{remark}

\newaliascnt{rem}{definition}

\aliascntresetthe{rem}

\newaliascnt{not}{definition}

\aliascntresetthe{not}

\newaliascnt{exa}{definition}

\aliascntresetthe{exa}

\newaliascnt{lem}{definition}
\newtheorem{lem}[lem]{Lemma}
\aliascntresetthe{lem}

\theoremstyle{break}

\newtheorem{thm}{Theorem}

\newaliascnt{prop}{definition}
\newtheorem{prop}[prop]{Proposition}
\aliascntresetthe{prop}

\newaliascnt{cor}{definition}
\newtheorem{cor}[cor]{Corollary}
\aliascntresetthe{cor}

\usepackage[nameinlink,noabbrev,capitalise]{cleveref}
\crefname{rem}{Remark}{Remarks}
\crefname{definition}{Definition}{Definitions}
\crefname{exa}{Example}{Examples}
\crefname{lem}{Lemma}{Lemmas}
\crefname{thm}{Theorem}{Theorems}
\crefname{prop}{Proposition}{Propositions}
\crefname{cor}{Corollary}{Corollaries}
\crefname{section}{Section}{Sections}
\crefname{subsection}{Subsection}{Subsections}
\crefname{figure}{Figure}{Figures}

\usepackage{enumitem}

\numberwithin{figure}{section}			

\DeclareMathOperator{\SL}{SL}
\DeclareMathOperator{\vol}{vol}

\newcommand{\saddleconnection}{
	\begin{tikzpicture}[baseline=-0.1cm]
		\fill (0,0) circle (0.05cm);
		\fill (0.5,0) circle (0.05cm);
		\draw (0,0) -- (0.5,0);
	\end{tikzpicture}	
}

\newcommand{\scloop}{
	\begin{tikzpicture}[baseline=-0.1cm]
		\fill (0,0) circle (0.05cm);
		\draw (0,0) arc (180:-180:0.15);
	\end{tikzpicture}	
}

\newcommand{\cylinder}{
	\begin{tikzpicture}[baseline=-0.1cm]
		\fill (-0.3,0) circle (0.05cm);
		\fill (0.5,0) circle (0.05cm);
		\draw (0,0) arc (0:360:0.15);
		\draw (0.5,0) arc (0:90:0.15);
		\draw (0.5,0) arc (0:-90:0.15);
		\draw(-0.15,-0.15) -- (0.35,-0.15);
		\draw(-0.15,0.15) -- (0.35,0.15);		
	\end{tikzpicture}	
}

\newcommand{\loopcylinder}{
	\begin{tikzpicture}[baseline=-0.1cm]
		\fill (0,0) circle (0.05cm);
		\draw (0,0) to[out=120,in=180] (0,0.3) to[out=0,in=60] (0,0);
		\draw (0,0) to[out=120,in=180] (0,-0.3) to[out=0,in=60] (0,0);
		\draw (0,0.3) to[out=10,in=90] (0.4,0) to[out=-90,in=-10] (0,-0.3);
		\draw (0,0) to[out=45, in=90] (0.2,0) to[out=-90,in=-45] (0,0);
	\end{tikzpicture}	
}

\newcommand{\labelled}{\textnormal{labelled}}

\newcommand{\fixedzeros}{\textnormal{fixed zeros}}
\newcommand{\anyzeros}{\textnormal{any zeros}}
\newcommand{\fixedzero}{\textnormal{fixed zero}}
\newcommand{\anyzero}{\textnormal{any zero}}
\newcommand{\onefixedzero}{\textnormal{one fixed zero}}
\newcommand{\dominant}{\textnormal{dominant}}
\newcommand{\nondominant}{\textnormal{non-dominant}}
\newcommand{\any}{\textnormal{any}}
\newcommand{\hyp}{\textnormal{hyp}}
\newcommand{\odd}{\textnormal{odd}}
\newcommand{\even}{\textnormal{even}}
\newcommand{\nonhyp}{\textnormal{non-hyp}}
\newcommand{\oddeven}{\textnormal{odd/even}}
\newcommand{\oen}{\textnormal{odd/even/non-hyp}}
\newcommand{\gminusone}{g\!-\!1}
\newcommand{\twogminustwo}{2g\!-\!2}

\newcommand{\param}{{\mathchoice{\mkern1mu\mbox{\raise2.2pt\hbox{$
					\centerdot$}}
			\mkern1mu}{\mkern1mu\mbox{\raise2.2pt\hbox{$\centerdot$}}\mkern1mu}{
			\mkern1.5mu\centerdot\mkern1.5mu}{\mkern1.5mu\centerdot\mkern1.5mu}}
}

\title{Large-genus asymptotics \\of saddle connection \\Siegel--Veech constants}
\author{Anja Randecker}
\date{August 14, 2025}

\begin{document}
	
	\maketitle
	
	\begin{abstract}
		Siegel--Veech constants are powerful tools for counting saddle connections on a translation surface. Their computation can be involved, most famously with recursive formulas that use intricate combinatorics or intersection theory. From these formulas, asymptotics of Siegel--Veech constants for growing genus can be extracted. We extend the known asymptotics to all strata and to all multiplicities of saddle-connections between distinct zeros and of loops.
	\end{abstract}
	
	For closed hyperbolic surfaces, Huber has shown that the number of closed geodesics up to a given length grows exponentially in the length bound~\cite{huber_59,huber_61}. Mirzakhani refined this result by considering only simple closed curves: For a hyperbolic surface of genus $g$ with $n$ punctures, the number of simple closed geodesics up to length $L$ grows as~$L^{6g-6+2n}$ and the constant in the growth can be explicitly determined \cite{mirzakhani_08}.
	
	In this article, we equip a topological surface with marked points not with a hyperbolic metric with punctures but with a flat metric outside of the marked points. Then we can still count closed geodesics (which leads to cylinder counting) but also geodesics from a marked point to another marked point or itself. The latter case can be interpreted as geodesics relative to the marked points and
	these relative geodesics are called
	\emph{saddle connections}. Their number grows as $L^2$; in particular, the growth rate does not depend any more on the genus \cite{masur_88,masur_90}.
	However, the constant in the growth does depend on the genus but also on more properties, encoded in the choice of a connected component of a stratum (see \cref{sec:translation_surfaces} for definitions).
	In general, it does not depend on the flat metric of the surface, though, as Eskin and Masur have shown.
	
	\begin{thm}[Generic number of saddle connections \cite{eskin_masur_01}]
		Let $\mathcal{H}$ be a connected component of a stratum of unit-area translation surfaces. Then there exists a constant $c(\mathcal{H})$ such that for almost every $(X,\omega) \in \mathcal{H}$, we have
		\begin{equation*}
			\lim_{L \to \infty} \frac{ \left| \{\textnormal{saddle connections on } (X,\omega) \textnormal{ of length at most } L \} \right| }{\pi L^2}
			= c(\mathcal{H})
			.
		\end{equation*}
	\end{thm}
	
	The constant $c(\mathcal{H})$ is called the \emph{saddle connection Siegel--Veech constant} of $\mathcal{H}$ and fascinatingly, it first appeared as a constant in a slightly different counting problem.
	
	\begin{thm}[Siegel--Veech formula for saddle connections {\cite{veech_98}}] 
		\label{thm:Siegel-Veech_formula}
		Let $\mathcal{H}$ be a connected component of a stratum of unit-area translation surfaces. Then 
		for every integrable function $f \colon \mathbb{R}^2 \to \mathbb{R}$ with compact support, we have
		\begin{equation*}
			\int\limits_{\mathcal{H}} \sum_{\substack{v \textnormal{ holonomy vector} \\ \textnormal{of a saddle connection}}} f(v) = c(\mathcal{H}) \cdot \int\limits_{\mathbb{R}^2} f \cdot \vol(\mathcal{H})
			.
		\end{equation*}
	\end{thm}
	
	In fact, both theorems are much stronger than stated here. In particular, we can count not only saddle connections but saddle connections with additional properties, such as having the same start and end point (see \cref{sec:Siegel-Veech_constants} for a more general statement).
	
	There is a deep interplay between Siegel--Veech constants and volumes of strata that already shines through in \cref{thm:Siegel-Veech_formula}. 
	Therefore, the progress on determining volumes of strata has gone hand in hand with the progress on Siegel--Veech constants in the last decade.
	
	In this article, we give an overview of the large-genus asymptotics of Siegel--Veech constants.
	Specifically, we extend the known asymptotics for saddle connections between distinct zeros to non-connected strata and explicitly determine the Siegel--Veech constants for higher multiplicities. We also bound the Siegel--Veech constants for loops for any multiplicity and therefore prove that the Siegel--Veech constant for any multiplicity is given by the one for multiplicity $1$.
	In fact, often times in applications, it is only necessary to known that certain types of saddle connections are relatively unlikely, so for some of the values, we settle for upper bounds where too many cases would have been to be distinguished.
	
	There are also several other versions of Siegel--Veech constants beyond the scope of this article, notably \emph{cylinder Siegel--Veech constants} and more generally \emph{area Siegel--Veech constants}. The latter are particularly relevant as they are related to Lyapunov exponents~\cite{eskin_kontsevich_zorich_14} and there equally has been a lot of progress, for example in \cite{chen_moeller_zagier_18,sauvaget_18,aggarwal_19,chen_moeller_sauvaget_zagier_20,aulicino_calderon_matheus_salter_schmoll_24}.
	Furthermore, there is an even more involved and still less developed, parallel story for large-genus asymptotics for saddle connection Siegel--Veech constants for quadratic differentials. Some conjectures and first results can be found in \cite{delecroix_goujard_zograf_zorich_21,chen_moeller_sauvaget_23,aggarwal_delecroix_goujard_zograf_zorich_20,aggarwal_21,yang_zagier_zhang_20,duryev_goujard_yakovlev_25}.

	\paragraph{Organization.}
	In \cref{sec:translation_surfaces}, we start with some vocabulary and facts about translation surfaces, in particular on the volume of strata.
	The method of Siegel--Veech constants is explained in \cref{sec:Siegel-Veech_constants}.
	We recall the strategy from \cite{eskin_masur_zorich_03} and their formulas for Siegel--Veech constants for connected strata in \cref{sec:EMZ} and extensions to non-connected strata in \cref{sec:non-connected_strata}.
	In \cref{sec:lemmas}, we provide several lemmas that are used in the subsequent calculations.
	The values of Siegel--Veech constants are then recalled (where known) and determined in \cref{sec:known_values_non-loops,sec:known_values_loops}, with the special cases of $\mathcal{H}^\hyp(2g-2)$ and~$\mathcal{H}^\hyp(g-1,g-1)$ treated in \cref{sec:values_minimal,sec:values_g-1}.
	The article is summarized in a table of values in~\cref{sec:table}.

	\paragraph{Acknowledgements.}
	This work was funded by the Deutsche Forschungsgemeinschaft (DFG, German Research Foundation) -- 441856315. The author thanks Howie Masur, Kasra Rafi, and Anton Zorich for several discussions about Siegel--Veech constants in the last years. Special thanks go to Maurice Reichert for discussions and thorough proofreading of a previous version of this article.

	\section[Strata of translation surfaces and Masur--Veech measure]{{\textls[-35]{\spaceskip 0.2em Strata of translation surfaces and Masur--Veech measure}}} \label{sec:translation_surfaces}
	
	A \emph{translation surface} $(X,\omega)$ is a closed surface $X$, equipped with
	a non-zero holomorphic Abelian differential $\omega$, equivalently a holomorphic $1$--form $\omega$.
	Integrating the $1$--form yields a metric on the surface which is flat except in the zeros. Therefore, the zeros of the~$1$--form are also called \emph{singularities} of the translation structure and for a zero of order~$m$, we have an angle of $(m+1) \cdot 2 \pi$ around the singularity.
	Hence, for a translation surface of genus $g$ with zeros of orders $m_1, \ldots, m_\ell$, it holds $2g-2 = \sum_{n=1}^\ell m_n$.
	
	Every closed geodesic on a translation surface comes in a family, yielding a \emph{cylinder}. Any cylinder is bounded by a collection of \emph{saddle connections}, that is, of geodesic segments that start and end at a zero and do not contain a zero in their interior. Saddle connections do not necessarily have to start and end at the same zero but if they do, we call them~\emph{loops}.
	Saddle connections can be
	lifted to the universal cover and developed into the plane; the corresponding vector in $\mathbb{R}^2$ is called the \emph{holonomy vector} of a saddle connection.
	The set of all holonomy vectors of a given translation surface is discrete in~$\mathbb{R}^2$.
	
	\bigskip
	
	Using the flat structure of a translation surface, it is easy to see that there exists an action of $\SL(2,\mathbb{R})$ on the space of all translation surfaces. Formally, we can define it by~$A \cdot (X, \omega) = (X, A \cdot \omega)$ for every $A \in \SL(2,\mathbb{R})$.
	
	\bigskip
	
	For a fixed genus $g$ and a partition $2g-2 = m_1 + \ldots + m_\ell$, the \emph{stratum} $\mathcal{H}(m_1, \ldots, m_\ell)$ is the space of all translation surfaces of genus $g$ with zeros of order $m_1, \ldots, m_\ell$.
	Using an integer basis of relative homology, that is, a convenient set of saddle connections, we can show that each stratum locally is homeomorphic to $\mathbb{R}^{2(2g+\ell-1)}$. Pulling back the Lebesgue measure on~$\mathbb{R}^{2(2g+\ell-1)}$ to $\mathcal{H}(m_1, \ldots, m_\ell)$, we obtain an $\SL(2,\mathbb{R})$--invariant measure $\mu_\mathrm{MV}$.
	The measure $\mu_\mathrm{MV}$ is infinite but 
	we can use it to define another measure $\mu$ on the codimension--$1$ subset of translation surfaces in $\mathcal{H}(m_1, \ldots, m_\ell)$ with area $1$.
	Specifically, for a set $A$ of unit-area translation surfaces, we define $\mu(A)$ to be the normalized measure of the cone of $A$ in the whole stratum, that is,
	\begin{equation*}
		\mu (A) \coloneqq 2(2g+\ell -1) \cdot \mu_\mathrm{MV}( \{ r \cdot (X,\omega) : r \in (0,1], (X,\omega) \in A)\} )
		.
	\end{equation*}
	Then the \emph{Masur--Veech measure}\,\footnote{We give in to the usual ambiguity in notation by referring to $\mu$ as well as to $\mu_\mathrm{MV}$ as \emph{Masur--Veech measure} and to $\mathcal{H}(m_1, \ldots, m_\ell)$ as well as to its codimension--$1$ subset of unit-area translation surfaces as \emph{stratum}.} $\mu$ is a finite, $\SL(2,\mathbb{R})$--invariant measure which is ergodic on any connected component of the stratum \cite{masur_82,veech_82}.
	For the following, it will be relevant that by convention, the zeros in a stratum are always considered to be labelled when calculating volumes.
	
	\bigskip
	
	A stratum does not have to be connected, it can consist of up to three connected components \cite{kontsevich_zorich_03}. More specifically, for genus $g \geq 4$, the strata $\mathcal{H}(2g-2)$ and $\mathcal{H}(g-1, g-1)$ have a \emph{hyperelliptic component} and a stratum where all zeros are of even order has an \emph{even} and an \emph{odd connected component}\,\footnote{The terms ``even'' and ``odd'' refer to the parity of the so-called spin structure (see \cref{sec:non-connected_strata}).}. All other strata for $g\geq 4$ are~connected.
	
	Even for a connected stratum, it is difficult to calculate the Masur--Veech volume. Eskin and Okounkov developed an algorithm to calculate the volumes through intricate counting of torus covers~\cite{eskin_okounkov_01}
	and based on this, Eskin and Zorich conjectured in~2003 (but published in 2015) large-genus asymptotics for the volumes \cite{eskin_zorich_15}. After the most relevant cases of the principal stratum \cite{chen_moeller_zagier_18} and the minimal stratum \cite{sauvaget_18} have been proven with the use of intersection theory, Aggarwal in \cite{aggarwal_20} settled the conjecture for the rest of the cases by following the strategy in \cite{eskin_okounkov_01} and showing that~indeed
	\begin{equation} \label{eq:volume}
		\mu(\mathcal{H}( m_1, \ldots, m_\ell) ) 
		= \frac{4}{\prod_{n=1}^\ell (m_n + 1)} \cdot \left( 1+ O\left(\frac1g \right) \right)
	\end{equation}
	where $O(\sfrac1g)$ is always meant with respect to $g \to \infty$ from now on.
	
	For our purposes, the error term $(1 + O(\sfrac1g ))$ for the volume estimates will be sufficient but more precise error terms can be found in the work of Chen, Möller, Sauvaget, and Zagier~\cite[Theorem 1.5]{chen_moeller_sauvaget_zagier_20} and of Sauvaget \cite{sauvaget_21}. Furthermore, Chen, Möller, Sauvaget, and Zagier have also proven another conjecture from~\cite{eskin_masur_01} that the volumes of the odd and even components are comparable \cite[Theorem 1.6]{chen_moeller_sauvaget_zagier_20}, that is,
	\begin{equation} \label{eq:volume_odd_even}
		\frac{\mu(\mathcal{H}^\odd )}{\mu(\mathcal{H}^\even )}
		= 1+ O\left( \frac1g \right)
	\end{equation}
	under a technical assumption which was subsequently proven by Costantini, Möller, and Zachhuber in \cite{costantini_moeller_zachhuber_24}.
	
	Already before the recent leap in results on Masur--Veech volumes, Athreya, Eskin, and Zorich have shown that the volume of any hyperelliptic component is negligible in comparison to the whole stratum.
	Specifically, they calculated the volumes of certain strata of meromorphic \emph{quadratic differentials} with zeros and at most simple poles (which are not labelled).
	As every hyperelliptic translation surface with a unique zero of order~$2g-2$ can be obtained as a double cover of a quadratic differential with one zero of order~$2g-3$ and $2g+1$ simple poles,
	the volume of $\mathcal{H}^\hyp (2g-2)$ corresponds to the volume of $\mathcal{Q}(2g-3, (-1)^{2g+1})$ up to a factor of $(2g+1)!$ which corresponds to the possible choices of labelling the simple poles in $\mathcal{Q}(2g-3, (-1)^{2g+1})$ (see~\cite[Remark~1.2]{athreya_eskin_zorich_16}). With this correspondence, \cite[Theorem 1]{athreya_eskin_zorich_16} yields
	\begin{equation} \label{eq:volume_hyperelliptic_minimal}
		\mu(\mathcal{H}^\hyp (2g-2))
		= \frac{2 \cdot \pi^{2g}}{(2g+1)!} \cdot \frac{(2g-3)!!}{(2g-2)!!}
	\end{equation}
	where $n!! = n (n-2) (n-4) \cdots$ is the product of all even resp.\ odd numbers up to $n$.
	
	Similarly, the volume of $\mathcal{H}^\hyp (g-1, g-1)$ corresponds to the volume of the stratum $\mathcal{Q}(2g-2, (-1)^{2g+2})$ of quadratic differentials up to a factor of $\sfrac{(2g+2)!}{2}$ for labelling the simple poles in $\mathcal{Q}(2g-2, (-1)^{2g+2})$ and forgetting the names of the two zeros in~$\mathcal{H}^\hyp (g-1,g-1)$. Therefore, \cite[Theorem 1]{athreya_eskin_zorich_16} yields
	\begin{equation} \label{eq:volume_hyperelliptic_two}
		\mu(\mathcal{H}^\hyp (g-1,g-1))
		= \frac{8 \cdot \pi^{2g}}{(2g+2)!} \cdot \frac{(2g-2)!!}{(2g-1)!!}
		.
	\end{equation}
	With Stirling's approximation and $\frac{(2g)!!}{(2g-1)!!} = \sqrt{\pi g} \cdot (1 + O (\sfrac1g))$, we can write these~as
	\begin{align}
		\mu(\mathcal{H}^\hyp (2g-2))
		& = \frac{1}{2\pi} \cdot \left( \frac{\pi e}{2} \right)^{2g} \cdot \frac{1}{g^{2g+2}} \cdot \left( 1 + O \left( \frac1g \right) \right) \label{eq:volume_hyperelliptic_minimal_Stirling} \\
		\mu(\mathcal{H}^\hyp (g-1,g-1))
		& = \frac12 \cdot \left( \frac{\pi e}{2} \right)^{2g} \cdot \frac{1}{g^{2g+3}} \cdot \left( 1 + O\left( \frac1g \right)\right)  \label{eq:volume_hyperelliptic_two_Stirling}
	\end{align}
	which shows that the volumes are negligible compared to the volume of the whole stratum
	which is $\frac{4}{2g-1} \cdot (1+ O(\sfrac1g))$ for $\mathcal{H}(2g-2)$ and $\frac{4}{g^2} \cdot (1+ O(\sfrac1g))$ for $\mathcal{H}(g-1, g-1)$.

	\section{Siegel--Veech constants} \label{sec:Siegel-Veech_constants}
	
	In \cite{siegel_45}, Siegel famously studied averages of functions in the setup of unimodular lattices.
	Recall that the space of unimodular lattices in $\mathbb{R}^d$ can be identified with $\SL(d,\mathbb{R})/\SL(d,\mathbb{Z})$.
	For a unimodular lattice $\Lambda \subseteq \mathbb{R}^d$ and an integrable function $f \colon \mathbb{R}^d \to \mathbb{R}$ with compact support, the \emph{Siegel transform} of~$f$ is given by $\widehat{f}(\Lambda) = \sum_{v \in \Lambda \setminus \{0\}} f(v)$ and the \emph{Siegel formula} states $\int_{\SL(d,\mathbb{R})/ \SL(d,\mathbb{Z}) } \widehat{f} = \int_{\mathbb{R}^d} f$ where the measure on the left-hand side is the Haar measure on
	$\SL(d,\mathbb{R})/\SL(d,\mathbb{Z})$ 
	and the measure on the right-hand side is the Lebesgue~measure.
	
	According to the acknowledgements in \cite{veech_98} (and the corresponding anecdote), Veech attended a talk of Margulis at CIRM in Marseille that referenced the Siegel formula which inspired Veech to formulate the following analogue statement for translation surfaces. Fix a connected component $\mathcal{H}$ of a stratum of unit-area translation surfaces. Furthermore, fix a configuration $\mathcal{C}$ of saddle connections, that is, fix the multiplicity of the saddle connection, the order of the zeros that it connects, and so on	(see \cref{sec:EMZ} for more details on configurations). Let $V_{\mathcal{C}}(X, \omega)$ denote the set of holonomy vectors of saddle connections on~$(X,\omega) \in \mathcal{H}$ which are in the configuration $\mathcal{C}$. The corresponding \emph{Siegel--Veech transform} of an integrable function $f \colon \mathbb{R}^2 \to \mathbb{R}$
	is $\widehat{f}_{\mathcal{C}} \colon \mathcal{H} \to \mathbb{R}, (X, \omega) \mapsto \sum_{v \in V_{\mathcal{C} (X, \omega)}} f(v)$.
	
	\begin{thm}[Siegel--Veech formula {\cite{veech_98}}]
		Let $\mathcal{H}$ be the connected component of a stratum of unit-area translation surfaces and~$\mathcal{C}$ a configuration of saddle connections. Then there exists a constant~$c(\mathcal{C}, \mathcal{H})$ such that for every integrable $f \colon \mathbb{R}^2 \to \mathbb{R}$ with compact support, we have
		\begin{equation*}
			\frac{\int_{\mathcal{H}} \widehat{f}_{\mathcal{C}}}{\mu(\mathcal{H})} = c(\mathcal{C}, \mathcal{H}) \cdot \int_{\mathbb{R}^2} f
		\end{equation*}
		where the measure on the left-hand side is the Masur--Veech measure on $\mathcal{H}$ and the measure on the right-hand side is the Lebesgue measure.
	\end{thm}
	
	A key argument of the proof is that the linear functional $f \mapsto \int_{\mathcal{H}} \widehat{f}_{\mathcal{C}}$ is $\SL(2,\mathbb{R})$--invariant and therefore $\int_\mathcal{H} \widehat{f}_\mathcal{C}$ has to be a multiple of $\int_{\mathbb{R}^2} f$.
	
	By choosing $f$ to be the indicator function of the ball $B(0,L)$ of radius~$L$, we obtain that the average number of saddle connections in the configuration $\mathcal{C}$ of length up to~$L$ for translation surfaces in~$\mathcal{H}$ is $c(\mathcal{C}, \mathcal{H}) \cdot \mu(\mathcal{H}) \cdot \pi L^2$. 
	
	What makes the Siegel--Veech constant $c(\mathcal{C}, \mathcal{H})$ even more interesting is that it does not only govern the average number of saddle connections in $\mathcal{C}$ over the stratum $\mathcal{H}$ but also generically the growth rate of saddle connections on a given translation surface in~$\mathcal{H}$.
	
	\begin{thm}[Siegel--Veech constant as constant in quadratic growth \cite{eskin_masur_01}]
		Let $\mathcal{H}$ be the connected component of a stratum of unit-area translation surfaces and $\mathcal{C}$ a configuration of saddle connections with Siegel--Veech constant $c(\mathcal{C}, \mathcal{H})$. Then for almost every $(X,\omega) \in \mathcal{H}$, we have
		\begin{equation*}
			\lim_{L \to \infty} \frac{ \left| V_{\mathcal{C}}(X, \omega) \cap B(0, L) \right| }{\pi L^2}
			= c(\mathcal{C}, \mathcal{H})
			.
		\end{equation*}
	\end{thm}
	
	Both of these results can be stated more generally than for configurations of saddle connections, to include for example cylinders or other $\SL(2,\mathbb{R})$--invariant setups, even if they do not have a direct geometric interpretation. Nevertheless, we focus on configurations of saddle connections in the following.
	
	\bigskip
	
	The Siegel--Veech formula does not immediately provide a way to calculate the Siegel--Veech constant for a given connected component of a stratum and a configuration. However, Eskin, Masur, and Zorich developed a powerful technique to describe Siegel-Veech constants recursively in \cite{eskin_masur_zorich_03}.
	The recursion can get quite involved but many values of Siegel--Veech constants for small genera have been determined this way~in~\cite{eskin_masur_zorich_03}.
	
	However, we will put our focus on a different aspect: the asymptotic behaviour of Siegel--Veech constants for growing genus. For this, large-genus asymptotics of Masur--Veech volumes are relevant and hence the progress went parallely to the calculation of volumes. In particular,  Chen, Möller, and Zagier \cite{chen_moeller_zagier_18} computed large-genus asymptotics of saddle connection Siegel--Veech constants for principal strata, as well as Zorich in an appendix to \cite{aggarwal_20} and Aggarwal~\cite{aggarwal_19} for connected strata 
	and Chen, Möller, Sauvaget, and Zagier \cite{chen_moeller_sauvaget_zagier_20} for connected components of strata
	(see \cref{sec:known_values_non-loops,sec:known_values_loops} for detailed statements).

	\section{Technique and formulas from Eskin--Masur--Zorich} \label{sec:EMZ}
	
	To perform explicit calculations later, we recall here the setup and the technique from~\cite{eskin_masur_zorich_03} as well as the two main formulas for connected strata.
	
	Consider a saddle connection $\gamma$ on a translation surface. The \emph{multiplicity} of $\gamma$ is the maximal number of saddle connections $\gamma_1 = \gamma, \gamma_2, \ldots, \gamma_p$ which are pairwise homologous. On a generic translation surface, being homologous is equivalent to having the same holonomy vector \cite[Proposition 3.1]{eskin_masur_zorich_03}. Furthermore, if two saddle connections between distinct zeros are homologous, they both start at the same zero and both end at the same zero.
	For loops, this is not necessarily true, for instance when two loops bound a common cylinder, they do not need to start at the same zero.
	Note, however, that two loops which bound a common cylinder are counted only once for multiplicity.
	
	The strategy in \cite{eskin_masur_zorich_03} is to cut a translation surface with a set of $p$ homologous saddle connections (which also bound~$q$ cylinders) into $p$ subsurfaces with boundary and~$q$ cylinders. The boundary components of each of the $p$ subsurfaces are then manipulated and identified to obtain translation surfaces from lower-dimensional strata $\mathcal{H}_i$.
	
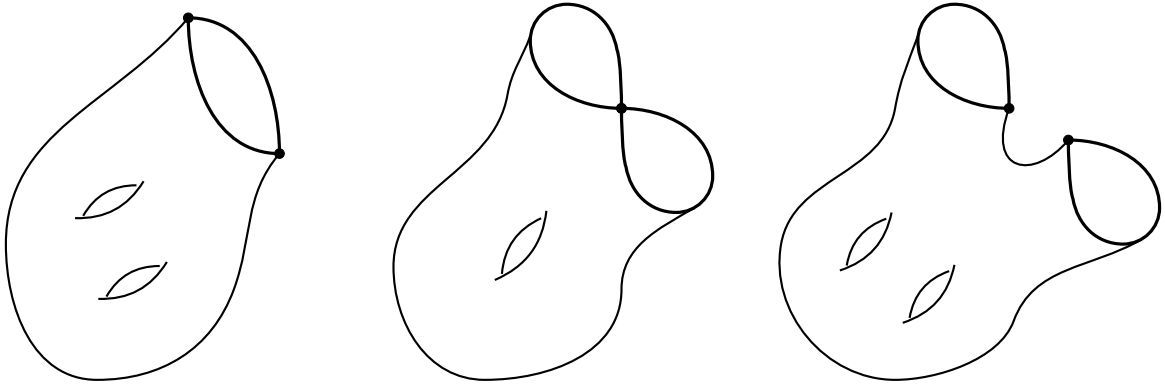
\begin{figure}
		\begin{center}
			\begin{tikzpicture}[scale=0.6]
				\newcommand\hole{
					\draw[very thick] (0,0) to[out=180,in=-90]
					(-2,1.5) to[out=90,in=180]
					(-1.2,2.3) to[out=0,in=100]
					(-0.1,1.3) to[out=-70,in=90,looseness=0.5] (0,0);
					\draw[fill] (0,0) circle (3pt);
				}
				
				\newcommand\genus{
					\draw[thick] ((1,0.05) to[bend left]  (2.5,0.05);
					\draw[thick] ((0.8,0.1) to[bend right] (2.7,0.05);
				}
				
				\begin{scope}[xshift=-8.5cm]
					\draw[very thick] (-1,2) to[out=0,in=90]
					(1,-1) to[out=180,in=-90] (-1,2);
					\draw[fill] (-1,2) circle (3pt);
					\draw[fill] (1,-1) circle (3pt);
					
					\draw[thick] (-1,2) to[out=-130,in=90]
					(-5,-3) to[out=-90, in=180]
					(-3,-6) to[out=0, in=-110]
					(0,-4) to[out=70, in=-130]
					(1,-1);
					
					\begin{scope}[scale=0.9, rotate=30, xshift=-5.5cm, yshift=-0.5cm]
						\genus
					\end{scope}
					\begin{scope}[scale=0.9, rotate=30, xshift=-6cm, yshift=-2.5cm]
						\genus
					\end{scope}
				\end{scope}
				
				\hole
				\begin{scope}[rotate=180]
					\hole
				\end{scope}
				\draw[thick] (-2,1.6) to[out=-110,in=80]
				(-2.5,0.3) to[out=-100, in=90]
				(-5,-3.5) to[out=-90, in=180]
				(-3,-6) to[out=0, in=-90]
				(0,-4) to[out=90, in=-150]
				(1.6,-2.2);
				\begin{scope}[rotate=55, xshift=-5.5cm]
					\genus
				\end{scope}
				
				\begin{scope}[xshift=8.5cm]
					\hole
					\begin{scope}[rotate=180, xshift=-1.3cm, yshift=0.7cm]
						\hole
					\end{scope}
					\draw[thick] (0,0) .. controls +(-110:1.5cm) and +(-130:1.2cm) .. (1.3,-0.7);
					\draw[thick] (-2,1.6) to[out=-110,in=80]
					(-2.5,0) to[out=-100, in=80]
					(-5,-3) to[out=-100, in=180]
					(-2.5,-6) to[out=0, in=-110, looseness=0.8]
					(0.1,-4.7) to[out=70, in=-150]
					(2.9,-2.9);
					\begin{scope}[scale=0.9, rotate=50, xshift=-6.5cm, yshift=0.5cm]
						\genus
					\end{scope}
					\begin{scope}[scale=0.9, rotate=50, xshift=-6.5cm, yshift=-1.5cm]
						\genus
					\end{scope}
				\end{scope}
			\end{tikzpicture}
			\caption{The three types of surfaces with boundary: obtained by a slit construction, by a figure-eight construction, and by a two-hole construction (from left to right).}
			\label{fig:constructions}
		\end{center}
	\end{figure}

	In fact, it is more convenient to describe the process reversely: Consider one translation surface from each stratum $\mathcal{H}_i$ and perform one of the following constructions, depending on the configuration that we want to realize.
	\begin{itemize}
		\item By a \emph{slit construction},\footnote{A special case is also called ``breaking up a zero'' in \cite[Section 8.1]{eskin_masur_zorich_03}.} we obtain from $S \in \mathcal{H}_i = \mathcal{H}(a_i, \ldots)$ a surface with one boundary component as on the left of \cref{fig:constructions}.
		The surface with boundary has the same zeros as $S$ except that the distinguished zero of order $a_i$ is now  split up into two points which are connected by 
		two geodesic segments with the same holonomy vector that form the boundary of the surface. The angles at these two points are~$(a_i' +1) \cdot 2 \pi$ and $(a_i'' +1) \cdot 2 \pi$ where $a_i' + a_i '' = a_i$.
		\item By a \emph{figure-eight construction}, we obtain from $S \in \mathcal{H}_i = \mathcal{H}(a_i, \ldots)$ a surface with one boundary component as in the middle of \cref{fig:constructions}.
		The surface with boundary has the same zeros as $S$ except that the zero of order $a_i$ is now a point on the boundary, dividing it into two geodesic segments with the same holonomy vector. The angles between the two geodesic segments are $(a_i' +1) \cdot 2 \pi$ and $(a_i'' +1) \cdot 2 \pi$ where $a_i' + a_i '' = a_i$.
		\item By a \emph{two-hole construction},\footnote{This construction consists of ``creating a pair of holes'' and then ``transporting a hole'' as described in \cite[Sections 11 and 12]{eskin_masur_zorich_03}. It is also closely related to their ``parallelogram construction''.} we obtain from $S \in \mathcal{H}_k = \mathcal{H}(b_k', b_k'', \ldots)$ a surface with two boundary components as on the right of \cref{fig:constructions}.
		The surface with boundary has the same zeros as $S$ except that the two zeros of order $b_k'$ and $b_k''$ are now points on each of the boundary components. The angles at these two points are~$(b_k' + \sfrac32) \cdot 2 \pi$ and $(b_k'' + \sfrac32) \cdot 2 \pi$.
	\end{itemize}
	
	The boundaries of the surfaces that we obtain are then glued, either directly between two surfaces with boundary or with one of the $q$ cylinders in between.
	
	Note that slit constructions can only be combined with other slit constructions, whereas figure-eight constructions and two-hole constructions can be combined and can also be combined with cylinders. The first case yields saddle connections between distinct zeros whereas the latter yields loops.
	
	In the case of saddle connections between distinct zeros, 
	gluing the subsurfaces with boundaries yields two new zeros whose order can be determined by the data $a_i = a_i +, a_i''$ as $m_1 = \sum_{i=1}^p a_i' +p-1$ and $m_2 = \sum_{i=1}^p a_i'' +p-1$. The genus $g$ of surfaces from $\mathcal{H}$ can be determined by the genera~$g_i$ of surfaces from $\mathcal{H}_i$ as $g = \sum_{i=1}^p g_i$.
	
	In the case of loops, the number of new zeros is the number of two-hole constructions and cylinders and the order of the new zeros can be calculated from the information on the gluing and the~$a_i$ and~$(b_k', b_k'')$. As the gluing produces an additional genus, here we have~$g =\sum_{i=1}^p g_i +1$.
	
	\bigskip

	We describe a \emph{configuration} $\mathcal{C}$ of a given saddle connection $\gamma$ by specifying part or all of the following data: 
	\begin{itemize}
		\item the multiplicity of $\gamma$,
		\item the order $m_1$ of the zero $z_1$ where $\gamma$ starts and the order $m_2$ of the zero $z_2$ where $\gamma$ ends,
		\item whether $z_1 = z_2$ or not, that is, whether $\gamma$ is a loop or not,
		\item whether the zeros $z_1$ and $z_2$ are fixed, or not,
		\item the strata $\mathcal{H}_1, \ldots, \mathcal{H}_p$, up to cyclic permutation,
		\item whether the zeros in $\mathcal{H}_1, \ldots, \mathcal{H}_p$ are \emph{labelled}, that is, whether we specify which zero from $\mathcal{H}$ is contained in which of the $\mathcal{H}_i$,
		\item if $\gamma$ is a loop, on which $\mathcal{H}_1, \ldots, \mathcal{H}_p$ we perform figure-eight constructions and on which two-hole constructions,
		\item the partitions $a_i = a_i' + a_i''$,
		\item if $\gamma$ is a loop, the number $q$ of saddle connections in the homology class of $\gamma$ that bound a cylinder.
	\end{itemize}
	Note that the Siegel--Veech constants for a fixed multiplicity $p \geq 2$ are for the count of homology classes of saddle connections, not of the saddle connections themselves.
	In the following, we will not refer to Siegel--Veech constants as $c(\mathcal{C}, \mathcal{H})$ but give them more telling (but less consistent) names.
	
	\bigskip
	
	Although the order in which the subsurfaces obtained from the $\mathcal{H}_i$ are glued matters, up to cyclic permutation, there can be additional symmetries. For a given configuration, we define the two symmetry groups $\Gamma_-$ and $\Gamma$. If for a given ordered set of strata $\mathcal{H}_1, \ldots, \mathcal{H}_p$ and the data of $a_i = a_i' + a_i''$ or $(b_k', b_k'')$ for each $\mathcal{H}_i$ or $\mathcal{H}_k$, the data, read backwards (including exchanging of $a_i'$ and $a_i''$ and of $b_k'$ and $b_k''$), is the same as the original data up to cyclic permutation, we define that $\Gamma_-$ is $\mathbb{Z}/2\mathbb{Z}$. Otherwise,~$\Gamma_-$ is the trivial group. Furthermore, $\Gamma$ is the group of translational symmetries of the~data, that is, it describes whether we can cyclically permute the $\mathcal{H}_1, \ldots, \mathcal{H}_p$. These symmetry groups can differ, depending on whether we consider labelled or unlabelled configurations.
	
	\bigskip
	
	As every translation surface in $\mathcal{H}$ with saddle connections in a given configuration can be obtained by the gluing process described above, the Siegel--Veech constant of a configuration can be expressed by combinatorial factors and ratios of volumes of $\mathcal{H}$ and~$\mathcal{H}_1 \times \ldots \times \mathcal{H}_p$.
	To adjust for the different normalization factors in the volumes of the $\mathcal{H}_1, \ldots, \mathcal{H}_p$, we also need to take into account the dimensions of $\mathcal{H}_1, \ldots, \mathcal{H}_p$. For this, let $d$ be the real dimension of $\mathcal{H}$ and~$d_i$ the real dimension of $\mathcal{H}_i$.
	We want to emphasize here that for any $\mathcal{H}_i$ or $\mathcal{H}_k$,~$a_i = 0$ or~$b_k' = 0$ and/or $b_k'' = 0$ gives rise to a marked point and has to be considered when determining the dimension. It holds $d = \sum_{i=1}^p d_i + 2q + 2$ where $q$ is the number of cylinders (if cylinders are allowed, that is, in the case of loops).
	
	With this notation, we recall here now the two main formulas, the first from \cite[Section 8.5]{eskin_masur_zorich_03}.
	
	\begin{thm}[Formula for Siegel--Veech constant for distinct zeros \cite{eskin_masur_zorich_03}]
		\label{thm:emz_distinct_labelled_zeros}
		Let $\mathcal{H}$ be a connected stratum of unit-area translation surfaces.
		We consider the configuration of a saddle connection of multiplicity $p$ between a fixed zero of order $m_1$ and a different, fixed zero of order $m_2$.
		
		For fixed choices of $\mathcal{H}_i$ (including labelling the zeros) and of $a_i = a_i'+a_i''$, the Siegel--Veech constant is
		\begin{align*}
			& c(m_1 \saddleconnection m_2, p, \labelled, \mathcal{H}) \\
			& = \frac{1}{|\Gamma| \cdot |\Gamma_-|} \cdot \prod_{i=1}^p (a_i +1 ) \cdot \frac{1}{2^{p-1}} 
			\cdot \frac{\prod_{i=1}^p (\sfrac{d_i}{2}-1)!}{(\sfrac{d}{2}-2)!} 
			\cdot \frac{\prod_{i=1}^p \mu(\mathcal{H}_i)}{\mu(\mathcal{H})} 
			.
		\end{align*}
	\end{thm}
	
	For multiplicity $p=1$, we necessarily have $a_1' = m_1$ and $a_1'' = m_2$ and hence $\mathcal{H}_1$ is the (not necessarily connected) stratum of unit-area translation surfaces which have one less zero of order~$m_1$, one less zero of order $m_2$, and one more zero of order $a_1 = m_1+m_2$ than translation surfaces in $\mathcal{H}$.
	Furthermore, $\sfrac{d_1}{2} -1 = 2g + \ell -3 = \sfrac{d}{2} -2$ and hence, the Siegel--Veech constant simplifies to
	\begin{equation} \label{eq:SV_distinct_multiplicity_one}
		c(m_1 \saddleconnection m_2, p=1, \fixedzeros, \mathcal{H}) 
		= (m_1 + m_2 +1 ) \cdot \frac{\mu(\mathcal{H}_1)}{\mu(\mathcal{H})}
	\end{equation}
	for the sole eligible labelled configuration as well as for the unlabelled configuration \cite[Section 8.3]{eskin_masur_zorich_03}.
	
	The second main formula, from \cite[Section 13.3]{eskin_masur_zorich_03}, is for loops.
	
	\begin{thm}[Formula for Siegel--Veech constant for loops \cite{eskin_masur_zorich_03}]
		\label{thm:emz_loops_labelled_zeros}
		Let $\mathcal{H}$ be a connected stratum of unit-area translation surfaces.
		We consider the configuration of a loop of multiplicity $p$ at a fixed zero of order $m$.
		
		For fixed choices of $\mathcal{H}_i$ (including labelling the zeros), position and number of figure-eight constructions, two-hole constructions, and cylinders, and of $a_i = a_i'+a_i''$ or $(b_k', b_k'')$, the Siegel--Veech constant is
		\begin{align*}
			& c(m \scloop, p, \labelled, \mathcal{H}) \\
			& = \frac{1}{|\Gamma| \cdot |\Gamma_-|} \cdot \prod (a_i +1 ) \cdot \prod (b_k' +1)(b_k'' +1) \cdot \frac{1}{2^{p-1}} 
			\cdot \frac{\prod_{i=1}^p (\sfrac{d_i}{2}-1)!}{(\sfrac{d}{2}-2)!} 
			\cdot \frac{\prod_{i=1}^p \mu(\mathcal{H}_i)}{\mu(\mathcal{H})}
			.
		\end{align*}
	\end{thm}

	We want to point out here that \cref{eq:volume} will be applied for $\mathcal{H}$ in the following calculations but is not necessarily useful for the $\mathcal{H}_i$ as $g \to \infty$ does not imply $g_i \to \infty$ for every $i = 1, \ldots, p$. For example, we have (see \cite{zorich_02})
	\begin{equation}
		\label{eq:volume_H0}
		\mu(\mathcal{H}(0)) = \frac{\pi^2}{3}
		.
	\end{equation}
	
	However, we can use \cref{eq:volume} also in the cases where $g_i \nrightarrow \infty$ but will do so explicitly with a worse error term, that is
	\begin{equation}
		\label{eq:volume_low-genus}
		\mu(\mathcal{H}( m_1, \ldots, m_\ell) ) 
		= \frac{4}{\prod_{n=1}^\ell (m_n + 1)} \cdot O(1)
		.
	\end{equation}

	\section{Siegel--Veech constants for non-connected strata}
	\label{sec:non-connected_strata}
	
	Most of the known large-genus asymptotics for Siegel--Veech constants have been derived for strata which are connected. The reason for the distinction between connected and non-connected strata is that for non-connected strata, we have to determine whether a surface which is obtained from gluing subsurfaces with boundary has odd or even spin structure and/or whether it is hyperelliptic. Hence, more case analysis is needed~in~this~situation.
	
	Eskin, Masur, and Zorich also considered these situations and gave explicit formulas which involve volumes of connected components of the $\mathcal{H}_i$ instead of volumes of the whole strata.
	We use their arguments, together with the volume estimates for the connected components from \cref{eq:volume_odd_even,eq:volume_hyperelliptic_minimal_Stirling,eq:volume_hyperelliptic_two_Stirling}, to show that an approximate version of the formulas from~\cref{thm:emz_distinct_labelled_zeros,thm:emz_loops_labelled_zeros} also hold for non-connected strata and for the non-hyperelliptic components of such~strata.
	
	\bigskip
	
	We start with the exceptional components which are the \emph{hyperelliptic components} of non-connected strata. These only appear for strata of the form $\mathcal{H}(2g-2)$ and~$\mathcal{H}(g-1, g-1)$ for $g \geq 3$ as these strata contain some translation surfaces that admit a hyperelliptic involution. For $g \in \{1,2\}$, all translation surfaces are hyperelliptic, hence the strata~$\mathcal{H}(0)$, $\mathcal{H}(1,1)$, and $\mathcal{H}(2)$ are connected and coincide with their hyperelliptic component.
	
	We first consider saddle connections between two different zeros. By definition, these cannot exist on surfaces from $\mathcal{H}^\hyp(2g-2)$, only on surfaces from $\mathcal{H}(g-1,g-1)$.
	To obtain a hyperelliptic translation surface in $\mathcal{H}$ from gluing subsurfaces with boundary, all the translation surfaces from $\mathcal{H}_i$ also have to be hyperelliptic. This also shows that there is only one choice for the partition~$a_i = a_i' + a_i''$. This means that there are fewer labelled configurations than in the general case; concretely, for every unlabelled configuration, there exists exactly one labelled configuration.
	
	The following theorem collects the statements from \cite[Corollary 10.4 and Formulas 10.1 and 10.2]{eskin_masur_zorich_03} on Siegel--Veech constants for saddle connections between different zeros in $\mathcal{H}^\hyp(g-1,g-1)$.

	\begin{thm}[{\textls[-9]{Saddle connections between distinct zeros in $\mathcal{H}^\hyp(\gminusone,\gminusone)$ \cite{eskin_masur_zorich_03}}}]
		\label{thm:emz_hyperelliptic_distinct}
		Let $\mathcal{H} = \mathcal{H}^\hyp (g-1,g-1)$ for some genus $g \geq 3$.
		For almost every translation surface in $\mathcal{H}$, there are no saddle connections of multiplicity greater than $2$.
		
		The Siegel--Veech constant for the configuration of a saddle connection between two different zeros
		of multiplicity $1$ is
		\begin{equation*}
			c(\gminusone \saddleconnection \gminusone, p=1, \mathcal{H}^\hyp(g-1,g-1))
			= 
			(2g-1) \cdot \frac{\mu(\mathcal{H}^\hyp(2g-2))}{\mu(\mathcal{H}^\hyp(g-1,g-1))} 
		\end{equation*}
		and of multiplicity $2$ with a fixed partition $g= g_1+g_2$ is
		\begin{align*}
			& c(\gminusone \saddleconnection \gminusone, p=2, \labelled, \mathcal{H}^\hyp(g-1,g-1)) \\
			& = \frac{(2g_1-1)(2g_2-1)}{2 \cdot |\Gamma|} 
			\cdot \frac{(2g_1-1)! \cdot (2g_2-1)!}{(2g-1)!}
			\cdot \frac{\mu(\mathcal{H}^\hyp(2g_1-2)) \cdot \mu(\mathcal{H}^\hyp(2g_2-2)) }{\mu(\mathcal{H}^\hyp(g-1,g-1))} 
		\end{align*}
		with $|\Gamma| = 2$ if and only if $g_1 = g_2$.
	\end{thm}
	
	Similar arguments can also be used for loops, now in $\mathcal{H}^\hyp(g-1, g-1)$ as well as in~$\mathcal{H}^\hyp(2g-2)$.
	The following theorem recalls the statements from \cite[Lemma~14.6 and Formula 14.2]{eskin_masur_zorich_03} for $\mathcal{H}^\hyp(g-1,g-1)$.
	
	\begin{thm}[Loops in $\mathcal{H}^\hyp(g-1,g-1)$ \cite{eskin_masur_zorich_03}]
		\label{thm:emz_hyperelliptic_gminusone_gminusone_loops}
		Let $\mathcal{H} = \mathcal{H}^\hyp(g-1,g-1)$ for some genus $g \geq 3$.
		For almost every translation surface in $\mathcal{H}$, there are no loops of multiplicity greater than $2$.
		
		The Siegel--Veech constant for the configuration of a loop at any zero of multiplicity~$1$~is
		\begin{equation*}
			c( \gminusone \scloop, p=1, \anyzero, \mathcal{H}^\hyp(g-1,g-1))
			= 
			\frac{g-1}{2g-1} \cdot \frac{\mu(\mathcal{H}^\hyp(g-2,g-2))}{\mu(\mathcal{H}^\hyp(g-1,g-1))} 
		\end{equation*}
		and of multiplicity $2$ with a fixed partition $g= g_1+g_2$ is
		\begin{align*}
			& c(\twogminustwo \scloop, p=2, \anyzero, \mathcal{H}^\hyp(g-1,g-1)) \\
			& = \frac{g_1 \cdot g_2}{2 \cdot |\Gamma|} 
			\cdot \frac{(2g_1)! \cdot (2g_2)!}{(2g-1)!}
			\cdot \frac{\mu(\mathcal{H}^\hyp(g_1-1, g_1-1)) \cdot \mu(\mathcal{H}^\hyp(g_2-1, g_2-1)) }{\mu(\mathcal{H}^\hyp(g-1,g-1))} 
		\end{align*}
		with $|\Gamma| = 2$ if and only if $g_1 = g_2$.
	\end{thm}
	
	We also recall the statements from \cite[Lemma 14.5 and Formula 14.1]{eskin_masur_zorich_03} for $\mathcal{H}(2g-2)$.
	
	\begin{thm}[Loops in $\mathcal{H}^\hyp(2g-2)$ \cite{eskin_masur_zorich_03}]
		\label{thm:emz_hyperelliptic_twogminustwo_loops}
		Let $\mathcal{H} = \mathcal{H}^\hyp(2g-2)$ for some genus $g \geq 3$.
		For almost every translation surface in $\mathcal{H}$, there are no loops of multiplicity greater than $2$.
		
		The Siegel--Veech constant for the configuration of a loop of multiplicity $1$ and $\mathcal{H}_1 = \mathcal{H}^\hyp(g-1,g-1)$ is
		\begin{equation*}
			c(\twogminustwo \scloop, p=1, \mathcal{H}^\hyp(2g-2))
			= \frac{g-1}{2} \cdot \frac{\mu(\mathcal{H}^\hyp(g-2, g-2))}{\mu(\mathcal{H}^\hyp(2g-2))}
			,
		\end{equation*}
		of multiplicity $1$ and $\mathcal{H}_1 = \mathcal{H}^\hyp(2g-4)$ is
		\begin{equation*}
			c(\twogminustwo \scloop, p=1, \mathcal{H}^\hyp(2g-2))
			= \frac12 \cdot \frac{2g-3}{2g-2} \cdot \frac{\mu(\mathcal{H}^\hyp(2g-4))}{\mu(\mathcal{H}^\hyp(2g-2))}
			,
		\end{equation*}
		and of multiplicity $2$ with a fixed partition $g-1= g_1+g_2$ is
		\begin{align*}
			& c(\twogminustwo \scloop , p=2, \mathcal{H}^\hyp(2g-2)) \\
			& = \frac{(2g_1-1) \cdot 2g_2}{8} 
			\cdot \frac{(2g_1-1)! \cdot (2g_2)!}{(2g-2)!}
			\cdot \frac{\mu(\mathcal{H}^\hyp(2g_1-2)) \cdot \mu(\mathcal{H}^\hyp(g_2-1, g_2-1)) }{\mu(\mathcal{H}^\hyp(2g-2))} 
			.
		\end{align*}
	\end{thm}
	
	We turn now to the other connected components. If translation surfaces in a given stratum have only zeros of even order, then these surfaces have a spin structure whose parity determines whether they belong to the \emph{odd} or \emph{even component} (as long as they do not have a hyperelliptic involution and belong to the hyperelliptic component).
	In particular, strata of the type $\mathcal{H}(2g-2)$ and $\mathcal{H}(g-1,g-1)$ for odd $g$ have three connected components if $g\geq 3$.  Strata of the type  $\mathcal{H}(g-1,g-1)$ for even $g \geq 4$ have two connected components and the complement of the hyperelliptic component is called \emph{non-hyperelliptic component}.
	As detailed in \cite[Sections 10 and 14]{eskin_masur_zorich_03}, for the odd/even/non-hyperelliptic components, we have the same combinatorics as for connected strata but we have to make sure that we do not produce hyperelliptic translation surfaces and to consider whether the surfaces from~$\mathcal{H}_i$ have a spin structure and which connected components of $\mathcal{H}_i$ are eligible to obtain translation surfaces from $\mathcal{H}^\oddeven$.
	
	Therefore, we have to replace in the formulas from \cref{thm:emz_distinct_labelled_zeros,thm:emz_loops_labelled_zeros} the term $\mu(\mathcal{H})$ by $\mu(\mathcal{H}^\oen)$ and also replace the terms $\mu(\mathcal{H}_i)$ accordingly.
	For this, recall that \cref{eq:volume_hyperelliptic_minimal_Stirling,eq:volume_hyperelliptic_two_Stirling} imply
	\begin{equation}
		\label{eq:volume_hyperelliptic_approximate}
		\mu(\mathcal{H}^\hyp) = \mu(\mathcal{H}) \cdot O \left( \frac1g \right)
	\end{equation}
	and hence with \cref{eq:volume_odd_even}, we have
	\begin{align}
		\mu(\mathcal{H}^\nonhyp) & = \mu(\mathcal{H}) \cdot \left( 1 + O \left( \frac1g \right) \right), \label{eq:volume_nonhyp} \\
		\mu(\mathcal{H}^\odd) & = \frac12 \cdot \mu(\mathcal{H}) \cdot \left( 1 + O \left( \frac1g \right) \right), \label{eq:volume_odd} \\
		\mu(\mathcal{H}^\even) & = \frac12 \cdot \mu(\mathcal{H}) \cdot \left( 1 + O \left( \frac1g \right) \right), \label{eq:volume_even}	
	\end{align}
	whenever these connected components exist.
	Note that for the $\mathcal{H}_i$, the same formulas with error term $O(1)$ are more useful (cf.~the remarks before~\cref{eq:volume_low-genus}).
	
	When determining Siegel--Veech constants in \cref{sec:known_values_non-loops,sec:known_values_loops}, we will consider separately dominant and non-dominant configurations. A configuration is called \emph{dominant} if all but one $\mathcal{H}_i$ are equal to $\mathcal{H}(0)$ or $\mathcal{H}(0,0)$. In the non-dominant cases, an error term of~$O(1)$ is sufficient, hence we will work with the following two propositions.
	
	\begin{prop}[{{\textls[-28]{\spaceskip 0.17em Saddle connections between distinct zeros in non-connected strata}}}]
		\label{prop:distinct_non-connected_strata}
		Let $\mathcal{H}$ be a (not necessarily connected) stratum
		and $\mathcal{H}^\oen$ the odd/even/non-hyperelliptic components of $\mathcal{H}$ whenever they exist.
		We consider the configuration of a saddle connection of multiplicity $p$ between a fixed zero of order $m_1$ and a different, fixed zero of order $m_2$.
		
		For fixed choices of $\mathcal{H}_i$(including labelling the zeros) and of $a_i = a_i'+a_i''$, the Siegel--Veech constant for a dominant configuration (where all $\mathcal{H}_i$ but $\mathcal{H}_1$ are equal to $\mathcal{H}(0)$)~is
		\begin{align*}
			& c(m_1 \saddleconnection m_2, p, \labelled, \dominant, \mathcal{H} / \mathcal{H}^\oen) \\
			= & \frac{(m_1+1)(m_2+1)}{|\Gamma| \cdot |\Gamma_-|} 
			\cdot \left( \frac{\pi^2}{6} \right)^{p-1}
			\cdot \frac{\prod_{i=1}^p (\sfrac{d_i}{2}-1)!}{(\sfrac{d}{2}-2)!} 
			\cdot \left(1 + O \left( \frac1g \right)\right)
		\end{align*}
		and the Siegel--Veech constant for a non-dominant configuration is
		\begin{align*}
			& c(m_1 \saddleconnection m_2, p, \labelled, \nondominant, \mathcal{H} / \mathcal{H}^\oen) \\
			= & (m_1+1)(m_2+1) 
			\cdot \frac{\prod_{i=1}^p (\sfrac{d_i}{2}-1)!}{(\sfrac{d}{2}-2)!} 
			\cdot O \left( 1 \right)^p
			.
		\end{align*}
		
		\begin{proof}
			For a dominant configuration where all $\mathcal{H}_i$ but $\mathcal{H}_1$ are equal to $\mathcal{H}(0)$,
			we have that~$\mathcal{H}_1$ contains all zeros of $\mathcal{H}$ but the fixed ones of order $m_1$ and $m_2$ whereas it additionally contains a zero of order $a_1$. Hence, with \cref{eq:volume,eq:volume_H0}, we have
			\begin{align*}
				& \prod_{i=1}^p (a_i +1 ) \cdot \frac{1}{2^{p-1}}
				\cdot \frac{\prod_{i=1}^p \mu(\mathcal{H}_i)}{\mu(\mathcal{H})} \\
				& = (a_1+1) \cdot \frac{1}{2^{p-1}} \cdot \frac{4 \cdot (m_1+1)(m_2+1) \cdot \mu(\mathcal{H}(0))^{p-1}}{4 \cdot (a_1 +1)} 
				\cdot \left( 1+ O \left( \frac1g \right) \right)^2 \\
				& = \left( \frac{\pi^2}{6} \right)^{p-1} \cdot (m_1+1)(m_2+1)
				\cdot \left( 1 + O \left( \frac1g \right) \right)
				.
			\end{align*}
			
			For a non-dominant configuration, as every zero in $\mathcal{H}$ is either one of the two fixed zeros of order $m_1$ and $m_2$ or a zero in one of the $\mathcal{H}_i$, and as the $\mathcal{H}_i$ have as only additional zeros the ones of order $a_i$, we have with \cref{eq:volume_low-genus}
			\begin{align*}
				& \prod_{i=1}^p (a_i +1 ) \cdot \frac{1}{2^{p-1}}
				\cdot \frac{\prod_{i=1}^p \mu(\mathcal{H}_i)}{\mu(\mathcal{H})} \\
				& = \prod_{i=1}^p (a_i +1 ) \cdot \frac{1}{2^{p-1}}
				\cdot \frac{4^{p} \cdot (m_1 + 1) (m_2 +1)}{4 \cdot \prod_{i=1}^p (a_i +1)} \cdot O(1)^p \\
				& = (m_1+1)(m_2+1) \cdot O(1)^p
				.
			\end{align*}
			
			For connected strata, plugging these estimates into \cref{thm:emz_distinct_labelled_zeros} (and disregarding the factor $\frac{1}{|\Gamma| \cdot |\Gamma_-|}$ for non-dominant configurations) yields the statements.
			
			Let $\mathcal{H}$ be now a non-connected stratum.
			If $\mathcal{H}$ contains a zero of odd order, the surfaces from $\mathcal{H}$ do not have a spin structure, so we only need to consider whether surfaces are hyperelliptic or not.
			Generically, surfaces can only be hyperelliptic if all $\mathcal{H}_i$ are hyperelliptic.
			This means that for $\mathcal{H}^\nonhyp$, we have to replace the term $\frac{\prod_{i=1}^p \mu(\mathcal{H}_i)}{\mu(\mathcal{H})}$ in \cref{thm:emz_distinct_labelled_zeros} by
			\begin{equation*}
				\frac{\prod_{i=1}^p \mu(\mathcal{H}_i) - \prod_{i=1}^p \mu(\mathcal{H}^\hyp_i)}{\mu(\mathcal{H}^\nonhyp)}
				= \frac{\prod_{i=1}^p \mu(\mathcal{H}_i)}{\mu(\mathcal{H})}
				\cdot \left( 1 + O\left( \frac1g \right) \right)
			\end{equation*}
			where we can use \cref{eq:volume_nonhyp} and the stronger error bounds from \cref{eq:volume_hyperelliptic_approximate} as~$g\to \infty$ implies~$g_i \to \infty$ for at least one $g_i$.
			With this, the statements in the case~$\mathcal{H}^\nonhyp$ follow with the same estimates as for connected strata.
			
			If $\mathcal{H}$ contains only zeros of even order, then the surfaces from $\mathcal{H}$ have a spin structure.
			Furthermore, the surfaces from all of the $\mathcal{H}_i$ have also only zeros of even order and hence a spin structure.
			The parity of the spin structure of the surfaces obtained by gluing is determined by the sum of the parities of the surfaces from $\mathcal{H}_i$.
			Note that hyperelliptic surfaces are sometimes odd and sometimes even which means that the hyperelliptic components are sometimes eligible. The details can be found in \cite[Formulas~10.6--10.8]{eskin_masur_zorich_03} but do not matter for our rough estimates.
			In summary, for $\mathcal{H}^\oddeven$, we have to replace the term $\frac{\prod_{i=1}^p \mu(\mathcal{H}_i)}{\mu(\mathcal{H})}$ in \cref{thm:emz_distinct_labelled_zeros} by
			\begin{equation*}
				\sum_{\substack{\phi_1, \ldots, \phi_p \in \{\odd, \even\} \\ \phi_1 + \ldots + \phi_p = \ast}}
				\frac{	\prod_{i=1}^p \mu( \mathcal{H}_i^{\phi_i} [\cup \mathcal{H}_i^\hyp] )}{\mu(\mathcal{H}^\oddeven)}
				\left[- \frac{	\prod_{i=1}^p \mu( \mathcal{H}_i^\hyp )}{\mu(\mathcal{H}^\oddeven)} \right]
				.
			\end{equation*}
			
			For dominant configurations, all but one of the $\mathcal{H}_i$ are equal to $\mathcal{H}(0)$ and as $\mathcal{H}(0)$ is a connected stratum, we have only one choice for the connected components of such~$\mathcal{H}_i$. Furthermore, for $\mathcal{H}_1$, only surfaces from one of the two connected components $\mathcal{H}_1^\odd$ and~$\mathcal{H}_1^\even$ produce surfaces with the right parity -- which component this is, depends on the genus. Therefore, the volume term simplifies to
			\begin{equation*}
				\frac{ \mu( \mathcal{H}_1^\oddeven [\cup \mathcal{H}_1^\hyp] ) \cdot \mu(\mathcal{H}(0))^{p-1}}{\mu(\mathcal{H}^\oddeven)}
				= \left( \frac{\pi^2}{3} \right)^{p-1} \cdot \frac{ \sfrac12 \cdot \mu( \mathcal{H}_1)}{ \sfrac12 \cdot \mu(\mathcal{H})}
				\cdot \left( 1 + O \left( \frac1g \right) \right)
			\end{equation*}
			which, in combination with \cref{eq:volume}, proves the statement for dominant configurations.
			For non-dominant configurations, 
			we actually have up to $2^{p-1}$ eligible combinations of connected components and can use the weaker error terms from \cref{eq:volume_hyperelliptic_approximate,eq:volume_nonhyp,eq:volume_odd,eq:volume_even} to obtain for the volume term
			\begin{equation*}
				2^{p-1} \cdot 
				\frac{	\prod_{i=1}^p \sfrac12 \cdot \mu( \mathcal{H}_i )}{\sfrac12 \cdot \mu(\mathcal{H})}
				\cdot O(1)^p
				= \frac{ \prod_{i=1}^p \mu( \mathcal{H}_i)}{\mu(\mathcal{H})} \cdot O(1)^p
			\end{equation*}
			from which the statement for non-dominant configurations follows as above with $O(1)^p \cdot O(1)^p = O(1)^p$.
			
			We consider now the whole non-connected stratum $\mathcal{H}$. The Siegel--Veech constant can be determined by averaging the Siegel--Veech constants for the connected components as	
			\begin{equation*}
				c( \param, \mathcal{H} )
				= c( \param, \mathcal{H}^\hyp ) \cdot \frac{\mu (\mathcal{H}^\hyp ) }{ \mu( \mathcal{H})}
				+ c( \param, \mathcal{H}^\nonhyp ) \cdot \frac{\mu (\mathcal{H}^\even ) }{ \mu( \mathcal{H})}
			\end{equation*}
			or
			\begin{equation*}
				c( \param, \mathcal{H} )
				= c( \param, \mathcal{H}^\hyp ) \cdot \frac{\mu (\mathcal{H}^\hyp ) }{ \mu( \mathcal{H})}
				+ c( \param, \mathcal{H}^\odd ) \cdot \frac{\mu (\mathcal{H}^\odd ) }{ \mu( \mathcal{H})}
				+ c( \param, \mathcal{H}^\even ) \cdot \frac{\mu (\mathcal{H}^\even ) }{ \mu( \mathcal{H})}
				.
			\end{equation*}
			
			We can read from \cref{thm:emz_hyperelliptic_distinct} and will also see in detail in \cref{prop:distinct_hyperelliptic_g-1} (which is independent of the \cref{sec:known_values_non-loops,sec:known_values_loops}) that the Siegel--Veech constants for the hyperelliptic components are at most order $g^2$ larger than the Siegel--Veech constants for the odd/even/non-hyperelliptic components. As the volumes of the hyperelliptic components are exponentially smaller than the volume of the odd/even/non-hyperelliptic components, we can disregard the term for the hyperelliptic components.
			Hence, we have
			\begin{equation*}
				c( \param, \mathcal{H} )
				= c( \param, \mathcal{H}^\oen )
				\cdot \left( 1 + O \left( \frac1g \right) \right)
			\end{equation*}
			which shows the statements for non-connected strata $\mathcal{H}$.
		\end{proof}
	\end{prop}
	
	We also have analogous formulas for loops.
	
	\begin{prop}[Loops in non-connected strata]
		\label{prop:loops_non-connected_strata}
		Let $\mathcal{H}$ be a (not necessarily connected) stratum
		and $\mathcal{H}^\oen$ the odd/even/non-hyperelliptic components of $\mathcal{H}$ whenever they exist.
		We consider the configuration of a loop of multiplicity $p$ at a fixed zero of order $m$.
		
		For fixed choices of $\mathcal{H}_i$ (including labelling the zeros), position and number of figure-eight constructions, two-hole constructions, and cylinders, and of $a_i = a_i'+a_i''$ or $(b_k', b_k'')$, 
		let $m, m_2, \ldots, m_n$ be the order of the zeros that are obtained by gluing.
		Then the Siegel--Veech constant for a dominant configuration (where all $\mathcal{H}_i$ but $\mathcal{H}_1$ are equal to~$\mathcal{H}(0)$ or~$\mathcal{H}(0,0)$) is
		\begin{align*}
			& c(m \scloop , p, \labelled, \dominant, \mathcal{H} / \mathcal{H}^\oen) \\
			= & \frac{(m+1) \cdots (m_n+1)}{|\Gamma| \cdot |\Gamma_-|}
			\cdot \left( \frac{\pi^2}{6} \right)^{p-1}
			\cdot \frac{\prod_{i=1}^p (\sfrac{d_i}{2} - 1)!}{(\sfrac{d}{2} - 2)!}
			\cdot \left( 1 + O \left( \frac1g \right) \right)
		\end{align*}
		and the Siegel--Veech constant for a non-dominant configuration is
		\begin{align*}
			& c(m \scloop , p, \labelled, \nondominant, \mathcal{H} / \mathcal{H}^\oen) \\
			= & (m+1) \cdots (m_n+1)
			\cdot \frac{\prod_{i=1}^p (\sfrac{d_i}{2}-1)!}{(\sfrac{d}{2}-2)!} 
			\cdot O \left( 1 \right)^p
			.
		\end{align*}
		
		\begin{proof}
			The proof follows the same strategy as the proof of \cref{prop:distinct_non-connected_strata}.
			In particular, as $\mu(\mathcal{H}(0)) = \mu(\mathcal{H}(0,0))$, we have again
			\begin{align*}
				& \prod (a_i +1 ) \prod (b_k'+1)(b_k''+1)
				\cdot \frac{1}{2^{p-1}} \cdot \frac{\prod_{i=1}^p \mu(\mathcal{H}_i)}{\mu(\mathcal{H})} \\
				& = \left( \frac{\pi^2}{6} \right)^{p-1}
				\cdot (m+1) \cdots (m_n+1)
				\cdot \left( 1+ O \left( \frac1g \right) \right)
			\end{align*}
			for a dominant configuration and
			\begin{equation*}
				\prod (a_i +1 ) \cdot \prod (b_k'+1)(b_k''+1) \cdot \frac{1}{2^{p-1}}
				\cdot \frac{\prod_{i=1}^p \mu(\mathcal{H}_i)}{\mu(\mathcal{H})}
				= (m_1+1) \cdots (m_{n}+1) \cdot O(1)^p
			\end{equation*}
			for a non-dominant configuration, as every zero in $\mathcal{H}$ is either one of the fixed zeros of order $m_1, \ldots, m_n$ or a zero in one of the $\mathcal{H}_i$, and as the $\mathcal{H}_i$ have as only additional zeros the ones of order $a_i$ or $(b_k', b_k'')$.
			
			For connected strata, plugging these estimates into \cref{thm:emz_loops_labelled_zeros} (and disregarding the factor $\frac{1}{|\Gamma| \cdot |\Gamma_-|}$ for non-dominant configurations) yields the statements.
			
			Let $\mathcal{H}$ now be a non-connected stratum.
			If $\mathcal{H}$ contains a zero of odd order, we need to exclude only hyperelliptic surfaces. As we cannot construct all hyperelliptic surfaces this way, we have to remove even only a fraction of the volume. In particular, for $\mathcal{H}^\nonhyp$, we have to replace the term~$\frac{\prod_{i=1}^p \mu(\mathcal{H}_i)}{\mu(\mathcal{H})}$ in \cref{thm:emz_loops_labelled_zeros} by
			\begin{equation*}
				\frac{\prod_{i=1}^p \mu(\mathcal{H}_i) - \prod_{i=1}^p \mu(\mathcal{H}^\hyp_i) \cdot O(1)}{\mu(\mathcal{H}^\nonhyp)}
				= \frac{\prod_{i=1}^p \mu(\mathcal{H}_i)}{\mu(\mathcal{H})}
				\cdot \left( 1 + O\left( \frac1g \right) \right)
			\end{equation*}
			which proves the statements in the case $\mathcal{H}^\nonhyp$.
			
			If $\mathcal{H}$ contains only zeros of even order, the surfaces from $\mathcal{H}$ all have a spin structure. 
			However, for loops, there can exist eligible $\mathcal{H}_k$ without a spin structure. This happens exactly in the case with two-hole constructions where $b_k'$ and $b_k''$ are odd. Then surfaces from $\mathcal{H}_k$ do not have a spin structure but still half of them give a surface in $\mathcal{H}$ with the correct parity of the spin structure (see \cite[Lemma 14.4]{eskin_masur_zorich_03}).
			Hence for $\mathcal{H}^\oddeven$, we have to replace the term~$\frac{\prod_{i=1}^p \mu(\mathcal{H}_i)}{\mu(\mathcal{H})}$ in \cref{thm:emz_loops_labelled_zeros} by
			\begin{equation*}
				\frac12 \cdot \frac{\prod_{i=1}^p \mu(\mathcal{H}_i))}{\mu(\mathcal{H}^\oddeven)}
				= \frac{\prod_{i=1}^p \mu(\mathcal{H}_i)}{\mu(\mathcal{H})}
				\cdot \left( 1 + O\left( \frac1g \right) \right)
				.
			\end{equation*}
			If the surfaces from all $\mathcal{H}_i$ have a spin structure, we have to replace the term~$\frac{\prod_{i=1}^p \mu(\mathcal{H}_i)}{\mu(\mathcal{H})}$ in \cref{thm:emz_loops_labelled_zeros} by
			\begin{equation*}
				\sum_{\substack{\phi_1, \ldots, \phi_p \in \{\odd, \even\} \\ \phi_1 + \ldots + \phi_p = \ast}}
				\frac{	\prod_{i=1}^p \mu( \mathcal{H}_i^{\phi_i} [\cup \mathcal{H}_i^\hyp] )}{\mu(\mathcal{H}^\oddeven)}
				\left[- \frac{	\prod_{i=1}^p \mu( \mathcal{H}_i^\hyp )}{\mu(\mathcal{H}^\oddeven)} \cdot O(1) \right]
				.
			\end{equation*}
			From here, we can conclude literally as in the proof of \cref{prop:distinct_non-connected_strata} that the statements also hold true for $\mathcal{H}^\oddeven$.
			
			Also the statement for the whole non-connected stratum $\mathcal{H}$ follows as in the proof of \cref{prop:distinct_non-connected_strata}.
		\end{proof}
	\end{prop}

	\section{Some combinatorial lemmas} \label{sec:lemmas}
	
	As the formulas for Siegel--Veech constants are all built in the same way, we collect some arguments in this section that will be used repeatedly.
	
	For instance, in any Siegel--Veech constant, we have to consider the dimension term~$\frac{\prod_{i=1}^p (\sfrac{d_i}{2} -1)!}{(\sfrac{d}{2}-2)!}$.
	Keep in mind that $\sum_{i=1}^p (\sfrac{d_i}{2} -1)$ is equal to $\sfrac{d}{2}-p-1$ (for saddle \linebreak connections between distinct zeros) or to $\sfrac{d}{2} - p -q -1$ (for loops), hence the number of factors in~$\prod_{i=1}^p (\sfrac{d_i}{2} -1)!$ is always the same for given $d$, $p$, and potentially $q$. Therefore, the dimension term becomes largest if one of the $\sfrac{d_i}{2}$ is as large as possible.
	
	Such a maximal dimension term corresponds to dominant configurations as we will see later. For non-dominant configurations, we use the next lemma (inspired by \cite[Lemma 2.3]{aggarwal_19}) to replace the factors $(\sfrac{d_i}{2} -1)!$ by simpler terms.
	
	\begin{lem}
		\label{lem:comparison_product_factorials_of_two_partitions}
		Let $p \geq 1$, $r \geq 0$ and $A,B \geq 0$. Let $A = A_1 + \ldots + A_p$ and $B = B_1 + \ldots + B_p$ be partitions such that $0 \leq B_i \leq A_i$ for every $i = 1, \ldots, p$.
		Furthermore, we assume~$B_i < B$ for every $i = 1, \ldots, p$.
		Then we have
		\begin{equation*}
			\frac{\prod_{i=1}^p (A_i+r)!}{(A+r-1)!}
			\leq \frac{\prod_{i=1}^p (B_i+r)!}{(B+r-1)!}
			.
		\end{equation*}
		
		\begin{proof}
			For fixed $i$, we have
			\begin{equation*}
				\frac{ (A_i +r)!}{ (B_i +r)!}
				= (B_i + r +1) \cdot \ldots \cdot (A_i + r)
				= \prod_{j=1}^{A_i-B_i} (B_i + r + j)
				.
			\end{equation*}
			
			As $B_i + r + j \leq B -1 + r + k$ for any $j \leq k$ and $\sum_{i=1}^p A_i - B_i = A -B$, it holds
			\begin{equation*}
				\prod_{i=1}^p \frac{ (A_i +r)!}{ (B_i +r)!}
				\leq \prod_{i=1}^p \prod_{j=1}^{A_i-B_i} (B -1 + r + j)
				\leq \prod_{k=1}^{A-B} (B + r -1 +k)
				= \frac{(A+r-1)!}{(B+r-1)!}
				.
				\qedhere
			\end{equation*}
		\end{proof}
	\end{lem}

	The next lemma will be used to simplify the dimension term even more.
	
	\begin{lem}
		\label{lem:product_binomials}
		Let $p \geq 1$ and $r\geq 0$. Furthermore, let $\ell_i, a_i \geq 0$ for every $i=1, \ldots, p$.
		Then
		\begin{equation*}
			\frac{\prod_{i=1}^p (2\ell_i + a_i + r)! }{ \prod_{i=1}^p \ell_i! }
			\leq 2^{rp} \cdot \binom{\sum_{i=1}^p (2\ell_i + a_i)}{\sum_{i=1}^p \ell_i}
			\cdot \prod_{i=1}^p (\ell_i + a_i + r)!
		\end{equation*}
		
		\begin{proof}
			Note first that for $A_1, A_2, \ell_1, \ell_2 \geq 0$, the term $\binom{A_1}{\ell_1} \cdot \binom{A_2}{\ell_2}$ denotes the number of choices of $\ell_1$ elements from a set of $A_1$ elements and simultaneously of $\ell_2$ elements from a set of $A_2$ elements. Each such choice is also a choice of $\ell_1 + \ell_2$ elements from a set of~$A_1 + A_2$ elements. Hence
			$\binom{A_1}{\ell_1} \cdot \binom{A_2}{\ell_2} \leq \binom{A_1 + A_2}{ \ell_1 + \ell_2}$.
			
			Using additionally $\binom{A_1+1}{\ell_1} = \frac{A_1+1}{A_1+1 - \ell_1} \cdot \binom{A_1}{\ell_1}$, we obtain
			\begin{align*}
				\prod_{i=1}^p \binom{2\ell_i + a_i + r}{\ell_i}
				& = \prod_{i=1}^p \frac{2\ell_i +a_i +1}{\ell_i +a_i +1} \cdot \ldots \cdot \frac{2\ell_i +a_i + r}{\ell_i +a_i + r}  \cdot \binom{2\ell_i + a_i}{\ell_i} \\
				& \leq 2^{rp} \cdot \prod_{i=1}^p \binom{2\ell_i + a_i}{\ell_i} \\
				& \leq 2^{rp} \cdot \binom{\sum_{i=1}^p (2\ell_i + a_i)}{\sum_{i=1}^p \ell_i}
				.
			\end{align*}
			Multiplying both sides by $\prod_{i=1}^p (\ell_i + a_i + r)!$ yields the statement.
		\end{proof}
	\end{lem}
	
	The following statement will be used in the proof of \cref{lem:sum_product_factorials_over_all_partitions}.
	
	\begin{lem}
		\label{lem:ingredient_for_sum_over_all_partitions}
		Let $p \geq 2$, $r\geq 0$, and $A \geq r$.
		Let $A = A_1 + \ldots + A_p$ be a partition of $A$ such that $A_i \leq A -3$ for every $i = 1, \ldots, p$.
		Then
		\begin{equation*}
			\prod_{i=1}^p (A_i + 3 + r)!
			\leq (3 + r)!^{p-1} \cdot (A + r)!
			.
		\end{equation*}
		
		\begin{proof}
			We prove the statement by induction on $p$. For $p = 2$, we have
			\begin{equation*}
				\prod_{i=1}^2 (A_i + 3 +r)!
				= (A_1 + 3 + r)! \cdot (A_2 + 3+ r)!
				\leq (A -3 +3 + r)! \cdot (3 + r)!
				= (3 + r)! \cdot (A+r)!
			\end{equation*}
			by choosing the largest possibility as upper bound.
			
			Now assume that the statement is true for some $p-1 \geq 2$ for all choices of $A$. Then
			\begin{equation*}
				\prod_{i=1}^p (A_i +3 + r)!
				= (A_p +3 + r)! \cdot \prod_{i=1}^{p-1} (A_i +3 + r)!
				\leq (3 + r)!^{p-2} \cdot (A_p +3 + r)! \cdot \left( \sum_{i=1}^{p-1} A_i + r \right)!
			\end{equation*}
			which is taking the largest value when $\sum_{i=1}^{p-1} A_i$ is the largest which is being equal to~$A$.~Hence
			\begin{equation*}
				\prod_{i=1}^p (A_i + 3 + r)!
				\leq (3 + r)!^{p-2} \cdot (0 + 3 +r)! \cdot (A + r)!
				= (3 + r)!^{p-1} \cdot (A+r)!
				. \qedhere
			\end{equation*}
		\end{proof}
	\end{lem}

	After considering dimension terms for specific configurations, we have to sum them up over all configurations. Aggarwal does this in \cite{aggarwal_19} by using the following crucial lemma.
	
	\begin{lem}[{\cite[Lemma 2.7]{aggarwal_19}}]
		There exists a constant $C>0$ such that for any~$p \geq 1$ and $M_1 , M_2, L \geq 0$, the sum over all ordered partitions of $M_1 = a_1' + \ldots + a_p'$, $M_2 = a_1'' + \ldots + a_p''$, and $L = \ell_1 + \ldots + \ell_p$ of the terms $\prod_{i=1}^p (a_i' + a_i'' + \ell_i +1)!$ is bounded as
		\begin{equation*}
			\sum_{M_1 = a_1' + \ldots + a_p'}
			\thickspace
			\sum_{M_2 = a_1'' + \ldots + a_p''}
			\thickspace
			\sum_{L = \ell_1 + \ldots + \ell_p}
			\thickspace
			\prod_{i=1}^p (a_i' + a_i'' + \ell_i +1)!
			\leq C^p \cdot (M_1 + M_2 + L +1)!
			.
		\end{equation*}
	\end{lem}
	
	We adapt Aggarwal's original proof of the lemma, so that we get the following versions which will give smaller upper bounds for non-dominant configurations.
	
	\begin{lem}
		\label{lem:sum_product_factorials_over_all_partitions}
		There exists a constant $C > 0$ such that the following holds.
		
		Let $p \geq 1$, $n\geq 1$ and $M_1, \ldots, M_n, L \geq 0$.
		We consider all choices of $a_1', \ldots, a_p'$, $a_1'', \ldots, a_p''$, $\ell_1, \ldots, \ell_p \geq 0$ such that $L = \ell_1 + \ldots + \ell_p$
		and each $a_i'$ and $a_i''$ is (once and for all) assigned to be a block in an ordered partition of $M_1, \ldots, M_n$ such that $a_i'$ and~$a_i''$ are assigned to different $M_1, \ldots, M_n$ for every $i=1,\ldots,p$.
		Furthermore, we assume $a_i' + a_i'' + \ell_i < M_1 + \ldots + M_n + L$ for every $i= 1, \ldots, p$.
		
		Then the sum of all terms $\prod_{i=1}^p (a_i' + a_i'' + \ell_i +1)!$ over all choices is bounded as
		\begin{equation*}
			\sum_{\substack{\text{partition} \\ \text{of } M_1}}
			\ldots
			\sum_{\substack{\text{partition} \\ \text{of } M_n}}
			\thickspace
			\sum_{L = \ell_1 + \ldots + \ell_p}
			\prod_{i=1}^p (a_i' + a_i'' + \ell_i +1)!
			\leq C^p \cdot (M_1 + \ldots + M_n + L)!
			.
		\end{equation*}
		
		Under the same conditions, we also have
		\begin{equation*}
			\sum_{\substack{\text{partition} \\ \text{of } M_1}}
			\ldots
			\sum_{\substack{\text{partition} \\ \text{of } M_n}}
			\thickspace
			\sum_{L = \ell_1 + \ldots + \ell_p}
			\prod_{i=1}^p (a_i' + a_i'' + \ell_i +2)!
			\leq C^p \cdot (M_1 + \ldots + M_n + L +1)!
			.
		\end{equation*}
		
		\begin{proof}
			Let $k' \colon \{1, \ldots, p\} \to \{1, \ldots, n\}$ be such that $a_i'$ is a block of the partition $M_{k'(i)}$ and $k'' \colon \{1, \ldots, p\} \to \{1, \ldots, n\}$ such that $a_i''$ is a block of the partition $M_{k''(i)}$. By assumption, $k'(i) \neq k''(i)$.
			We stress again that the assignments $k'$ and $k''$ are fixed.
			
			We start with a proof of the first bound.
			
			Let $D \coloneqq M_1 + \ldots + M_n+L$. For every partition of $D = d_1 + \ldots + d_p$, let $N(d_1, \ldots, d_p)$ be the number of ways to write this partition of $D$ as a tuple of partitions of $M_1, \ldots, M_n$, and~$L$, that is, as $d_i = a_i' + a_i'' + \ell_i$ for every $i = 1, \ldots, p$. Then we can write
			\begin{equation*}
				\sum_{\substack{\text{partition} \\ \text{of } M_1}}
				\ldots
				\sum_{\substack{\text{partition} \\ \text{of } M_n}}
				\thickspace
				\sum_{L = \ell_1 + \ldots + \ell_p}
				\prod_{i=1}^p (a_i' + a_i'' + \ell_i +1)!
				=
				\sum_{D = d_1 + \ldots + d_p}
				N(d_1, \ldots, d_p) \prod_{i=1}^p (d_i +1)!
				.
			\end{equation*}
			To determine the latter sum, we split it up by the number $s$ of $d_i$ which are not equal to~$0$.
			Because of the condition $d_i = a_i' + a_i'' + \ell_i < D$ for every $i=1,\ldots, p$, there are no partitions with $s=1$.
			
			Let now $s$ be fixed and every such partition $D= d_1 + \ldots + d_p$ written as $D = e_1 + \ldots + e_s$ where the $e_1, \ldots, e_s$ are equal to the $d_1, \ldots, d_p$ which are not $0$. We consider
			\begin{equation*}
				\sum_{D = e_1 + \ldots + e_s}
				N(e_1, \ldots, e_s) \prod_{i=1}^s (e_i +1)!
			\end{equation*}
			where $N(e_1, \ldots, e_s)$ is analogously the number of ways in which the partition $D = e_1 + \ldots + e_s$ can be written as a tuple of partitions of $M_1, \ldots, M_n$, and $L$ (into $s$ blocks).
			
			We distinguish between five cases, depending on the maximum of the $e_i$.
			
			Case 1: If one of the $e_i$ is equal to $D-(s-1)$, then all other $e_i$ are equal to $1$. Hence, there are $s$ such partitions	
			and for each of them, $N(e_1, \ldots, e_s) \leq s^2 \cdot 3^{s-1}$ where the first factor counts the ways to partition $e_i = D-(s-1) = M_1 + \ldots + M_n + L -(s-1)$ into $a_i' + a_i'' + \ell_i$ with $a_i' \geq M_{k'(i)}-(s-1)$, $a_i'' \geq M_{k''(i)} -(s-1)$, and $\ell_i = e_i - a_i' - a_i''$ and the second factor counts the ways to partition $1 = a_i' + a_i'' + \ell_i$ (which we have to do $s-1$ times).
			
			Case 2: If one of the $e_i$ is equal to $D-s$, then another $e_i$ is equal to $2$ and all other~$e_i$ are equal to $1$. Hence, there are $s (s-1)$ such partitions and for each of them, $N(e_1, \ldots, e_s) \leq (s+1)^2 \cdot 6 \cdot 3^{s-2}$ where the first factor counts the ways to partition $e_i = D-s$ into $a_i' + a_i'' + \ell_i$ with $a_i' \geq M_{k'(i)}-s$, $a_i'' \geq M_{k''(i)} -s$, and $\ell_i = e_i - a_i' - a_i''$, the second factor counts the ways to partition $2 = a_i' + a_i'' + \ell_i$, and the third factor counts the ways to partition $1 = a_i' + a_i'' + \ell_i$.
			
			Case 3: If one of the $e_i$ is equal to $D-s-1$ and one other $e_i$ is equal to $3$,
			then all other $e_i$ are equal to $1$.
			Hence, there are $s (s-1)$ such partitions and for each of them, $N(e_1, \ldots, e_s) \leq (s+2)^2 \cdot 10 \cdot 3^{s-2}$ where the first factor counts the ways to partition $e_i = D-s-1$ into $a_i' + a_i'' + \ell_i$ with $a_i' \geq M_{k'(i)}-s-1$, $a_i'' \geq M_{k''(i)} -s-1$, and $\ell_i = e_i - a_i' - a_i''$, the second factor counts the ways to partition $3 = a_i' + a_i'' + \ell_i$, and the third factor counts the ways to partition $1 = a_i' + a_i'' + \ell_i$.
			
			Case 4: If one of the $e_i$ is equal to $D-s-1$	and two other $e_i$ are equal to $2$,
			then all other $e_i$ are equal to $1$.
			Hence, there are $s (s-1) (s-2)$ such partitions and for each of them, 
			$N(e_1, \ldots, e_s) \leq (s+2)^2 \cdot 6^2 \cdot 3^{s-3}$ where the first factor counts the ways to partition $e_i = D-s-1$ into $a_i' + a_i'' + \ell_i$ with $a_i' \geq M_{k'(i)}-s-1$, $a_i'' \geq M_{k''(i)} -s-1$, and $\ell_i = e_i - a_i' - a_i''$, the second factor counts the ways to partition $2 = a_i' + a_i'' + \ell_i$, and the third factor counts the ways to partition $1 = a_i' + a_i'' + \ell_i$.
			
			Case 5: Now assume that all of the $e_i$ are at most $D-s-2$.
			For each $e_i$, we have at most $(e_i+1)^2$ partitions of $e_i = a_i' + a_i'' + \ell_i$. Hence
			\begin{equation*}
				N(e_1, \ldots, e_s) \prod_{i=1}^s (e_i +1)!
				\leq \prod_{i=1}^s (e_i +1)^2 \cdot \prod_{i=1}^s (e_i +1)!
				\leq \prod_{i=1}^s (e_i +3)!
				.
			\end{equation*}
			Now we apply \cref{lem:ingredient_for_sum_over_all_partitions} with $A_i = e_i -1$, $A = A_1 + \ldots + A_s = D - s$, and $r=1$. As $A_i \leq (D-s-2) -1 = A-3$ for every $i=1,\ldots, s$, this yields
			\begin{equation*}
				\prod_{i=1}^s (e_i +3)!
				\leq 4!^{s-1} \cdot (D-s+1)!
				.
			\end{equation*}
			The number of partitions $D = e_1 + \ldots + e_s$ where each $e_i$ is at most $D-s-2$ can be bounded by the number of choices of $e_1, \ldots, e_{s-1} \leq D-s-2$, that is, by $(D-s-2)^{s-1}$.
			
			Adding up the five cases, we have
			\begin{align*}
				\sum_{D = e_1 + \ldots + e_s}
				& N(e_1, \ldots, e_s) \prod_{i=1}^s (e_i +1)! \\
				& = s \cdot s^2 \cdot 3^{s-1} \cdot (D - s +2)! \cdot 2^{s-1} \\
				& + s(s-1) \cdot (s+1)^2 \cdot 6 \cdot 3^{s-2} \cdot (D - s +1)! \cdot 3! \cdot 2^{s-2} \\
				& + s(s-1) \cdot (s+2)^2 \cdot 10 \cdot 3^{s-2} \cdot (D-s)! \cdot 4! \cdot 2^{s-2} \\
				& + s(s-1)(s-2) \cdot (s+2)^2 \cdot 6^2 \cdot 3^{s-3} \cdot (D-s)! \cdot 3!^2 \cdot 2^{s-3} \\
				& + (D-s-2)^{s-1} \cdot 4!^{s-1} \cdot (D-s+1)! \\
				& \leq C_1^s \cdot D!
			\end{align*}
			for some constant $C_1$ which does not depend on $s$.
			
			We go back now to our original problem where $s$ and the choice of the indices whose blocks in the partition $D = d_1 + \ldots + d_p$ are not equal to $0$ is not fixed. As there are at most $2^p$ such choices, we have
			\begin{align*}
				\sum_{D = d_1 + \ldots + d_p}
				N(d_1, \ldots, d_p) \prod_{i=1}^p (d_i +1)!
				& \leq 2^p \cdot \max_{2 \leq s \leq p}
				\sum_{D = e_1 + \ldots + e_s}
				N(e_1, \ldots, e_s) \prod_{i=1}^s (e_i +1)! \\
				& \leq 2^p \cdot \max_{2 \leq s \leq p} C^s \cdot D! \\
				& = (2C_1)^p \cdot D!
				.
			\end{align*}
			Hence, the first bound is proven.
			
			For the second bound, we can follow the same arguments (with $d_i+2$ and $e_i+2$ instead of $d_i+1$ and $e_i+1$) up to and including Case 4.
			To replace Case 5, we again assume that all of the $e_i$ are at most $D-s-2$.
			For each $e_i$, we have at most $(e_i+1)^2$ partitions of $e_i = a_i' + a_i'' + \ell_i$. Hence
			\begin{equation*}
				N(e_1, \ldots, e_s) \prod_{i=1}^s (e_i +2)!
				\leq \prod_{i=1}^s (e_i +4)!
				.
			\end{equation*}
			We apply now \cref{lem:ingredient_for_sum_over_all_partitions} with $A_i = e_i -1$, $A = A_1 + \ldots + A_s = D - s$, and $r=2$. As $A_i \leq (D-s-2) -1 = A-3$ for every $i=1,\ldots, s$, this yields
			\begin{equation*}
				\prod_{i=1}^s (e_i +4)!
				\leq 5!^{s-1} \cdot (D-s+2)!
				.
			\end{equation*}
			The number of partitions $D = e_1 + \ldots + e_s$ where each $e_i$ is at most $D-s-2$ is again bounded by $(D-s-2)^{s-1}$.
			
			Adding up the five cases, we have
			\begin{align*}
				\sum_{D = e_1 + \ldots + e_s}
				& N(e_1, \ldots, e_s) \prod_{i=1}^s (e_i +2)! \\
				& = s \cdot s^2 \cdot 3^{s-1} \cdot (D - s +3)! \cdot 3^{s-1} \\
				& + s(s-1) \cdot (s+1)^2 \cdot 6 \cdot 3^{s-2} \cdot (D - s +2)! \cdot 4! \cdot 3!^{s-2} \\
				& + s(s-1) \cdot (s+2)^2 \cdot 10 \cdot 3^{s-2} \cdot (D-s+1)! \cdot 5! \cdot 3!^{s-2} \\
				& + s(s-1)(s-2) \cdot (s+2)^2 \cdot 6^2 \cdot 3^{s-3} \cdot (D-s+1)! \cdot 4!^2 \cdot 3!^{s-3} \\
				& + (D-s-2)^{s-1} \cdot 5!^{s-1} \cdot (D-s+2)! \\
				& \leq C_1^s \cdot (D +1)!
			\end{align*}
			for some constant $C_1^s$ which does not depend on $s$.
			
			With the same arguments as above for summing over all choices, we obtain a constant~$C$ as in the statement. Hence, the second bound is proven.
		\end{proof}
	\end{lem}
	
	In the calculations for the hyperelliptic components of $\mathcal{H}(2g-2)$ and $\mathcal{H}(g-1,g-1)$, 
	the dimension terms simplify but we have to consider more complicated volume terms.
	The sum over all configurations for multiplicity $2$ is then a sum over all partitions of~$g = g_1 + g_2$. To exemplify  the arguments in the lemma after next, we consider first a simpler sum.
	
	\begin{lem} \label{lem:partition_hyperelliptic}
		It holds
		\begin{equation*}
			\sum_{g_1 + g_2 = g} \frac{g^{\sfrac32}}{(g_1 \cdot g_2)^{\sfrac32}}
			= 2 \cdot \zeta(\sfrac32) \cdot \left( 1 + O \left( \frac{1}{g^{\sfrac14}} \right) \right)
			.
		\end{equation*}
		
		\begin{proof}
			For $g$ large enough, we can use the symmetry in $g_1$ and $g_2$ to obtain
			\begin{align*}
				\sum_{g_1 + g_2 = g} \frac{g^{\sfrac32}}{(g_1 \cdot g_2)^{\sfrac32}}
				& =\sum_{n=1}^{g-1} \frac{g^{\sfrac32}}{(n(g-n))^{\sfrac32}} \\
				& = 2 \cdot \sum_{n=1}^{\sfrac{g}{2}} \frac{1}{n^{\sfrac32} (1-\sfrac{n}{g})^{\sfrac32}} + O\left( \frac{1}{g^{\sfrac32}} \right) \\
				& = 2 \cdot \left( \sum_{n=1}^{g^{\sfrac56}} \frac{1}{n^{\sfrac32} (1-\sfrac{n}{g})^{\sfrac32}} + \sum_{n= g^{\sfrac56}}^{\sfrac{g}{2}} \frac{1}{n^{\sfrac32} (1-\sfrac{n}{g})^{\sfrac32}} \right) + O\left( \frac{1}{g^{\sfrac32}} \right) 
				.
			\end{align*}
			For $1 \leq n \leq g^{\sfrac56}$, it holds $1 - g^{-\sfrac16} \leq 1 - \sfrac{n}{g} \leq 1$, hence $\frac{1}{(1-\sfrac{n}{g})^{\sfrac32}} = 1 + O\left( \sfrac{1}{g^{\sfrac14}} \right)$ for all summands in the first sum.
			For $g^{\sfrac56} \leq n \leq \sfrac{g}{2}$, it holds $g^{\sfrac56} \cdot (1 - \sfrac12) \leq n(1 - \sfrac{n}{g})$, hence $\frac{1}{(n(1-\sfrac{n}{g}))^{\sfrac32}} = O\left( \sfrac{1}{g^{\sfrac54}} \right)$ for all summands in the second sum.
			Therefore, we obtain
			\begin{align*}
				\sum_{g_1 + g_2 = g} \frac{g^{\sfrac32}}{(g_1 \cdot g_2)^{\sfrac32}}
				& = 2 \cdot \sum_{n=1}^{g^{\sfrac56}} \frac{1 + O\left( \sfrac{1}{g^{\sfrac14}} \right) }{n^{\sfrac32}} 
				+ 2 \cdot \sum_{n= g^{\sfrac56}}^{\sfrac{g}{2}} O\left( \frac{1}{g^{\sfrac54}} \right) 
				+ O \left( \frac{1}{g^{\sfrac32}} \right) \\
				& = 2 \cdot \left( 1 + O\left( \frac{1}{g^{\sfrac14}} \right) \right) 
				\cdot \left( \zeta(\sfrac32) - \sum_{n= g^{\sfrac56}}^\infty \frac{1}{n^{\sfrac32}} \right)
				+ O(g) \cdot O\left( \frac{1}{g^{\sfrac54}} \right)
			\end{align*}
			where we use the definition of $\zeta(\sfrac32) = \sum_{n=1}^\infty \sfrac{1}{n^{\sfrac32}}$. 	
			We can estimate the error term related to this as
			\begin{equation*}
				\sum_{n= g^{\sfrac56}}^\infty \frac{1}{n^{\sfrac32}}
				= \int_{g^{\sfrac56}}^\infty \frac{1}{x^{\sfrac32}} \, dx + O\left( \frac{1}{(g^{\sfrac56} ) ^{\sfrac32}} \right)
				= \left[ -2 x^{-\sfrac12} \right]_{g^{\sfrac56}}^\infty + O\left( \frac{1}{g^{\sfrac54}} \right)
				= \frac{2}{g^{\sfrac{5}{12}}} + O\left( \frac{1}{g^{\sfrac54}} \right)
				.
			\end{equation*}
			
			Summarizing, we obtain
			\begin{align*}
				\sum_{g_1 + g_2 = g} \frac{g^{\sfrac32}}{(g_1 \cdot g_2)^{\sfrac32}}
				& = 2 \cdot \left( 1 + O\left( \frac{1}{g^{\sfrac14}} \right) \right) 
				\cdot \left( \zeta(\sfrac32) + O \left( \frac{1}{g^{\sfrac{5}{12}}}  \right)\right)
				+ O\left( \frac{1}{g^{\sfrac14}} \right) \\
				& = 2 \cdot \zeta(\sfrac32) \cdot \left( 1 + O\left( \frac{1}{g^{\sfrac14}} \right) \right)
			\end{align*}
			which was the desired statement.
		\end{proof}
	\end{lem}
	
	We apply now the arguments from the previous lemma to the sums that appear in calculations of the Siegel--Veech constants for $\mathcal{H}^\hyp(2g-2)$ and $\mathcal{H}^\hyp(g-1,g-1)$.
	
	\begin{lem} \label{lem:sum_with_double_factorials}
		It holds
		\begin{align*}
			\sum_{g_1 + g_2 = g} \! \frac{g^{\sfrac32}}{(2g_1 \! + \! 1) \! \cdot \! (2g_2 \! + \! 1)}
			\cdot \frac{(2g_1 \! - \! 1)!!}{(2g_1)!!}
			\cdot \frac{(2g_2 \! - \! 1)!!}{(2g_2)!!}
			& = \sqrt{\pi} \cdot \left( \frac{1}{2} - \frac{1}{\pi} \right)
			\cdot \left( \! 1 \! + \! O \! \left( \! \frac{1}{g^{\sfrac14}} \! \right) \! \right)
			\\
			\sum_{g_1 + g_2 = g-1} \! \frac{g^{\sfrac32}}{(2g_1 \! + \! 2) \! \cdot \! (2g_2 \! + \! 2)}
			\cdot \frac{(2g_1)!!}{(2g_1 \! + \! 1)!!}
			\cdot \frac{(2g_2)!!}{(2g_2 \! + \! 1)!!}
			& = \frac{\sqrt{\pi}}{4} 
			\cdot \left( \frac{\pi^2}{4} - 1 \right)		
			\cdot \left( \! 1 \! + \! O \! \left( \! \frac{1}{g^{\sfrac14}} \! \right) \! \right)
			\\
			\sum_{g_1 + g_2 = g-1}
			\! \frac{g^{\sfrac32}}{(2g_1 \! + \! 1) \! \cdot \! (2g_2 \! + \! 2)}
			\cdot \frac{(2g_1 \! - \! 1)!!}{(2g_1)!!}
			\cdot \frac{(2g_2)!!}{(2g_2 \! + \! 1)!!}
			& =  \frac{\sqrt{\pi}}{4} \! \cdot \! \left( \! \frac{3\pi}{4} \! - \! 1 \! - \! \frac{1}{\pi} \! \right)
			\! \cdot \! \left( \! 1 \! + \! O \! \left( \! \frac{1}{g^{\sfrac14}} \! \right) \! \right)
			. 
		\end{align*}
		
		\begin{proof}
			First note that the estimate $\sfrac{(2g-1)!!}{(2g)!!} = \sfrac{1}{\sqrt{\pi g}} \cdot (1 + O(\sfrac{1}{g})$ and a comparison with the sum in \cref{lem:partition_hyperelliptic} shows that all three sums converge for $g \to \infty$.
			
			For the first sum, by symmetry, we obtain
			\begin{align*}
				& \sum_{g_1 + g_2 = g} \frac{g^{\sfrac32}}{(2g_1+1) \cdot (2g_2+1)}
				\cdot \frac{(2g_1-1)!!}{(2g_1)!!}
				\cdot \frac{(2g_2-1)!!}{(2g_2)!!} \\
				& = 2 \sum_{n = 1}^{\sfrac{g}{2}}
				\frac{(2n-1)!!}{(2n+1) \cdot (2n)!!}
				\cdot \frac{g^{\sfrac32} \cdot (2g-2n-1)!!}{(2g -2n +1) \cdot (2g -2n)!!}
				+ O \left( \frac{1}{g^{\sfrac32}} \right)
				.
			\end{align*}
			As $\sfrac{(2g-1)!!}{(2g)!!} = \sfrac{1}{\sqrt{\pi g}} \cdot (1 + O(\sfrac{1}{g})$, we are in the same situation as in \cref{lem:partition_hyperelliptic} and obtain that the sum is equal to
			\begin{equation*}
				2 \cdot \frac{g^{\sfrac32}}{2g \! \cdot \! \sqrt{\pi g}} 
				\cdot \sum_{n = 1}^{\infty} \frac{(2n-1)!!}{(2n+1) \! \cdot \! (2n)!!}
				\cdot \left( \! 1 \! + \! O \! \left( \frac{1}{g^{\sfrac14}} \right) \! \! \right)
				= \frac{1}{\sqrt{\pi}}
				\cdot \sum_{n = 1}^{\infty} \frac{(2n-1)!!}{(2n+1) \! \cdot \! (2n)!!}
				\cdot \left( \! 1 \! + \! O \! \left( \frac{1}{g^{\sfrac14}} \right) \! \! \right)
				\! .
			\end{equation*}
			We consider now the series starting with $n=0$ (where we assign $(-1)!! = 0!! = 1$) and express the double factorials through the gamma function:
			\begin{align*}
				\sum_{n = 0}^{\infty} \frac{(2n-1)!!}{(2n+1) \cdot (2n)!!}
				& =  \sum_{n = 0}^{\infty} \frac{1}{2n+1} \cdot \frac{ 2^n \cdot \Gamma(n+ \sfrac12) \cdot \Gamma(\sfrac12)^{-1}}{ 2^n \cdot \Gamma(n+1)}
				\intertext{Using $\Gamma(\sfrac12)^2 = \pi$ and expressing $\frac{\Gamma(n+\sfrac12) \cdot \Gamma(\sfrac12)}{\Gamma(n+1)} = B(n+\sfrac12, \sfrac12) = \int_0^1 t^{n-\sfrac12}(1-t)^{-\sfrac12}\,dt$ as beta function yields}
				& =  \sum_{n = 0}^{\infty} \frac{1}{2n+1} \cdot \frac{1}{\pi} \cdot \frac{ \Gamma(n+ \sfrac12) \cdot \Gamma(\sfrac12)}{\Gamma(n+1)} \\
				& = \frac{1}{\pi} \cdot \sum_{n = 0}^{\infty} \frac{1}{2n+1} \cdot \int_0^1 t^{n- \sfrac12} (1-t)^{-\sfrac12} \, dt
				\intertext{As the original series converges, Fubini's theorem yields}
				& =  \frac{1}{\pi} \cdot \int_0^1 \frac{1}{t} \cdot \frac{1}{\sqrt{1-t}} \cdot \sum_{n = 0}^{\infty} \frac{\sqrt{t}^{2n+1}}{2n+1} \, dt \\
				& =  \frac{1}{\pi} \cdot \int_0^1 \frac{\tanh^{-1}(\sqrt{t})}{t \cdot \sqrt{1-t}} \, dt
				\intertext{Substituting $t \mapsto \tanh^2(x)$ yields}
				& =  \frac{1}{\pi} \cdot \int_0^\infty \frac{x}{\tanh^2(x) \cdot \sqrt{1-\tanh^2(x))}} \cdot 2 \tanh(x) \cdot (1-\tanh^2(x)) \, dx \\
				& =  \frac{2}{\pi} \cdot \int_0^\infty \frac{x \cdot \sqrt{1 - \tanh^2(x)}}{\tanh(x)} \, dx \\
				& = \frac{2}{\pi} \cdot \int_0^\infty \frac{x}{\sinh(x)} \, dx
				= \frac{2}{\pi} \cdot \int_0^\infty \frac{2x}{e^x - e^{-x}} \, dx
				= \frac{4}{\pi} \cdot \int_0^\infty \frac{x \cdot e^{-x}}{1- e^{-2x}} \, dx \\
				\intertext{Using an expression as geometric series and then again Fubini's theorem yields}
				& = \frac{4}{\pi} \cdot \int_0^\infty x \cdot e^{-x} \cdot \sum_{k=0}^\infty e^{-2kx} \, dx \\
				& = \frac{4}{\pi} \cdot \sum_{k=0}^\infty  \int_0^\infty x \cdot e^{-(2k+1)x} \, dx
				\intertext{We substitute $x \mapsto \frac{t}{2k+1}$ and obtain}
				& = \frac{4}{\pi} \cdot \sum_{k=0}^\infty  \int_0^\infty \frac{t}{2k+1} \cdot e^{-t} \cdot \frac{1}{2k+1} \, dt \\
				& = \frac{4}{\pi} \cdot \sum_{k=0}^\infty \frac{1}{(2k+1)^2} \cdot \int_0^\infty t \cdot e^{-t} \cdot \, dt \\
				& = \frac{4}{\pi} \cdot \sum_{k=0}^\infty \frac{1}{(2k+1)^2} \cdot \Gamma(2) \\
				& = \frac{4}{\pi} \cdot \left( \sum_{k=0}^\infty \frac{1}{(2k+1)^2} + \sum_{k=0}^\infty \frac{1}{(2k)^2} - \sum_{k=0}^\infty \frac{1}{(2k)^2} \right) \cdot 1! \\
				& = \frac{4}{\pi} \cdot \left( \sum_{k=0}^\infty \frac{1}{k^2} - \frac{1}{4} \cdot \sum_{k=0}^\infty \frac{1}{k^2} \right) \\
				& = \frac{4}{\pi} \cdot \frac{3}{4} \cdot \sum_{k=0}^\infty \frac{1}{k^2}
				= \frac{3}{\pi} \cdot \zeta(2)
				= \frac{3}{\pi} \cdot \frac{\pi^2}{6}
				= \frac{\pi}{2}
				.
			\end{align*}
			Hence we obtain that the first sum is equal to
			\begin{equation*}
				\frac{1}{\sqrt{\pi}} \cdot 
				\left( \sum_{n = 0}^{\infty} \frac{(2n-1)!!}{(2n+1) \cdot (2n)!!} -1 \right)
				\cdot \left(1 + O \left( \frac{1}{g^{\sfrac14}} \right) \right)
				= \frac{1}{\sqrt{\pi}} \cdot \left( \frac{\pi}{2} - 1 \right) 
				\cdot \left(1 + O \left( \frac{1}{g^{\sfrac14}} \right) \right)
				.
			\end{equation*}
			
			With similar arguments and $\sfrac{(2g-2)!!}{(2g-1)!!} = \sfrac{\sqrt{\pi g}}{2g} \cdot (1 + O(\sfrac{1}{g})$, we obtain for the second sum
			\begin{align*}
				& \sum_{g_1 + g_2 = g-1} \frac{g^{\sfrac32}}{(2g_1+2) \cdot (2g_2+2)}
				\cdot \frac{(2g_1)!!}{(2g_1+1)!!}
				\cdot \frac{(2g_2)!!}{(2g_2+1)!!} \\
				& = 2 \cdot \sum_{n = 1}^{\sfrac{g}{2}}
				\frac{(2n)!!}{(2n+2) \cdot (2n+1)!!}
				\cdot \frac{g^{\sfrac32} \cdot (2g-2n-2)!!}{(2g -2n) \cdot (2g -2n-1)!!}
				+ O \left( \frac{1}{g^{\sfrac32}} \right) \\
				& = 2 \cdot \frac{g^{\sfrac32} \cdot \sqrt{\pi g}}{2g \cdot 2g} 
				\cdot \sum_{n = 1}^{\infty}
				\frac{(2n)!!}{(2n+2) \cdot (2n+1)!!}
				\cdot \left( 1 + O \left( \frac{1}{g^{\sfrac14}} \right) \right)
				.
			\end{align*}
			Considering now the series starting with $n=0$ and performing similar steps as for the first sum, in particular using the two expressions of the beta function $B(n+1, \sfrac12)$, yields
			\begin{align*}
				\sum_{n = 0}^{\infty} \frac{(2n)!!}{(2n+2) \cdot (2n+1)!!}
				& =  \sum_{n = 0}^{\infty} \frac{1}{2(n+1)} \cdot \frac{ 2^n \cdot \Gamma(n+ 1) \cdot \Gamma(\sfrac12)}{ 2^{n+1} \cdot \Gamma(n+\sfrac32)} \\
				& = \frac14 \cdot  \sum_{n = 0}^{\infty} \frac{1}{n+1} \cdot \int_0^1 t^{n} (1-t)^{-\sfrac12} \, dt \\
				& = \frac14 \cdot  \int_0^1 \frac{1}{t \cdot \sqrt{1-t}} \cdot \sum_{n = 0}^{\infty} \frac{t^{n+1}}{n+1} \, dt \\
				& =  \frac14 \cdot \int_0^1 \frac{-\log(1-t)}{t \cdot \sqrt{1-t}} \, dt
				\intertext{Substituting $t \mapsto 1 - e^{-2x}$, we obtain}
				& =  \frac14 \cdot \int_0^\infty \frac{-(-2x)}{(1-e^{-2x}) \cdot \sqrt{e^{-2x}}} \cdot 2 e^{-2x} \, dx \\
				& = \int_0^\infty \frac{x \cdot  e^{-x}}{1-e^{-2x}} \, dx
				.
			\end{align*}
			Note that we evaluated this integral as $\sfrac34 \cdot \zeta(2) = \sfrac{\pi^2}{8}$ in the previous calculation and hence the second sum is equal to
			\begin{equation*}
				\frac{\sqrt{\pi}}{2} 
				\cdot \left( \sum_{n = 0}^{\infty}
				\frac{(2n)!!}{(2n+2) \cdot (2n+1)!!} - \frac12 \right)
				\cdot \left( \! 1 \! + \! O \! \left( \frac{1}{g^{\sfrac14}} \right) \! \right)
				= \frac{\sqrt{\pi}}{2} 
				\cdot \left( \frac{\pi^2}{8} - \frac{1}{2} \right)		
				\cdot \left( \! 1 \! + \! O \! \left( \frac{1}{g^{\sfrac14}} \right) \! \right)
				.
			\end{equation*}
			
			The last sum in the statement, we split into two sums where in the first one, we set~$n = g_1$ and in the second, we set $n = g_2$. Then we have
			\begin{align*}
				& \sum_{g_1 + g_2 = g-1}
				\frac{g^{\sfrac32}}{(2g_1+1) \cdot (2g_2+2)}
				\cdot \frac{(2g_1-1)!!}{(2g_1)!!}
				\cdot \frac{(2g_2)!!}{(2g_2+1)!!} \\
				& = \sum_{n = 1}^{\sfrac{g}{2}}
				\frac{(2n-1)!!}{(2n+1) \cdot (2n)!!}
				\cdot \frac{g^{\sfrac32} \cdot (2g-2n-2)!!}{(2g -2n) \cdot (2g -2n -1)!!} \\
				& + \sum_{n = 1}^{\sfrac{g}{2}}
				\frac{(2n)!!}{(2n+2) \cdot (2n+1)!!}
				\cdot \frac{g^{\sfrac32} \cdot (2g-2n-3)!!}{(2g -2n-1) \cdot (2g - 2n -2)!!} 
				+ O \left( \frac{1}{g^{\sfrac32}} \right)
				.
			\end{align*}
			As $\sfrac{(2g)!!}{(2g+1)!!} = \sfrac{\sqrt{\pi g}}{2g} \cdot (1 + O(\sfrac{1}{g})$ and $\sfrac{(2g-1)!!}{(2g)!!} = \sfrac{1}{\sqrt{\pi g}} \cdot (1 + O(\sfrac{1}{g})$, we are in the same situation as in \cref{lem:partition_hyperelliptic} and obtain that the sum is equal to
			\begin{equation*}
				\left( \frac{\sqrt{\pi}}{4} \cdot \sum_{n = 1}^{\infty}	\frac{(2n-1)!!}{(2n+1) \cdot (2n)!!}
				+ \frac{1}{2 \sqrt{\pi}} \cdot \sum_{n = 1}^{\infty} \frac{(2n)!!}{(2n+2) \cdot (2n+1)!!} \right)
				\cdot \left( 1 +  O \left( \frac{1}{g^{\sfrac14}} \right) \right)
				.
			\end{equation*}
			Note that the series are exactly the ones that we have calculated for the first and second sum and adding up the values yields
			\begin{equation*}
				\frac{\sqrt{\pi}}{4} \cdot \left( \frac{\pi}{2} -1 \right) + \frac{1}{2 \sqrt{\pi}} \cdot \left( \frac{\pi^2}{8} - \frac12 \right)
				= \frac{\sqrt{\pi}}{4} \cdot \left( \frac{3\pi}{4} - 1 - \frac{1}{\pi} \right)
				.
				\qedhere
			\end{equation*}
		\end{proof}
	\end{lem}
	
	The penultimate lemma is more a remark which justifies 
	that the error term $(1 + O(\sfrac1g) \cdot O(1)^p)$ is good enough, even if we intend to sum up over all $p \leq g$.
	
	\begin{lem}
		For fixed $n \in \mathbb{N}$, we have
		\begin{equation*}
			\left( 1 + O \left( \frac1g \right) \right)^n = 1 + O \left( \frac1g \right)
			.
		\end{equation*}
		
		Furthermore, it holds
		\begin{equation*}
			\sum_{p=1}^g \frac{1}{g^{2p-2}} \cdot \left( 1 + O \left( \frac1g \right) \cdot O(1)^p \right)
			= 1 + O \left( \frac1g \right)
			.
		\end{equation*}
		
		\begin{proof}
			For the first estimate, note that
			\begin{equation*}
				\left( 1 + O \left( \frac1g \right) \right)^n 
				= 1 + O \left( \frac{n}{g} \right)
				= 1 + O \left( \frac1g \right)
				.
			\end{equation*}
			
			For the second estimate, there exists a constant $C > 0$ such that
			\begin{align*}
				\sum_{p=1}^g \frac{1}{g^{2p-2}} \cdot \left( 1 + O \left( \frac1g \right) \cdot O(1)^p \right)
				& \leq 1 + \frac{C}{g} + \frac{1}{g^2} + \frac{C^2}{g^3} + \ldots + \frac{1}{g^{2g-2}} + \frac{C^g}{g^{2g-1}} \\
				& \leq 1 + \frac{C}{g} + \frac{1}{g^2} + \frac{1}{g} + \ldots + \frac{1}{g^{2g-2}} + \frac{1}{g^{g-1}} \\
				& = 1 + O \left( \frac1g \right)
			\end{align*}
			where we use that $C < g$ for $g$ large enough.
		\end{proof}
	\end{lem}
	
	While the volume estimate from \cref{eq:volume_low-genus} is the main source for the error term $(1+ O(\sfrac1g) \cdot O(1)^p)$, another source is more out of convenience and is explained in the final~lemma.
	
	\begin{lem}
		\label{lem:cancelling_factorials}
		For $p \leq g$, it holds
		\begin{equation*}
			\frac{(2g-2p)!}{(2g)!}
			= \frac{1}{(2g)^{2p}}
			\cdot \left( 1 + O \left( \frac1g \right) \cdot O(1)^p \right)
			.
		\end{equation*}
		
		\begin{proof}
			For $p=g$, the estimate follows immediately as the error term becomes exponential in $g$.
			
			For $p < g$, using Stirling's formula $n! = \sqrt{2\pi n} \cdot n^n \cdot e^{-n} \cdot ( 1 + O(\sfrac1n))$
			and $(1+\sfrac{p}{n})^n = e^p \cdot ( 1 + O(\sfrac1n) \cdot O(1)^p)$, we obtain
			\begin{align*}
				& \frac{(2g-2p)! \cdot (2g)^{2p}}{(2g)!} \\
				& = \sqrt{ \frac{2 \pi (2g-2p)}{2 \pi \cdot 2g}}
				\cdot \frac{(2g-2p)^{2g-2p} \cdot e^{2g} \cdot (2g)^{2p}}{e^{2g-2p} \cdot (2g)^{2g}}
				\cdot \left(1 + O \left(\frac1g \right) \right) \cdot \left( 1 + O \left( \frac{1}{g-p} \right) \right) \\
				& = \sqrt{ \frac{ (2g-2p)}{2g}}
				\cdot \frac{(2g-2p)^{2g-2p} \cdot e^{2p}}{(2g)^{2g-2p}}
				\cdot \left(1 + O \left(\frac1g \right) \cdot O(1)^p \right) \\
				& = e^{2p}
				\cdot \left( \frac{2g-2p}{2g} \right)^{2g-2p+ \sfrac12}
				\cdot \left(1 + O \left(\frac1g \right) \cdot O(1)^p \right) \\
				& = e^{2p}
				\cdot \left( 1 - \frac{p}{g} \right)^{2g-2p+ \sfrac12}
				\cdot \left(1 + O \left(\frac1g \right) \cdot O(1)^p \right) \\
				& = e^{2p} \cdot (e^{-p})^2
				\cdot \left(1 + O \left(\frac1g \right) \cdot O(1)^p \right)
				\cdot \left(1 + O \left(\frac1g \right) \cdot O(1)^p \right) \\
				& = 1 + O \left(\frac1g \right) \cdot O(1)^p
				.
				\qedhere
			\end{align*}
		\end{proof}
	\end{lem}

	\section{Siegel--Veech constants for saddle connections between distinct zeros} \label{sec:known_values_non-loops}
	
	In this section, we recall and refine the main results on asymptotics of Siegel--Veech constants for saddle connections between distinct zeros. We also carry out the special case of principal strata for easier reference.
	
	\begin{thm}[Saddle connections between distinct, fixed zeros]
		\label{thm:distinct_fixed_zeros}
		Let $\mathcal{H}$ be a stratum or the odd/even/non-hyperelliptic component of a non-connected stratum.
		Consider the configuration of a saddle connection of multiplicity $p$ between a fixed zero of order~$m_1$ and a different, fixed zero of order $m_2$.
		There are no saddle connections of multiplicity greater than $\min\{m_1, m_2\}+1$.
		
		For multiplicity $p \leq \min\{m_1, m_2\} +1$, the Siegel--Veech constant is
		\begin{align*}
			& c(m_1 \saddleconnection m_2, p\leq \min\{m_1, m_2\} +1, \fixedzeros, \mathcal{H})  \\
			& = (m_1+1)(m_2+1) \cdot \left( \frac{\pi^2}{6} \right)^{p-1}
			\cdot \frac{1}{(2g + \ell - 3)^{2p-2}}
			\cdot \left( 1 + O \left( \frac1g \right) \cdot O(1)^p \right) 
			.
		\end{align*}
		
		\begin{proof}
			For the ease of notation, we denote by $\mathcal{H}$ the (possibly non-connected) stratum and show the statement for this case. The statement for $\mathcal{H}^\oen$ follows then from \cref{prop:distinct_non-connected_strata}.
			
			For a fixed choice of $a_i = a_i' + a_i''$, we have $\sum_{i=1}^p (a_i' +1) = m_1 + 1$ and $\sum_{i=1}^p (a_i'' + 1) = m_2 +1$. As all the $a_i'$ and $a_i''$ are non-negative, there can be no saddle connection of multiplicity greater than $\min\{m_1, m_2\}+1$.
			
			In the case $p = 1$, plugging \cref{eq:volume} into \cref{eq:SV_distinct_multiplicity_one} yields (see also the proof of Zorich in \cite[Corollary~1]{aggarwal_20})
			\begin{align*}
				c(m_1 \saddleconnection m_2, p=1, \fixedzeros, \mathcal{H}) 
				& = (m_1 + m_2 +1 ) \cdot \frac{4 \cdot (m_1+1)(m_2+1)}{4 \cdot (m_1+m_2+1)}
				\cdot \left( \! 1 \! + \! O \! \left( \frac1g \right) \! \right) \\
				& = (m_1+1)(m_2+1) \cdot \left( 1 + O \left( \frac1g \right) \right)
				.
			\end{align*}
			
			Now let $p$ be fixed with $2 \leq p \leq \min\{m_1, m_2\}+1$.
			
			We call a configuration \emph{dominant} when all but one $\mathcal{H}_i$ are equal to $\mathcal{H}(0)$.
			Then one $\mathcal{H}_i$ (say, $\mathcal{H}_1$) contains all zeros of $\mathcal{H}$ but the fixed zeros of order $m_1$ and~$m_2$ and additionally a zero of order $a_1 = m_1+m_2 -2p+2$. In particular,
			\begin{equation*}
				\sfrac{d_1}{2} -1 = 2(g-p+1) + (\ell-2 +1)  -1 -1
				= 2g -2p + \ell -1
				.
			\end{equation*}
			Note that $m_1$ and $m_2$ already determine $a_1'$ and $a_1''$, hence there is only one dominant configuration up to cyclic reordering of the $\mathcal{H}_i$.
			Furthermore, in the labelled situation, we have $|\Gamma| = |\Gamma_-| = 1$.
			Hence with \cref{prop:distinct_non-connected_strata} and \cref{lem:cancelling_factorials}, the Siegel--Veech constant for the dominant configuration is
			\begin{align*}
				& \frac{(m_1+1)(m_2+1)}{|\Gamma| \cdot |\Gamma_-|}
				\cdot \left( \frac{\pi^{2}}{6} \right)^{p-1}
				\cdot \frac{(\sfrac{d_1}{2}-1)!}{(\sfrac{d}{2}-2)!} 
				\cdot \left( 1 + O \left( \frac1g \right) \right) \\
				& = (m_1+1)(m_2+1) \cdot \left( \frac{\pi^2}{6} \right)^{p-1}
				\cdot \frac{(2g-2p+\ell-1)!}{(2g + \ell -3)!}
				\cdot \left( 1 + O \left( \frac1g \right) \right) \\
				& = (m_1+1)(m_2+1) \cdot \left( \frac{\pi^2}{6} \right)^{p-1}
				\cdot \frac{1}{(2g + \ell -3)^{2p-2}}
				\cdot \left( 1 + O \left( \frac1g \right) \cdot O(1)^p \right) 
				.
			\end{align*}
			
			We follow now the strategy of \cite[Proposition 3.4]{aggarwal_19} to show that the sum of the Siegel--Veech constants over all non-dominant configurations is of a lower order than the Siegel--Veech constant of the dominant configuration.
			
			From \cref{prop:distinct_non-connected_strata}, we have that for a non-dominant configuration with fixed choices of $\mathcal{H}_i$  (including labelling the zeros) and $a_i = a_i' + a_i''$, the Siegel--Veech constant is
			\begin{equation*}
				(m_1+1)(m_2+1)
				\cdot \frac{\prod_{i=1}^p (\sfrac{d_i}{2}-1)!}{(\sfrac{d}{2}-2)!} 
				\cdot O(1)^p
				.
			\end{equation*}
			
			We now consider the dimension term.
			For every $i=1,\ldots, p$, let $g_i$ be the genus of the surfaces from $\mathcal{H}_i$ and $\ell_i$ the number of zeros except the one of order $a_i$. Then we have $\sfrac{d_i}{2} -1 = 2g_i + \ell_i -1$ for every $i=1, \ldots, p$ as well as $\sum_{i=1}^p g_i = g$ and $\sum_{i=1}^p \ell_i = \ell - 2$.
			
			Recall that we consider only non-dominant configurations, hence we can use
			\begin{equation*}
				\ell_i + a_i < \sum_{i=1}^p (\ell_i + a_i) = \ell - 2 + m_1 + m_2 +2 - 2p
			\end{equation*}
			for every $i= 1, \ldots, p$ in the following.
			First, we apply \cref{lem:comparison_product_factorials_of_two_partitions} with $A_i = 2g_i + \ell_i -2$, $B_i = 2 \ell_i + a_i$, and $r= 1$
			to obtain
			\begin{align*}
				\frac{\prod_{i=1}^p (\sfrac{d_i}{2}-1)!}{(\sfrac{d}{2} - 2)!}
				& \leq \frac{1}{(2g+ \ell -2p-1)^{2p-1}}
				\cdot \frac{\prod_{i=1}^p (2g_i + \ell_i -2 +1)!}{(2g +\ell -2 -2p)!} \\
				& \leq \frac{1}{(2g+ \ell -2p-1)^{2p-1}}
				\cdot \frac{\prod_{i=1}^p (2\ell_i + a_i +1)!}{(2\ell -4 +m_1 +m_2 +2 -2p)!} 
				.
			\end{align*}
			
			For a fixed partition of $\ell -2 = \ell_1 + \ldots + \ell_p$, we have
			\begin{equation*}
				\binom{\ell-2}{\ell_1} \cdot \binom{\ell-\ell_1-2}{\ell_2} \cdots \binom{\ell - \ell_1 - \ldots - \ell_{p-1} -2}{\ell_p} 
				= \frac{(\ell-2)!}{\prod_{i=1}^p \ell_i!}
			\end{equation*}
			choices of how to label the zeros in $\mathcal{H}_1, \ldots, \mathcal{H}_p$. We multiply the dimension term with the number of these choices and apply \cref{lem:product_binomials} with $r=1$ to obtain
			\begin{align*}
				& \frac{(\ell-2)!}{\prod_{i=1}^p \ell_i!}
				\cdot \frac{\prod_{i=1}^p (\sfrac{d_i}{2}-1)!}{(\sfrac{d}{2} -2)!} \\
				& \leq \frac{1}{(2g+\ell-2p-1)^{2p-1}}
				\cdot \frac{(\ell-2)!}{(2\ell -4 + m_1 +m_2 +2 -2p)!}
				\cdot \frac{ \prod_{i=1}^p (2\ell_i + a_i +1)!}{ \prod_{i=1}^p \ell_i! } \\
				& \leq \frac{1}{(2g+\ell-2p-1)^{2p-1}}
				\cdot 2^p \cdot  \frac{(\ell-2)!}{(2\ell -4 + m_1 +m_2 +2 -2p)!} \\
				& \quad
				\cdot \binom{2\ell -4 +m_1+m_2 +2-2p}{\ell -2}
				\cdot \prod_{i=1}^p (\ell_i + a_i +1)! \\
				& \leq \frac{2^p}{(2g+\ell-2p-1)^{2p-1}}
				\cdot \frac{\prod_{i=1}^p (\ell_i + a_i +1)!}{(\ell-2+m_1+m_2 +2 -2p)!}
				.
			\end{align*}
			
			To take into account now all non-dominant configurations, we write $a_i = a_i'+a_i''$. Then the sum over all possible $a_i'+a_i''$ and $\ell_i$ is
			\begin{equation*}
				\sum_{\substack{m_1 +1 -p = \\ a_1' + \ldots + a_p'}}
				\enspace
				\sum_{\substack{m_2 +1 -p = \\ a_1''+ \ldots + a_p''}}
				\enspace
				\sum_{\substack{\ell -2 = \\ \ell_1 + \ldots + \ell_p}}
				\frac{(\ell-2)!}{\prod_{i=1}^p \ell_i!}
				\cdot \frac{\prod_{i=1}^p (\sfrac{d_i}{2}-1)!}{(\sfrac{d}{2} -2)!}
				.
			\end{equation*}
			
			We insert the bound from above and apply the first bound of \cref{lem:sum_product_factorials_over_all_partitions} to obtain
			\begin{align*}
				& \frac{2^p}{(2g \! + \! \ell \! - \! 2p \! - \! 1)^{2p-1}}
				\cdot \frac{1}{(\ell \! - \! 2 \! + \! m_1 \! + \! m_2 \! + \! 2 \!  - \! 2p)!}
				\sum_{\substack{m_1 +1 -p = \\ a_1' + \ldots + a_p'}}
				\
				\sum_{\substack{m_2 +1 -p = \\ a_1''+ \ldots + a_p''}}
				\
				\sum_{\substack{\ell -2 = \\ \ell_1 + \ldots + \ell_p}}
				\prod_{i=1}^p (\ell_i + a_i' + a_i'' +1)! \\
				& \leq \frac{2^p}{(2g+\ell-2p-1)^{2p-1}}
				\cdot \frac{1}{(\ell-2+m_1+m_2 +2 -2p)!}
				\cdot C^p \cdot (\ell-2+m_1+m_2 +2-2p)! \\
				& = (2C)^p \cdot \frac{1}{(2g+\ell-2p-1)^{2p-1}}
			\end{align*}
			for some constant $C$ which does not depend on $p$.
			
			Multiplying this term with $(m_1+1)(m_2+1) \cdot O(1)^p$ yields the Siegel--Veech constant for all non-dominant configurations at once.
			Altogether, this yields for the Siegel--Veech constant
			\begin{align*}
				& c(m_1 \saddleconnection m_2, p \leq \min\{m_1,m_2\}+1, \fixedzeros, \mathcal{H}) \\
				& = (m_1+1)(m_2+1) \cdot \left( \frac{\pi^2}{6} \right)^{p-1} \cdot \frac{1}{(2g + \ell -3)^{2p-2}}
				\cdot \left( 1 + O \left( \frac1g \right) \cdot O(1)^p \right) \\
				& \quad
				+ (m_1+1)(m_2+1) \cdot \frac{ 1 }{(2g+ \ell -2p-1)^{2p-1}} \cdot O(1)^p \\
				& = (m_1+1)(m_2+1) \cdot \left( \frac{\pi^2}{6} \right)^{p-1} \cdot \frac{1}{(2g + \ell -3)^{2p-2}}
				\cdot \left( 1 + O \left( \frac1g \right) \cdot O(1)^p \right)
			\end{align*}
			which finishes the proof.
		\end{proof}
	\end{thm}
	
	From the theorem above, we can also directly deduce through summation that the Siegel--Veech constant for counting all saddle connections of \emph{any} multiplicity between a fixed zero of order $m_1$ and a different, fixed zero of order $m_2$ is
	\begin{align}
		\label{eq:saddle_connections_fixed_all_multiplicities}
		& c(m_1 \saddleconnection m_2, p= \any, \fixedzeros, \mathcal{H}) \notag \\
		& = \sum_{p=1}^{\min\{m_1, m_2\} +1}
		p \cdot
		(m_1+1)(m_2+1) \cdot \left( \frac{\pi^2}{6} \right)^{p-1}
		\cdot \frac{1}{(2g + \ell -3)^{2p-2}}
		\cdot \left( 1 + O \left( \frac1g \right) \! \cdot \! O(1)^p \right) \notag \\
		& = (m_1 + 1)(m_2 + 1) \cdot \left(1 + O \left( \frac{1}{g} \right) \right)
		.
	\end{align}
	
	This has previously been shown in \cite[Theorem 1.2]{aggarwal_19}.
	
	\bigskip
	
	Note that it is actually more natural in our context to not count all the saddle connections of higher multiplicity but only one per homology class. With this variation, Chen, Möller, Sauvaget, and Zagier have determined that the above asymptotics of the Siegel--Veech constant over all multiplicities is actually \emph{exact}!
	That is, for any connected component $\mathcal{H}$ of a stratum and for the configuration of homology classes of saddle connections between a fixed zero of order $m_1$ and a different, fixed zero of order $m_2$, the Siegel--Veech constant is \cite[Theorem 1.3]{chen_moeller_sauvaget_zagier_20}
	\begin{equation}
		\label{eq:exact_up_to_homology}
		c(m_1 \saddleconnection m_2, p=\any, \fixedzeros, \textnormal{up to homology}, \mathcal{H}) = (m_1+1)(m_2+1)
		.
	\end{equation}
	
	\bigskip
	
	From \cref{thm:distinct_fixed_zeros}, we obtain the Siegel--Veech constants for non-fixed zeros by multiplying with the number of  possibilities to choose the zeros.
	For $m_1 \neq m_2$, this is the number of zeros of order $m_1$ multiplied with the number of zeros of order $m_2$. If $m_1 = m_2$, we have to subtract $1$ from the number of zeros of order $m_2$ and divide by $2$ to account for the symmetry of exchanging the two zeros. 
	
	\bigskip
	
	As a special case, we consider the \emph{principal stratum} $\mathcal{H}(1, \ldots, 1)$ with $2g-2$ zeros of order $1$.
	To obtain the Siegel--Veech constant for a saddle connection between any two distinct zeros, we consider the Siegel--Veech constant from \cref{thm:distinct_fixed_zeros} with $m_1 = m_2 = 1$ and multiply by the combinatorial factor of choices for labelling the zeros, that is
	\begin{align} \label{eq:mult_any_distinct_unlabelled_zeros_principal}
		c(1 \saddleconnection 1, p= 1, \anyzeros, \mathcal{H}(1,\ldots, 1)) 
		& =  \frac{(2g-2)(2g-1)}{2} \cdot 4 \cdot \left(1 + O \left( \frac{1}{g} \right) \right) \notag \\
		& = 8g^2 \cdot \left(1 + O \left( \frac{1}{g} \right) \right)
	\end{align}
	and
	\begin{align} \label{eq:mult_two_distinct_unlabelled_zeros_principal}
		c(1 \saddleconnection 1, p= 2, \anyzeros, \mathcal{H}(1,\ldots, 1)) 
		& =  \frac{(2g-2)(2g-1)}{2} \cdot 4 \cdot \frac{\pi^2}{6} \cdot \frac{1}{(4g)^2} \left(1 + O \left( \frac{1}{g} \right) \right) \notag \\
		& = \frac{\pi^2}{12} \cdot \left(1 + O \left( \frac{1}{g} \right) \right).
	\end{align}
	
	Formulas for these cases have also been given in \cite[Formulas 8.3 and 9.2]{eskin_masur_zorich_03}.

	\section{Siegel--Veech constants for loops}
	\label{sec:known_values_loops}
	
	In the previous section, we focused on saddle connections between two distinct zeros. Now we consider \emph{loops}, that is, saddle connections from a zero back to the same~zero.
	
	A loop could either bound a cylinder or not, giving different cases for the formula for the Siegel--Veech constant, even for a loop of multiplicity $1$. In particular, if the loop bounds a cylinder, the formula involves an additional term for the possible choices of the zero on the second boundary of the cylinder.
	
	We start with the Siegel--Veech constants for multiplicity $1$ which can also be found in Zorich's appendix to \cite{aggarwal_20} for the case of connected strata.
	
	\begin{prop}[Loops of multiplicity $1$ for fixed zeros]
		\label{prop:mult_1_loops_labelled_zeros}
		Let $\mathcal{H}$ be a stratum or the odd/even/non-hyperelliptic component of a non-conneted stratum.
		Consider the configuration of a loop of multiplicity~$1$ at a fixed zero of order $m$.
		
		If the loop does not bound a cylinder, the Siegel--Veech constant is
		\begin{equation*}
			c(m \scloop, \mathrm{no\ cylinder}, p=1, \fixedzero, \mathcal{H}) 
			= \frac{(m + 1)(m - 1)}{2} \cdot \left(1 + O \left( \frac{1}{g} \right) \right)
			.
		\end{equation*}
		
		If the loop bounds a cylinder, whose other boundary component contains the same zero, the Siegel--Veech constant is
		\begin{equation*}
			c(m \loopcylinder, p=1, \fixedzero, \mathcal{H}) 
			= \frac{(m + 1)(m - 1)}{2 \cdot (2g + \ell -3)} \cdot \left(1 + O \left( \frac{1}{g} \right) \right)
			.
		\end{equation*}
		
		If the loop bounds a cylinder, whose other boundary component contains a different, fixed zero of order $m_2$, the Siegel--Veech constant is
		\begin{equation*}
			c(m \cylinder m_2, p=1, \fixedzeros, \mathcal{H}) 
			= \frac{(m + 1)(m_2 + 1)}{2g + \ell -3} \cdot \left(1 + O \left( \frac{1}{g} \right) \right)
			.
		\end{equation*}
		
		\begin{proof}
			We essentially repeat here the proof of Zorich from \cite[Corollaries 3--5]{aggarwal_20} while including non-connected strata and their odd/even/non-hyperelliptic components.
			For the ease of notation, we denote by $\mathcal{H}$ the stratum and show the statements for this case.
			The statements for $\mathcal{H}^\oen$ follow then from \cref{prop:loops_non-connected_strata}.
			
			Note that for loops of multiplicity $1$, all configurations are dominant and $|\Gamma| = 1$, hence the Siegel--Veech constant is by \cref{prop:loops_non-connected_strata}
			\begin{equation*}
				\frac{(m+1) (m_2+1)}{|\Gamma_-|}
				\cdot \frac{(\sfrac{d_1}{2} - 1)!}{(\sfrac{d}{2} - 2)!}
				\cdot \left( 1 + O \left( \frac1g \right) \right)
			\end{equation*}
			where the term $m_2+1$ only appears if the loop bounds a cylinder and the other boundary component contains a different zero of order $m_2$.
			
			There are three types of configurations, corresponding to the three statements: a two-hole construction without a cylinder, a figure-eight construction with a cylinder, and a two-hole construction with a cylinder.
			
			If we perform a two-hole construction on surfaces from $\mathcal{H}_1$ and do not glue in a cylinder, we obtain a loop that does not bound a cylinder. We have then $b_1' + b_1'' = m -2$ and $\sfrac{d_1}{2} - 1 = 2g-2 + \ell +1 -2$. Furthermore, there are $m-1$ choices for the partition~$m-2 = b_1' + b_1''$ of which always two give the same configuration (except for~$b_1' = b_1''$) and we have $|\Gamma_-| = 1$ (except for $b_1' = b_1''$ where we have $|\Gamma_-| = 2$).
			Therefore, the Siegel--Veech constant is
			\begin{equation*}
				\frac{m-1}{2}
				\cdot (m+1)
				\cdot \frac{(2g + \ell -3)!}{(2g + \ell -3)!}
				\cdot \left( 1 + O \left( \frac1g \right) \right)
				= \frac{(m+1)(m-1)}{2}
				\cdot \left( 1 + O \left( \frac1g \right) \right)
				.
			\end{equation*}
			
			If we perform a figure-eight construction on surfaces from $\mathcal{H}_1$, we need to glue in a cylinder and hence obtain a loop which bounds a cylinder but its other boundary component is a loop at the same zero.
			We have then $a_1 = m -2$ and $\sfrac{d_1}{2} - 1 = 2g-2 + \ell -2$. Furthermore, there are $m-1$ choices for the partition $m-2 = a_1' + a_1''$ of which always two give the same configuration (except for $a_1' = a_1''$) and we have $|\Gamma_-| = 1$ (except for $a_1' = a_1''$ where we have $|\Gamma_-| = 2$).
			Therefore, the Siegel--Veech constant is
			\begin{equation*}
				\frac{m-1}{2}
				\cdot (m +1)
				\cdot \frac{(2g + \ell -4)!}{(2g + \ell -3)!}
				\cdot \left( 1 + O \left( \frac1g \right) \right)
				= \frac{(m+1)(m-1)}{2 \cdot (2g+ \ell-3)}
				\cdot \left( 1 + O \left( \frac1g \right) \right)
				.
			\end{equation*}
			
			If we perform a two-hole construction on surfaces from $\mathcal{H}_1$ and glue in a cylinder, we obtain two loops which bound a cylinder and are at different zeros. Let $m_2$ be the order of the second zero.
			We have then $b_1' = m -1$, $b_1'' = m_2 -1$, and $\sfrac{d_1}{2} - 1 = 2g-2 + \ell -2$. Furthermore, there are no choices for $(b_1', b_1'')$ and we have $|\Gamma_-| = 1$.
			Therefore, the Siegel--Veech constant is
			\begin{equation*}
				(m+1)(m_2+1)
				\cdot \frac{(2g + \ell -4)!}{(2g + \ell -3)!}
				\cdot \left( 1 + O \left( \frac1g \right) \right) \\
				= \frac{(m+1)(m_2+1)}{2g+ \ell-3}
				\cdot \left( 1 + O \left( \frac1g \right) \right)
				. \qedhere
			\end{equation*}
		\end{proof}
	\end{prop}
	
	To obtain the Siegel--Veech constant for loops of multiplicity $1$ without further specification, we have to add up the Siegel--Veech constants for all possible configurations. In particular, we have to consider all choices of zeros $z_2$ on the other boundary of a potential cylinder. With $\sum_{z_2 \neq z_1} (m_2 + 1) = (2g-2 -m) + (\ell -1)$, we have
	\begin{align} \label{eq:mult_one_loop_labelled}
		& c(m \scloop, p=1, \fixedzero, \mathcal{H}) \notag \\
		& = \left( \frac{(m + 1)(m - 1)}{2} 
		+ \frac{(m + 1)(m - 1)}{2 (2g + \ell -3)}
		+ \sum_{z_2 \neq z_1} \frac{(m + 1)(m_2 + 1)}{2g + \ell -3} \right)
		\cdot \left(1 + O \left( \frac{1}{g} \right) \right) \notag \\
		& = \left( \frac{(m + 1)(m - 1)}{2} 
		+ (m+1) \cdot \frac{2g + \ell -3 -m}{2g + \ell -3} \right)
		\cdot \left(1 + O \left( \frac{1}{g} \right) \right) \notag \\
		& = \left( \frac{(m + 1)(m - 1)}{2} 
		+ m+1 - \frac{(m + 1) \cdot m}{2g + \ell -3} \right)
		\cdot \left(1 + O \left( \frac{1}{g} \right) \right) \notag \\
		& = \frac{(m+1)^2}{2} \cdot \left(1 + O \left( \frac{1}{g} \right) \right)
		.
	\end{align}
	
	The Siegel--Veech constants for loops bounding a cylinder can also be expressed through \emph{area Siegel--Veech constants} \cite{vorobets_05} for which more large-genus asymptotics are known. However, as Siegel--Veech constants for loops \emph{not} bounding a cylinder are of higher order for all zeros which are not of order $1$, we will not follow this route in this article.
	
	\bigskip
	
	In the principal stratum $\mathcal{H}(1, \ldots, 1)$, for geometric reasons, any loop has to bound a cylinder and the other boundary component of the cylinder cannot contain the same zero. By \cref{prop:mult_1_loops_labelled_zeros}, the Siegel--Veech constant for a loop of multiplicity $1$ at a fixed zero $z_1$ is
	\begin{align*}
		c(1 \scloop, p=1, \fixedzero, \mathcal{H}(1,\ldots, 1)) 
		& = \sum_{z_2 \neq z_1} \frac{4}{2g + 2g-2 -3}
		\cdot \left( 1+ O \left( \frac1g \right)\right) \\
		& = (2g-3) \cdot \frac{4}{4g-5}
		\cdot \left( 1+ O \left( \frac1g \right)\right) \\
		& = 2 \cdot \left( 1+ O \left( \frac1g \right)\right)
	\end{align*}	
	and for any zero is (see also the formula in \cite[Section 13.5.1]{eskin_masur_zorich_03})
	\begin{align} \label{eq:loop_mult_1_simple_zero}
		c(1 \scloop, p =1, \anyzero, \mathcal{H}(1,\ldots, 1))
		& = \frac{(2g-2)(2g-3)}{2} \cdot \frac{2 \cdot 2}{4g -5} \cdot \left(1 + O \left( \frac{1}{g} \right) \right) \nonumber \\
		& = 2g  \cdot \left(1 + O \left( \frac{1}{g} \right) \right)
		.
	\end{align}
	
	We see here that for loops, we can not just sum over the Siegel--Veech constants for fixed zeros to obtain the Siegel--Veech constant for any zero. This is because 
	loops can be homologous although they start at different zeros.
	Hence when summing the Siegel--Veech constant for fixed zeros over all zeros, we count some loops twice.
	
	This is only one of the issues that makes calculating the Siegel--Veech constant for loops more complicated than for saddle connections between distinct zeros.
	Another one is the phenomenon that the order of other zeros play a role for the Siegel--Veech constant of a fixed zero.
	However, we can still show generally that the Siegel--Veech constant for multiplicity $1$ dominates the Siegel--Veech constants for higher multiplicities (see also \cite[Appendix A.1]{vallejos_24} for the case that the order of the largest zero in the stratum is bounded and \cite[Proposition~8.2.5]{reichert_24} for the case of the minimal stratum.).
	
	We first show bounds for the asymptotics in the situation where all zeros that appear on the homologous loops are specified.

	\begin{prop}[Loops of higher multiplicity with all zeros fixed]
		\label{prop:loops_all_zeros_fixed}
		Let $\mathcal{H}$ be a stratum or the odd/even/non-hyperelliptic component of a non-connected stratum.
		Consider the configuration of a loop of multiplicity $p$ at exactly $n$ fixed zeros of order $m_1, \ldots, m_{n}$.
		Let $M \coloneqq m_1 + \ldots + m_n$.
		There are no loops of multiplicity smaller than $\sfrac{n}{2}$ or greater than $\sfrac{M}{2}$.
		
		If $n =1$, the Siegel--Veech constant for multiplicity $p \leq \sfrac{m_1}{2}$ is
		\begin{equation*}
			c(m_1 \scloop, p \leq \sfrac{m_1}{2}, \fixedzero, \mathcal{H})
			= \frac12
			\cdot \left( \frac{\pi^2}{6} \right)^{p-1}
			\cdot \frac{(m_1 \! + \! 1)(m_1 \! - \! 2p \! + \! 1)}{(2g \! + \! \ell \! - \!3)^{2p-2}}
			\cdot \left( \! 1 \! + \! O \! \left( \frac1g \right) \! \cdot \! O(1)^p \! \right)
			\! .
		\end{equation*}
		
		If $n \geq 2$, the Siegel--Veech constant for multiplicity $p$ with $\sfrac{n}{2} \leq p \leq \sfrac{M}{2}$ is
		\begin{equation*}
			c(m_1, \ldots, m_p \scloop, \sfrac{n}{2} \leq p \leq \sfrac{M}{2}, \fixedzeros, \mathcal{H})
			= \frac{(m_1+1) \cdots (m_n+1)}{(2g + \ell - 3)^{2p-3 +n}} \cdot O(1)^p
			.
		\end{equation*}
		
		\begin{proof}
			For the ease of notation, we denote by $\mathcal{H}$ the (possibly non-connected) stratum and show the statement for this case. The statement for $\mathcal{H}^\oen$ follows then from \cref{prop:loops_non-connected_strata}.
			
			Note that $M = \sum (a_i +2) + \sum (b_k' + b_k'' +2) = \sum a_i + \sum (b_k' + b_k'') +2p$. As the $a_i$ and~$(b_k', b_k'')$ are non-negative, this implies $M \geq 2p$. Furthermore, 
			the number of new zeros is equal to the number of two-hole constructions and cylinders, hence $n \leq 2p$.
			
			For $p=1$, we therefore have either one or two zeros. Hence, \cref{prop:mult_1_loops_labelled_zeros} implies the desired statement.
			
			Now let $p$ be fixed with $\sfrac{n}{2} \leq p \leq \sfrac{M}{2}$ and let $p'$ be the number of two-hole constructions and $q$ the number of cylinders. Then we have $n = p' + q$.
			
			We call a configuration \emph{dominant} when all but one $\mathcal{H}_i$ are equal to $\mathcal{H}(0)$ or $\mathcal{H}(0,0)$.
			Then one $\mathcal{H}_i$ (say, $\mathcal{H}_1$) contains all zeros of $\mathcal{H}$ but the $n$ fixed zeros of orders $m_1, \ldots, m_n$ and additionally a zero of order $a_1 = M -2p$ or two zeros $b_1'$ and $b_1''$ with $b_1' + b_1'' = M-2p$. In particular, all other $a_i$ or $(b_k', b_k'')$ are equal to $0$ and $g_1 = g-p$.
			By \cref{prop:loops_non-connected_strata}, the Siegel--Veech constant for a dominant configuration is
			\begin{equation*}
				\frac{(m_1+1) \cdots (m_n+1)}{|\Gamma| \cdot |\Gamma_-|}
				\cdot \left( \frac{\pi^2}{6} \right)^{p-1}
				\cdot \frac{\prod_{i=1}^p (\sfrac{d_i}{2} - 1)!}{(\sfrac{d}{2} - 2)!}
				\cdot \left( 1 + O \left( \frac1g \right) \right)
				.
			\end{equation*}
			
			If we perform a figure-eight construction on $\mathcal{H}_1$, then the dimension term is
			\begin{align*}
				\frac{\prod_{i=1}^p (\sfrac{d_i}{2}-1)!}{(\sfrac{d}{2}-2)!}
				= \frac{(\sfrac{d_1}{2} -1)! \cdot 2^{p'}}{(2g + \ell -3)!}
				& = \frac{(2g-2p +\ell -n +1 -2)! \cdot 2^{p'}}{(2g + \ell -3)!} \\
				& = \frac{2^{p'}}{(2g + \ell -3)^{2p-2 +n}}
				\cdot \left( 1 + O \left( \frac1g \right) \cdot O(1)^p \right)
			\end{align*}
			and we have $M - 2p+1$ choices for the partition $M - 2p = a_1 = a_1' + a_1''$. Furthermore, we have at most $n!$ choices for the labelling of the $n$ zeros. 
			Hence, the Siegel--Veech constant for these configurations is
			\begin{equation*}
				(M-2p+1) \cdot \frac{(m_1+1) \cdots (m_n+1)}{(2g + \ell -3)^{2p-2 +n}} \cdot O(1)^p
				.
			\end{equation*}
			
			If we perform a two-hole construction on $\mathcal{H}_1$, then the dimension term is
			\begin{align*}
				\frac{\prod_{i=1}^p (\sfrac{d_i}{2}-1)!}{(\sfrac{d}{2}-2)!}
				& = \frac{(2g -2p +\ell -n +2 -2)! \cdot 2^{p'-1}}{(2g + \ell -3)!} \\
				& = \frac{2^{p'-1}}{(2g + \ell -3)^{2p-3 +n}}
				\cdot \left( 1 + O \left( \frac1g \right) \cdot O(1)^p \right)
				.
			\end{align*}
			We distinguish now the cases when $n=1$ and when $n>1$.
			
			If $n=1$, then we have no cylinders and we perform figure-eight constructions on all $\mathcal{H}_i$ but $\mathcal{H}_1$.
			Furthermore, we have $m_1-2p+1$ choices for the partition of $m_1 -2p = b_1' + b_1''$ where always two choices yield the same configuration (except for $b_1' = b_1''$).
			Also, we have~$|\Gamma| = |\Gamma_-| = 1$ (except for $b_1' = b_1''$ where we have $|\Gamma_-| = 2$).
			Hence, the Siegel--Veech constant for these configurations is
			\begin{align*}
				& \frac{m_1 - 2p + 1}{2}
				\cdot (m_1+ 1)
				\cdot \left( \frac{\pi^2}{6} \right)^{p-1}
				\cdot \frac{1}{(2g + \ell - 3)^{2p-2}}
				\cdot \left( 1 + O \left( \frac1g \right) \cdot O(1)^p \right) \\
				& =
				\frac12
				\cdot \left( \frac{\pi^2}{6} \right)^{p-1}
				\cdot \frac{(m_1+1)(m_1 -2p+1)}{(2g + \ell -3)^{2p-2}}
				\cdot \left( 1 + O \left( \frac1g \right) \cdot O(1)^p \right)
				.
			\end{align*}
			
			If $n>1$, then also $b_1'$ and $b_1''$ are fixed. However, we have up to $n!$ choices for the labelling of the $n$ zeros.
			Hence  the Siegel--Veech constant for these configurations is
			\begin{equation*}
				\frac{(m_1+1) \cdots (m_n+1)}{(2g + \ell -3)^{2p-3 +n}} \cdot O(1)^p
				.
			\end{equation*}
			
			Taking also into account that the number of choices for the number and position of cylinders and two-hole constructions is in $O(1)^p$, we have for all dominant configurations together the Siegel--Veech constants for $n=1$
			\begin{align*}
				& c(m_1 \scloop, p, \fixedzero, \dominant, \mathcal{H}) \\
				& = \frac12
				\cdot \left( \frac{\pi^2}{6} \right)^{p-1}
				\cdot \frac{(m_1+1)(m_1 -2p+1)}{(2g + \ell -3)^{2p-2}}
				\cdot \left( 1 + O \left( \frac1g \right) \cdot O(1)^p \right) \\
				& \quad + (m_1-2p+1) \cdot \frac{m_1+1}{(2g + \ell -3)^{2p-1}} \cdot O(1)^p \\
				& = \frac12
				\cdot \left( \frac{\pi^2}{6} \right)^{p-1}
				\cdot \frac{(m_1+1)(m_1 -2p+1)}{(2g + \ell -3)^{2p-2}}
				\cdot \left( 1 + O \left( \frac1g \right) \cdot O(1)^p \right)
			\end{align*}
			and for $n \geq 2$
			\begin{align*}
				& c(m_1, \ldots, m_n \scloop, p, \fixedzeros, \dominant, \mathcal{H}) \\
				& =		\frac{(m_1+1) \cdots (m_n+1)}{(2g + \ell -3)^{2p-3 +n}} \cdot O(1)^p
				+ 		(M-2p+1) \cdot \frac{(m_1+1) \cdots (m_n+1)}{(2g + \ell -3)^{2p-2 +n}} \cdot O(1)^p \\
				& = \frac{(m_1+1) \cdots (m_n+1)}{(2g + \ell -3)^{2p-3 +n}} \cdot O(1)^p
				.
			\end{align*}
			
			Note that in a dominant configuration, only one or two zeros can have an order greater than $2p-2$, hence not for every $n \geq 2$ and every choice of $m_1, \ldots, m_n$, dominant configurations exist.
			
			Consider now any configuration which is not a cyclic reordering of the dominant configuration.
			We fix again the number $p'$ of two-hole constructions and the number $q$ of cylinders as well as the position of the two-hole constructions and cylinders.
			For every $i=1,\ldots, p$, let $g_i$ be the genus of surfaces from $\mathcal{H}_i$ and let $\ell_i$ be the number of zeros which are not the ones of order $a_i$ or $(b_k', b_k'')$.
			
			From \cref{prop:loops_non-connected_strata}, we have that for fixed choices of $\mathcal{H}_i$ (including labelling the zeros) and $a_i = a_i' + a_i''$ and $(b_k', b_k'')$, the Siegel--Veech constant for a non-dominant configuration is
			\begin{equation*}
				(m+1) \cdots (m_n+1)
				\cdot \frac{\prod_{i=1}^p (\sfrac{d_i}{2}-1)!}{(\sfrac{d}{2}-2)!} 
				\cdot O \left( 1 \right)^p
				.
			\end{equation*}
			
			To be able to later apply \cref{lem:sum_product_factorials_over_all_partitions} and to highlight the similarities with the proof in \cref{thm:distinct_fixed_zeros}, we abuse notation and define $b_i' = a_i$ and $b_i'' = 0$ if we perform a figure-eight construction on $\mathcal{H}_i$.
			We consider now the dimension term together with the count of partitions of $a_i = a_i' + a_i''$, that is
			\begin{equation*}
				\prod (a_i +1) \cdot \frac{\prod_{i=1}^p (\sfrac{d_i}{2}-1)!}{(\sfrac{d}{2} -2)!}
				.
			\end{equation*}
			If we perform a figure-eight construction on surfaces from $\mathcal{H}_i$, we have
			\begin{equation*}
				(a_i +1) \cdot (\sfrac{d_i}{2}-1)!
				= (a_i +1) \cdot (2g_i + \ell_i +1 -2)!
				\leq (2g_i + \ell_i)!
			\end{equation*}
			and if we perform a two-hole construction on surfaces from $\mathcal{H}_i$, we have
			\begin{equation*}
				(\sfrac{d_i}{2}-1)!
				= (2g_i + \ell_i +2 -2)!
				= (2g_i + \ell_i)!
				.
			\end{equation*}
			
			Recall that we consider only non-dominant configurations, hence we can use 
			\begin{equation*}
				\ell_i + b_i' + b_i'' < \sum_{i=1}^p (\ell_i + b_i' + b_i'') = \ell - n + M - 2p
			\end{equation*}
			for every $i= 1, \ldots, p$ in the following.	
			First, we apply \cref{lem:comparison_product_factorials_of_two_partitions} with $A_i = 2g_i + \ell_i -2$, $B_i = 2\ell_i +b_i' + b_i''$, and $r=2$ to obtain
			\begin{align*}
				\prod (a_i +1) \cdot \frac{\prod_{i=1}^p (\sfrac{d_i}{2}-1)!}{(\sfrac{d}{2}-2)!}
				& \leq \frac{1}{(2g+\ell -n-2p)^{2p-2+n}} \cdot \frac{\prod_{i=1}^p (2g_i + \ell_i -2 +2)!}{(2g+\ell-n-2p-2+1)!} \\
				& \leq \frac{1}{(2g+\ell-n-2p)^{2p-2+n}} \cdot \frac{\prod_{i=1}^p (2\ell_i + b_i' + b_i'' +2)!}{(2\ell - 2n + M -2p +1)!}
				.
			\end{align*}
			
			For a fixed partition of $\ell -n = \ell_1 + \ldots + \ell_p$, we have $\frac{(\ell-n)!}{\prod_{i=1}^p \ell_i!}$ choices of how to label the zeros in $\mathcal{H}_1, \ldots, \mathcal{H}_p$.
			We multiply the above term with the number of these choices and apply \cref{lem:product_binomials} with $r=2$ to obtain
			\begin{align*}
				& \frac{(\ell-n)!}{\prod_{i=1}^p \ell_i!}
				\cdot \prod (a_i + 1)
				\cdot \frac{\prod_{i=1}^p (\sfrac{d_i}{2}-1)!}{(\sfrac{d}{2} -2)!} \\
				& \leq \frac{1}{(2g + \ell -n -2p)^{2p-2+n}}
				\cdot \frac{(\ell-n)!}{(2\ell - 2n + M -2p +1)!}
				\cdot \frac{\prod_{i=1}^p (2\ell_i + b_i' + b_i'' + 2)!}{\prod_{i=1}^p \ell_i!} \\
				& \leq \frac{1}{(2g + \ell -n -2p)^{2p-2+n}}
				\cdot \frac{(\ell-n)!}{(2\ell - 2n + M -2p +1)!} \\
				& \quad \cdot 2^{2p} \cdot \binom{2\ell -2n + M-2p}{\ell - n}
				\cdot \prod_{i=1}^p (b_i' + b_i'' + \ell_i + 2)! \\
				& \leq \frac{ 2^{2p}}{(2g + \ell -n -2p)^{2p -2 + n}}
				\cdot \frac{\prod_{i=1}^p (b_i' + b_i'' + \ell_i + 2)!}{(2\ell - 2n + M -2p +1) \cdot (\ell - n +M -2p)!} \\
				& \leq \frac{ 2^{2p} }{(2g + \ell -n -2p)^{2p -2 + n}}
				\cdot \frac{\prod_{i=1}^p (b_i' + b_i'' + \ell_i + 2)!}{(\ell - n +M -2p +1)!}
				.
			\end{align*}
			
			We want to sum this term now over all possible choices of $(b_i', b_i'')$ and $\ell_i$.
			For this, note that every $b_i'$ and $b_i''$ appears as a block in exactly one partition of the form $M_j \coloneqq m_j -k_j$ (where the $k_j$ is given by the gluing pattern of the subsurfaces with boundary). Furthermore, if $b_i'$ and $b_i''$ are not equal to $0$, they do not appear in the same partition for any $i = 1, \ldots, p$ -- except if we have exactly one two-hole construction and no cylinders.
			We exclude the latter case for a moment and can apply then the second bound from \cref{lem:sum_product_factorials_over_all_partitions} to obtain
			\begin{align*}
				& \sum_{\substack{\text{partition} \\ \text{of } M_1}}
				\ldots
				\sum_{\substack{\text{partition} \\ \text{of } M_n}}
				\thickspace
				\sum_{\substack{\ell - n = \\ \ell_1 + \ldots + \ell_p}}
				\frac{(\ell-n)!}{\prod_{i=1}^p \ell_i!}
				\cdot \prod (a_i +1)
				\cdot \frac{\prod_{i=1}^p (\sfrac{d_i}{2}-1)!}{(\sfrac{d}{2} -2)!} \\
				& \leq \sum_{\substack{\text{partition} \\ \text{of } M_1}}
				\ldots
				\sum_{\substack{\text{partition} \\ \text{of } M_n}}
				\thickspace
				\sum_{\substack{\ell - n = \\ \ell_1 + \ldots + \ell_p}}
				\frac{2^{2p} }{(2g + \ell -n -2p)^{2p + n -2}}
				\cdot \frac{\prod_{i=1}^p (b_i' + b_i'' + \ell_i + 2)!}{(\ell - n +M -2p +1)!} \\
				& \leq
				\frac{ 2^{2p} }{(2g \! + \! \ell \! - \! n \! - \! 2p)^{2p + n -2}}
				\cdot \frac{1}{(\ell \! - \! n \! + \! M \! - \! 2p \! + \! 1)!}
				\cdot C^p \cdot (M_1 + \ldots + M_n + \ell - n + 1)! \\
				& = (4C)^p
				\cdot \frac{1}{(2g + \ell -n -2p)^{2p -2 +n}}
			\end{align*}
			for some constant $C$ which does not depend on $p$.
			
			Multiplying this term with $(m+1)\cdots(m_n+1) \cdot O(1)^p$ and summing over all $O(1)^p$ possibilities to choose the number and position of cylinders and of two-hole constructions yields
			\begin{align*}
				c(m_1, ..., m_n \scloop, \sfrac{n}{2} \! \leq \! p \! \leq \! \sfrac{M}{2}, \fixedzeros, \nondominant, \mathcal{H})
				= \frac{(m_1 \! + \! 1) \cdots (m_n \! + \! 1)}{(2g \! + \! \ell \! - \! n \! - \! 2p)^{2p -2 +n}}
				\cdot O(1)^p
				.
			\end{align*}
			
			In particular, the order of the Siegel--Veech constant for the dominant configurations is greater than the one for the non-dominant configurations.
			
			We now consider the case $p' = 1$ and $q=0$. Assume that the two-hole construction is performed on surfaces from $\mathcal{H}_1$. We have to take into account that $b_1'$ and $b_1''$ belong to the same partition $M_1 = m_1 -2p = \sum_{i=1}^p$, hence we have to	multiply the term above not only with $(m+1) \cdot O(1)^p$ but also with the number of choices of $(b_1', b_1'')$ which is at most $m_1 - 2p + 1$. Hence the Siegel--Veech constant is bounded as
			\begin{align*}
				c(m_1 \scloop, p \leq \sfrac{m_1}{2}, \fixedzero, \nondominant, \mathcal{H})
				\leq (m_1 \! - \! 2p \! + \! 1) \cdot \frac{m_1 \! + \! 1}{(2g \! + \! \ell \! - \! 2p \! - \! 1)^{2p - 1}}
				\cdot O(1)^p
			\end{align*}
			which is of the same order as the error term of the Siegel--Veech constant for the dominant configurations.
		\end{proof}
	\end{prop}

	We can use this theorem now to deduce estimates for the Siegel--Veech constant where we fix only one zero and allow homologous loops at arbitrary other zeros.
	
	\begin{thm}[Loops of higher multiplicity with one fixed zero]
		\label{thm:loops_one_fixed_zero}
		Let $\mathcal{H}$ be a stratum or the odd/even/non-hyperelliptic component of a non-connected stratum.
		Consider the configuration of a loop of multiplicity $p$ at a fixed zero of order $m$.
		Then the Siegel--Veech constant is
		\begin{equation*}
			c(m \scloop, p, \onefixedzero, \mathcal{H})
			= \frac{(m-2p+1) \cdot (m +1)}{(2g+\ell-3)^{2p-2}} \cdot O(1)^p
			.
		\end{equation*}
		
		Furthermore, if $m \geq \sfrac{g}{C}$ for some constant $C > 0$, then the Siegel--Veech constant is
		\begin{equation*}
			c(m \scloop, p \leq \sfrac{m}{2}, \onefixedzero, \mathcal{H})
			=  \frac12 \cdot \left( \frac{\pi^2}{6} \right)^{p-1} \cdot \frac{(m \! +\! 1)(m \! -\! 2p \!+ \! 1)}{(2g \!+ \!\ell \!- \! 3)^{2p-2}}
			\cdot \left( \! 1 \! + \! O \! \left( \frac1g \right) \! \cdot \! O(1)^p \! \right)
			\! .
		\end{equation*}
		
		\begin{proof}
			We have to sum the Siegel--Veech constants for all $n \leq 2p$ and for all choices of zeros $z_2, \ldots, z_n$ of orders $m_2, \ldots, m_n$, that is
			\begin{align*}
				& c(m \scloop, p, \onefixedzero, \mathcal{H}) \\
				& \leq \frac12 \cdot \left( \frac{\pi^2}{6} \right)^{p-1} \cdot \frac{(m +1)(m-2p+1)}{(2g+\ell-3)^{2p-2}}
				\cdot \left( 1 + O \left( \frac1g \right) \cdot O(1)^p \right) \\
				& \quad + \sum_{n=2}^{2p} \sum_{z_2, \ldots, z_n} \frac{(m+1) (m_2 +1) \cdots (m_n+1)}{(2g + \ell -3)^{2p-3 +n}} \cdot O(1)^p
				.
			\end{align*}
			
			Using $\sum_{z_j} (m_j +1) = 2g-2 + \ell$, the second term in the sum above yields
			\begin{align*}
				\sum_{n=2}^{2p} \sum_{z_2, \ldots, z_n} \frac{(m+1) (m_2 +1) \cdots (m_n+1)}{(2g + \ell -3)^{2p-3 +n}} \cdot O(1)^p
				& \leq \sum_{n=2}^{2p} \frac{m+1}{(2g + \ell -3)^{2p-2}} \cdot O(1)^p \\
				& = (2p-2) \cdot \frac{m+1}{(2g + \ell -3)^{2p-2}} \cdot O(1)^p
				.
			\end{align*}
			
			Hence, the Siegel--Veech constant is of the claimed order.
			
			Furthermore, if $p \leq \sfrac{m}{2}$ and $m \geq \sfrac{g}{C}$ for some $C > 0$, the dominant configuration for~$n=1$, which provides the first term in the sum, exists, and is of a larger order than the second term in the sum.
		\end{proof}
	\end{thm}
	
	Applying the same argument of summing over all possible zeros once again for the fixed zero, yields the following statement on Siegel--Veech constants for non-fixed zeros.
	
	\begin{cor}[Loops of higher multiplicity at any zero]
		Let $\mathcal{H}$ be a stratum or the odd/even/non-hyperelliptic component of a non-connected stratum.
		Consider the configuration of a loop of multiplicity $p$ at any zero.
		Let $M$ be the largest order of zeros.
		
		For multiplicity $p$, the Siegel--Veech constant is
		\begin{equation*}
			c(\ast \scloop, p, \anyzeros, \mathcal{H})
			=  \frac{M-2p+1}{(2g+\ell-3)^{2p-3}} \cdot O(1)^p
		\end{equation*}
	\end{cor}
	
	With the arguments used here, we do not get exact growth rates for the Siegel--Veech constants for higher multiplicities for all strata, as is the case for multiplicity $1$.
	This is because the dominant configurations do not necessarily exist for every choice of $n$ zeros and hence the sum in the proof above over all zeros $z_2, \ldots, z_n$ depends strongly on the specific strata that we consider.
	
	However, the asymptotics are still sufficient to deduce the statement that the Siegel--Veech constant for any multiplicity is given by the Siegel--Veech constant for multiplicity~$1$.
	
	\begin{cor}[Loops of any multiplicity at a fixed zero]
		\label{cor:loop_any_mult}
		Let $\mathcal{H}$ be a stratum or the odd/even/non-hyperelliptic component of a non-connected stratum.
		Consider the configuration of a loop of any multiplicity at a fixed zero of order~$m$.
		Then the Siegel--Veech constant is
		\begin{equation*}
			c(m \scloop, p = \any, \onefixedzero, \mathcal{H})
			=  \frac{(m +1)^2}{2}
			\cdot \left( 1 + O \left( \frac1g \right) \right)
			.
		\end{equation*}
	\end{cor}
	
	Furthermore, if we fix a specific type of strata, we can refine the previous arguments to get exact growth rates. 
	For example, for strata with $\ell$ zeros of order $\sfrac{2g-2}{\ell}$ such as~$\mathcal{H}(2g-2)$ or $\mathcal{H}(g-1,g-1)$, we can not only use the second estimate in \cref{thm:loops_one_fixed_zero} but we also see from the proof of \cref{prop:loops_all_zeros_fixed} that the term of largest order comes from dominant configurations with only one zero involved, hence we can multiply the Siegel--Veech constant for a fixed zero with the number $\ell$ of zeros to obtain the Siegel--Veech constant for any zero as
	\begin{align*}
		& c(\sfrac{2g-2}{\ell} \scloop, p \leq \sfrac{g-1}{\ell}, \anyzeros, \mathcal{H}) \\
		& = \ell \cdot \frac12 \cdot \left( \frac{\pi^2}{6} \right)^{p-1} \cdot \frac{(\sfrac{2g-2}{\ell} +1)(\sfrac{2g-2}{\ell}-2p+1)}{(2g+\ell-3)^{2p-2}}
		\cdot \left( 1 + O \left( \frac1g \right) \cdot O(1)^p \right) \\
		& = \frac{1}{2 \ell} \cdot \left( \frac{\pi^2}{6} \right)^{p-1} \cdot \frac{1}{(2g + \ell -3)^{2p-4}}
		\cdot \left( 1 + O \left( \frac1g \right) \cdot O(1)^p \right)
		.
	\end{align*}
	
	For the principal stratum, neither the second estimate in \cref{thm:loops_one_fixed_zero} can be applied nor can we exclude this way that cylinders are counted multiple times.
	Hence, we carry out the arguments for the principal stratum separately.
	
	\begin{prop}[Loops of higher multiplicity at any zero in $\mathcal{H}(1,\ldots,1)$]
		\label{prop:mult_higher_loops_unlabelled_principal}
		Consider the principal stratum $\mathcal{H}(1,\ldots,1)$ and the configuration of a saddle connection from any zero to itself.
		
		The Siegel--Veech constant for multiplicity $p$ is
		\begin{equation*}
			c(1 \scloop, p, \anyzeros, \mathcal{H}(1,\ldots,1))
			= \frac{1}{2} \cdot \left(\frac{\pi^2}{3}\right)^{p-1} \cdot \frac{1}{(4g-5)^{2p-3}} \cdot \left( 1 +  O\left( \frac{1}{g} \right) \cdot O(1)^p \right)
			.
		\end{equation*}
		
		\begin{proof}
			To obtain a loop at a zero of order $1$,
			we only can glue a surface from a two-hole construction (with $b_i' = 0$ or $b_i'' = 0$) to a cylinder.
			Hence, to obtain a translation surface in the principal stratum $\mathcal{H}(1,\ldots,1)$ of genus~$g$ with $p$ homologous loops, we glue $p$ cylinders and $p$ subsurfaces obtained by two-hole constructions from principal strata $\mathcal{H}_i$ with $b_i' = b_i'' = 0$ (see \cite[Section 13.5]{eskin_masur_zorich_03}).
			
			From the proof of \cref{prop:loops_all_zeros_fixed}, we see that it is enough to consider dominant configurations. This means that one of the $\mathcal{H}_i$ (say, $\mathcal{H}_1$) has to be the principal stratum with $2g-2 -2p$ zeros of order $1$ and two marked points and hence
			\begin{equation*}
				\sfrac{d_1}{2}-1 = 2(g-p) + 2g-2p -2 = 4g-4p-2
			\end{equation*}
			whereas all other $\sfrac{d_i}{2} - 1
			= 2$. Furthermore, we have $|\Gamma| = |\Gamma_-| = 1$.
			
			By \cref{prop:loops_non-connected_strata}, the Siegel--Veech constant for such a configuration is
			\begin{align*}
				& c(1,\ldots, 1 \scloop, p, \labelled, \dominant, \mathcal{H}(1,\ldots,1)) \\
				& = \frac{ 2 \cdot \ldots \cdot 2}{|\Gamma| \cdot |\Gamma_-|}
				\cdot \left( \frac{\pi^2}{6} \right)^{p-1}
				\cdot \frac{\prod_{i=1}^p (\sfrac{d_i}{2} - 1)!}{(\sfrac{d}{2} - 2)!}
				\cdot \left( 1 + O \left( \frac1g \right) \right) \\
				& = \frac{2^{2p}}{1}
				\cdot \left( \frac{\pi^2}{6} \right)^{p-1}
				\cdot \frac{ (4g-4p -2)! \cdot 2^{p-1}}{(4g-5)!}
				\cdot \left(1 + O \left( \frac1g \right)\right) \\
				& = 2^{2p}
				\cdot \left( \frac{\pi^2}{3} \right)^{p-1}
				\cdot \frac{1}{(4g-5)^{4p-3}}
				\cdot \left(1 + O \left( \frac1g \right) \cdot O(1)^p \right)
				.
			\end{align*}
			
			To obtain the Siegel--Veech constant for the unlabelled situation, we have to take into account the $(2g-2)(2g-3)\ldots(2g-2p-1) = (2g-2)^{2p} \cdot (1+O(\sfrac1g))$ choices for the~$2p$ zeros. As always two choices yield the same configuration,
			we have
			\begin{align*}
				& c(1 \scloop, p, \anyzeros, \mathcal{H}(1,\ldots, 1)) \\
				& = \frac{(2g-2)^{2p}}{2} \cdot 2^{2p}
				\cdot \left( \frac{\pi^2}{3} \right)^{p-1}
				\cdot \frac{1}{(4g-5)^{4p-3}}
				\cdot \left(1 + O \left( \frac1g \right) \cdot O(1)^p \right) \\
				& = \frac{1}{2}
				\cdot \left( \frac{\pi^2}{3} \right)^{p-1}
				\cdot \frac{1}{(4g-5)^{2p-3}}
				\cdot \left(1 + O \left( \frac1g \right) \cdot O(1)^p \right)
				.
				\qedhere
			\end{align*}
		\end{proof}
	\end{prop}

	\section[Siegel--Veech constants for the hyperelliptic component of H(2g-2)]{Siegel--Veech constants for $\mathcal{H}^\hyp(2g-2)$} \label{sec:values_minimal}
	
	In the \emph{minimal stratum} $\mathcal{H}(2g-2)$, there only exist loops, no saddle connections between distinct zeros. Furthermore, there is no difference between the situation of fixed and non-fixed zeros.
	However, we do have to consider the three connected components $\mathcal{H}^\hyp (2g-2)$, $\mathcal{H}^\odd (2g-2)$, and $\mathcal{H}^\even (2g-2)$ separately.
	
	In the hyperelliptic component, the additional symmetry leads to the rare phenomenon that loops of multiplicity $1$ are \emph{not} more common than loops of multiplicity $2$.
	
	\begin{prop}[{\textls[-5]{Loops in the hyperelliptic component of the minimal stratum}}]
		\label{prop:loops_hyperelliptic_minimal}
		Consider the hyperelliptic component $\mathcal{H}^\hyp (2g-2)$ and the configuration of a loop. On a generic surface in $\mathcal{H}^\hyp (2g-2)$, all loops have multiplicity~$1$~or~$2$.
		
		The Siegel--Veech constant for multiplicity $1$ is
		\begin{align*}
			c( \twogminustwo \scloop, p=1, \mathcal{H}^\hyp (2g-2))
			& = \left(\frac{2}{\pi} + \frac{2}{\pi^2} \right) \cdot g^2 \cdot \left( 1 + O \left( \frac1g \right)\right)
			\intertext{and the Siegel--Veech constant for multiplicity $2$ is}
			c( \twogminustwo \scloop, p=2, \mathcal{H}^\hyp (2g-2))
			& = \left( \frac32 - \frac{2}{\pi} - \frac{2}{\pi^2} \right) \cdot g^2 \cdot \left( 1 + O \left( \frac{1}{g^{\sfrac14}} \right) \right)
			.
		\end{align*}
		
		\begin{proof}
			The formulas for this situation are collected in \cref{thm:emz_hyperelliptic_twogminustwo_loops}, as well as the statement that loops of multiplicity greater than $2$ do not exist on generic surfaces in~$\mathcal{H}^\hyp (2g-2)$.
			
			Adding up the first two formulas from \cref{thm:emz_hyperelliptic_twogminustwo_loops} which correspond to a loop of multiplicity $1$, obtained by either a two-hole construction or by a figure-eight construction and a cylinder,~yields
			\begin{equation*}
				c( \twogminustwo \scloop, p=1, \mathcal{H}^\hyp(2g-2)) \\
				= \frac{g -1}{2} \cdot \frac{\mu(\mathcal{H}^\hyp (g \! - \! 2,g \! - \! 2))}{\mu(\mathcal{H}^\hyp (2g \! - \! 2))}
				+ \frac{2g-3}{2(2g-2)} \cdot \frac{\mu(\mathcal{H}^\hyp (2g \! - \! 4))}{\mu(\mathcal{H}^\hyp (2g \! - \! 2))}
				.
			\end{equation*}
			From \cref{eq:volume_hyperelliptic_minimal,eq:volume_hyperelliptic_two} and $\frac{(2g-2)!!}{(2g-3)!!} = \sqrt{\pi (g-1)} \cdot (1 + O (\sfrac1g))$, we obtain
			\begin{align*}
				\frac{\mu(\mathcal{H}^\hyp (g-2,g-2))}{\mu(\mathcal{H}^\hyp (2g-2))}
				& = \frac{4}{\pi^2} \cdot (2g+1) \cdot \frac{(2g-4)!!}{(2g-3)!!} \cdot \frac{(2g-2)!!}{(2g-3)!!} \\
				& = \frac{4}{\pi^2} \cdot \frac{2g+1}{2g-2} \cdot \pi (g-1) \cdot \left( 1 + O \left( \frac1g \right) \right) \\
				& = \frac{4}{\pi} \cdot g \cdot \left( 1 + O \left(\frac1g \right) \right) , \\
				\frac{\mu(\mathcal{H}^\hyp (2g-4))}{\mu(\mathcal{H}^\hyp (2g-2))}
				& = \frac{1}{\pi^2} \cdot (2g+1) \cdot 2g \cdot \frac{(2g-5)!!}{(2g-4)!!} \cdot \frac{(2g-2)!!}{(2g-3)!!} \\
				& = \frac{4}{\pi^2} \cdot g^2 \cdot \frac{2g-2}{2g-3} \cdot \left( 1 + O \left( \frac1g \right) \right) \\
				& = \frac{4}{\pi^2} \cdot g^2 \cdot \left( 1 + O \left( \frac1g \right) \right)
			\end{align*}
			and hence
			\begin{align*}
				c( \twogminustwo \scloop, p=1, \mathcal{H}^\hyp(2g-2))
				& = \left( \frac{g-1}{2} \cdot \frac{4}{\pi} \cdot g
				+ \frac{1}{2} \cdot \frac{4}{\pi^2} \cdot g^2  \right)
				\cdot  \left( 1 + O \left( \frac1g \right)\right) \\
				& = \left( \frac{2}{\pi} + \frac{2}{\pi^2} \right) \cdot g^2 \cdot \left( 1 + O \left( \frac1g \right)\right)
				.
			\end{align*}
			
			The last formula in \cref{thm:emz_hyperelliptic_twogminustwo_loops} is for loops of multiplicity $2$,
			obtained from a figure-eight construction on surfaces from~$\mathcal{H}^\hyp(2g_1 -2)$ and a two-hole construction on surfaces from~$\mathcal{H}^\hyp(g_2-1,g_2-1)$. The Siegel--Veech constant for a fixed partition of $g-1 = g_1 + g_2$~is
			\begin{equation*}
				\frac{(2g_1 -1) \cdot 2g_2}{8} \cdot \frac{(2g_1-1)! \cdot (2g_2)!}{(2g-2)!} 
				\cdot \frac{\mu(\mathcal{H}^\hyp (2g_1-2)) 
					\cdot \mu(\mathcal{H}^\hyp (g_2-1, g_2-1))}{ \mu(\mathcal{H}^\hyp  (2g-2))}
				.
			\end{equation*}
			
			With \cref{eq:volume_hyperelliptic_minimal,eq:volume_hyperelliptic_two} and $\frac{(2g)!!}{(2g-1)!!} = \sqrt{\pi g} \cdot (1 + O (\sfrac1g))$, this is equal to
			\begin{align*}
				& \frac{(2g_1 \! -\! 1) \cdot 2g_2}{8} 
				\! \cdot \! \frac{(2g_1 \! - \! 1)!(2g_2)!}{(2g \! - \! 2)!} 
				\! \cdot \! \frac{8}{\pi^2}
				\! \cdot \! \frac{(2g \! + \! 1)!}{(2g_1 \! + \! 1)! \cdot (2g_2 \! + \! 2)!}
				\! \cdot \! \frac{(2g_1 \! - \! 3)!!}{(2g_1 \! - \! 2)!!}
				\! \cdot \! \frac{(2g_2 \! - \! 2)!!}{(2g_2 \! - \! 1)!!}
				\! \cdot \! \frac{(2g \! - \! 2)!!}{(2g \! - \! 3)!!}
				\\
				& = \frac{1}{\pi^2} \cdot \frac{(2g)^3}{(2g_1+1) \cdot (2g_2+2)}
				\cdot \frac{(2g_1-1)!!}{(2g_1)!!}
				\cdot \frac{(2g_2)!!}{(2g_2+1)!!}
				\cdot \frac{(2g)!!}{(2g-1)!!}
				\cdot \left( 1 + O \left( \frac1g \right) \right) \\
				& = \frac{8}{\pi^2} \cdot \frac{g^3}{(2g_1+1) \cdot (2g_2+2)}
				\cdot \frac{(2g_1-1)!!}{(2g_1)!!}
				\cdot \frac{(2g_2)!!}{(2g_2+1)!!}
				\cdot \sqrt{\pi g}
				\cdot \left( 1 + O \left( \frac1g \right) \right)
				.
			\end{align*}
			
			We sum over all partitions of $g_1 + g_2 = g-1$ and with the third estimate from \cref{lem:sum_with_double_factorials}, this yields
			\begin{align*}
				& c( \twogminustwo \scloop, p=2, \mathcal{H}^\hyp(2g-2)) \\
				& = \frac{8 \sqrt{\pi}}{\pi^2} \cdot g^2 \cdot \sum_{g_1 + g_2 = g-1} \frac{g^{\sfrac32}}{(2g_1+1) \cdot (2g_2+2)}
				\cdot \frac{(2g_1-1)!!}{(2g_1)!!} 
				\cdot \frac{(2g_2)!!}{(2g_2+1)!!}
				\cdot \left( 1+ O \left( \frac1g \right) \right)  \\
				& =  \frac{8 \sqrt{\pi}}{\pi^2} \cdot g^2
				\cdot \frac{\sqrt{\pi}}{4} \cdot \left( \frac{3 \pi}{4} - 1 - \frac{1}{\pi} \right)
				\cdot \left( 1 + O \left( \frac{1}{g^{\sfrac14}} \right) \right) \\
				& = \left( \frac32 - \frac{2}{\pi} -\frac{2}{\pi^2} \right)
				\cdot g^2
				\cdot \left( 1 + O \left( \frac{1}{g^{\sfrac14}} \right) \right)
				. \qedhere
			\end{align*}
		\end{proof}
	\end{prop}
	
	For the odd and the even component as well as for the whole stratum $\mathcal{H}(2g-2)$, 
	we can deduce the Siegel--Veech constants from \cref{thm:loops_one_fixed_zero} as
	\begin{align}
		\label{eq:loops_principal}
		& c(\twogminustwo \scloop, p, \mathcal{H}(2g-2) / \mathcal{H}^\oddeven(2g-2)) \notag \\
		& = \frac12 \cdot \left( \frac{\pi^2}{6} \right)^{p-1} \cdot \frac{1}{(2g - 2)^{2p-4}}
		\cdot \left( 1 + O \left( \frac1g \right) \cdot O(1)^p \right) 
		.
	\end{align}

	\section[Siegel--Veech constants for the hyperelliptic component of H(g-1,g-1)]{Siegel--Veech constants for $\mathcal{H}^\hyp(g-1,g-1)$} \label{sec:values_g-1}
	
	In the case of the stratum $\mathcal{H}(g-1,g-1)$, we again have to consider different components.
	
	We start with saddle connections between distinct zeros in the hyperelliptic component. Again, there is no distinction between the case of fixed and non-fixed zeros as saddle connections are not oriented.
	
	\pagebreak[3]
	
	\begin{prop}[Distinct zeros in the hyperelliptic component of $\mathcal{H}(g-1,g-1)$]
		\label{prop:distinct_hyperelliptic_g-1}
		Consider the hyperelliptic component $\mathcal{H}^\hyp (g-1,g-1)$ and the configuration of a saddle connection between distinct zeros. On a generic surface in $\mathcal{H}^\hyp (g-1,g-1)$, all saddle connections have multiplicity $1$ or $2$.
		The Siegel--Veech constant for multiplicity $1$ is
		\begin{align*}
			c( \gminusone \saddleconnection \gminusone, p=1, \mathcal{H}^\hyp (g-1,g-1))
			& = \frac{2}{\pi} \cdot g^2 \cdot \left( 1 + O \left( \frac1g \right) \right)
			\intertext{and the Siegel--Veech constant for multiplicity $2$ is}
			c( \gminusone \saddleconnection \gminusone, p=2, \mathcal{H}^\hyp (g-1,g-1))
			& = \left( 1 - \frac{2}{\pi} \right) \cdot g^2 \cdot \left( 1 + O \left( \frac{1}{g^{\sfrac14}} \right) \right)
			.
		\end{align*}
		
		\begin{proof}
			The formulas for this situation are collected in \cref{thm:emz_hyperelliptic_distinct}, as well as the statement that saddle connections of multiplicity greater than $2$ do not exist on generic surfaces in~$\mathcal{H}^\hyp (g-1, g-1)$.
			
			For saddle connections of multiplicity $1$, we can directly plug \cref{eq:volume_hyperelliptic_minimal_Stirling,eq:volume_hyperelliptic_two_Stirling} into the first formula in \cref{thm:emz_hyperelliptic_distinct} to obtain
			\begin{align*}
				c( \gminusone \saddleconnection \gminusone, p=1, \mathcal{H}^\hyp (g-1, g-1))
				& = (2g-1) \cdot \frac{\mu(\mathcal{H}^\hyp (2g-2))}{\mu( \mathcal{H}^\hyp (g-1,g-1))} \\
				& = (2g-1) \cdot \frac{g}{\pi} \cdot \left( 1 + O \left( \frac1g \right) \right) \\
				& = \frac{2}{\pi} \cdot g^2 \cdot \left( 1 + O \left( \frac1g \right) \right)
				.
			\end{align*}
			
			For saddle connections of multiplicity $2$, the second formula in \cref{thm:emz_hyperelliptic_distinct} yields for a fixed partition $g= g_1 + g_2$ the corresponding Siegel--Veech constant
			\begin{equation*}
				\frac{(2g_1-1)(2g_2-1)}{2} 
				\cdot \frac{(2g_1-1)! \cdot (2g_2-1)!}{(2g-1)!}
				\cdot \frac{\mu( \mathcal{H}^\hyp (2g_1-2)) \cdot \mu( \mathcal{H}^\hyp (2g_2-2))}{\mu( \mathcal{H}^\hyp (g-1,g-1))}
			\end{equation*}
			except for $g_1 = g_2$ where we have an additional factor of $\sfrac{1}{|\Gamma|} = \sfrac12$.
			
			With \cref{eq:volume_hyperelliptic_minimal,eq:volume_hyperelliptic_two} and $\frac{(2g-1)!!}{(2g)!!} = \sfrac{1}{\sqrt{\pi g}} \cdot (1 + O (\sfrac1g))$, this is equal to
			\begin{align*}
				& \frac{(2g_1 \!  - \! 1) (2g_2 \! - \! 1)}{2} \!  \cdot \!  \frac{(2g_1 \! - \! 1)!(2g_2 \! - \! 1)!}{(2g \! - \! 1)!} 
				\! \cdot \!  \frac12
				\! \cdot \! \frac{(2g \! + \! 2) !}{(2g_1 \! + \! 1) !  (2g_2 \! + \! 1) !}
				\! \cdot \! \frac{( 2g_1 \!\! - \!\! 3 ) !!}{( 2g_1 \!\! - \!\! 2 ) !!}
				\! \cdot \! \frac{( 2g_2 \!\! - \!\! 3 ) !!}{( 2g_2 \!\! - \!\! 2 ) !!}
				\! \cdot \! \frac{( 2g \!\! - \!\! 1 ) !!}{( 2g \!\! - \!\! 2 ) !!}
				\\
				& = \frac{1}{4} 
				\cdot \frac{(2g)^3 \cdot 2g}{(2g_1+1) \cdot (2g_2+1)}
				\cdot \frac{(2g_1-1)!!}{(2g_1)!!}
				\cdot \frac{(2g_2-1)!!}{(2g_2)!!}
				\cdot \frac{(2g-1)!!}{(2g)!!}
				\cdot \left( 1 + O \left( \frac1g \right) \right) \\
				& = 4 \cdot \frac{g^4}{(2g_1+1) \cdot (2g_2+1)}
				\cdot \frac{(2g_1-1)!!}{(2g_1)!!}
				\cdot \frac{(2g_2-1)!!}{(2g_2)!!}		
				\cdot \frac{1}{\sqrt{\pi g}}
				\cdot \left( 1 + O \left( \frac1g \right) \right)
				.
			\end{align*}
			
			Note that always two choices of a partition $g = g_1 + g_2$ yield the same configuration (except for $g_1 = g_2$). Hence when summing over all partitions, we obtain with the first estimate in \cref{lem:sum_with_double_factorials}
			\begin{align*}
				& c( \gminusone \saddleconnection \gminusone, p=2, \mathcal{H}^\hyp(g-1, g-1)) \\
				& = \frac12 \cdot \frac{4}{\sqrt{\pi}} \cdot g^2 \cdot
				\sum_{g_1 + g_2 = g} \frac{g^{\sfrac32}}{(2g_1+1) \cdot (2g_2+1)}
				\cdot \frac{(2g_1-1)!!}{(2g_1)!!}
				\cdot \frac{(2g_2-1)!!}{(2g_2)!!}
				\cdot \left( 1+ O \left( \frac1g \right) \right)
				\\
				& = \frac{2}{\sqrt{\pi}} \cdot g^2
				\cdot \sqrt{\pi} \cdot \left( \frac12 - \frac{1}{\pi} \right)
				\cdot \left( 1+ O \left( \frac{1}{g^{\sfrac14}} \right) \right)
				\\
				& = \left( 1 - \frac{2}{\pi} \right)
				\cdot g^2
				\cdot \left( 1 + O \left( \frac{1}{g^{\sfrac14}}\right) \right)
				. \qedhere
			\end{align*}
		\end{proof}
	\end{prop}
	
	Here, in the case of saddle connections between distinct zeros, we could have also deduced the Siegel--Veech constant for multiplicity $2$ from the one for multiplicity $1$ by \cref{eq:exact_up_to_homology}. 
	This actually yields a better error term of $1+ O(\sfrac1g)$ instead of $1+ O(\sfrac{1}{g^{\sfrac14}})$.

	\begin{prop}[Loops in the hyperelliptic component of $\mathcal{H}(g-1,g-1)$]
		\label{prop:loops_hyperelliptic_g-1}
		Consider the hyperelliptic component $\mathcal{H}^\hyp (g-1,g-1)$ and the configuration of a loop. On a generic surface in $\mathcal{H}^\hyp (g-1,g-1)$, all loops have multiplicity~$1$ or $2$.
		
		The Siegel--Veech constant for multiplicity $1$ is
		\begin{align*}
			c( \gminusone \scloop, p=1, \anyzeros, \mathcal{H}^\hyp (g-1,g-1))
			& = \frac{2}{\pi^2} \cdot g^2 \cdot \left( 1 + O \left( \frac1g \right) \right)
			\intertext{and the Siegel--Veech constant for multiplicity $2$ is}
			c( \gminusone\scloop , p=2, \anyzeros, \mathcal{H}^\hyp (g-1,g-1))
			& = \left(\frac12 - \frac{2}{\pi^2} \right) \cdot g^2 \cdot \left( 1 + O \left( \frac{1}{g^{\sfrac14}} \right) \right)
			.
		\end{align*}
		
		\begin{proof}
			The formulas for this situation are collected in \cref{thm:emz_hyperelliptic_gminusone_gminusone_loops}, as well as the statement that loops of multiplicity greater than $2$ do not exist on generic surfaces in~$\mathcal{H}^\hyp (g-1, g-1)$.
			
			For loops of multiplicity $1$, we can directly plug \cref{eq:volume_hyperelliptic_two} into the first formula in \cref{thm:emz_hyperelliptic_gminusone_gminusone_loops} to obtain
			\begin{align*}
				& c( \gminusone \scloop, p=1, \anyzeros, \mathcal{H}^\hyp (g-1, g-1)) \\
				& = \frac{g-1}{2g-1} \cdot \frac{\mu(\mathcal{H}^\hyp (g-2, g-2))}{\mu( \mathcal{H}^\hyp (g-1,g-1))} \\
				& = \frac{g-1}{2g-1} \cdot \frac{(2g+2)(2g+1)}{\pi^2}
				\cdot \frac{(2g-4)!!}{(2g-3)!!} \cdot \frac{(2g-1)!!}{(2g-2)!!}	\\
				& = \frac{2}{\pi^2} \cdot g^2 \cdot \left( 1 + O \left( \frac1g \right) \right)
				.
			\end{align*}
			
			For loops of multiplicity $2$, the second formula in \cref{thm:emz_hyperelliptic_gminusone_gminusone_loops} yields for a fixed partition $g -1 = g_1 + g_2$ the corresponding Siegel--Veech constant
			\begin{equation*}
				\frac{g_1 \cdot g_2}{2} 
				\cdot \frac{(2g_1)! \cdot (2g_2)!}{(2g-1)!}
				\cdot \frac{\mu( \mathcal{H}^\hyp (g_1-1, g_1-1)) \cdot \mu( \mathcal{H}^\hyp (g_2-1, g_2-1))}{\mu( \mathcal{H}^\hyp (g-1,g-1))}
			\end{equation*}
			except for $g_1 = g_2$ where we have an additional factor of $\sfrac{1}{|\Gamma|} = \sfrac12$.
			
			With \cref{eq:volume_hyperelliptic_two} and $\frac{(2g-1)!!}{(2g)!!} = \sfrac{1}{\sqrt{\pi g}} \cdot (1 + O (\sfrac1g))$, this is equal to
			\begin{align*}
				& \frac{2g_1 \cdot 2g_2}{8} \cdot \frac{(2g_1)!(2g_2)!}{(2g-1)!} 
				\cdot \frac{8}{\pi^2}
				\cdot \frac{(2g+2)!}{(2g_1+2)! \cdot (2g_2+2)!}
				\cdot \frac{(2g_1-2)!!}{(2g_1-1)!!}
				\cdot \frac{(2g_2-2)!!}{(2g_2-1)!!}
				\cdot \frac{(2g-1)!!}{(2g-2)!!}
				\\
				& = \frac{1}{\pi^2} 
				\cdot \frac{(2g)^3 \cdot 2g}{(2g_1 +2) \cdot (2g_2 +2)}
				\cdot \frac{(2g_1)!!}{(2g_1+1)!!}
				\cdot \frac{(2g_2)!!}{(2g_2+1)!!}
				\cdot \frac{(2g-1)!!}{(2g)!!}
				\cdot \left( 1 + O \left( \frac1g \right) \right) \\
				& = \frac{16}{\pi^2} \cdot \frac{g^4}{(2g_1 +2) \cdot (2g_2 +2)}
				\cdot \frac{(2g_1)!!}{(2g_1+1)!!}
				\cdot \frac{(2g_2)!!}{(2g_2+1)!!}
				\cdot \frac{1}{\sqrt{\pi g}}
				\cdot \left( 1 + O \left( \frac1g \right) \right)
				.
			\end{align*}
			
			Note that always two choices of a partition $g -1 = g_1 + g_2$ yield the same configuration (except for $g_1 = g_2$). Hence when summing over all partitions, we obtain with the second estimate from \cref{lem:sum_with_double_factorials}
			\begin{align*}
				& c( \gminusone \scloop, p=2, \anyzeros, \mathcal{H}^\hyp(g-1, g-1)) \\
				& = \frac12 \cdot \frac{16}{\pi^2 \cdot \sqrt{\pi}} \cdot g^2 \cdot
				\sum_{g_1 + g_2 = g-1} \frac{g^{\sfrac32}}{(2g_1+2) \cdot (2g_2+2)}
				\cdot \frac{(2g_1)!!}{(2g_1+1)!!}
				\cdot \frac{(2g_2)!!}{(2g_2+1)!!}
				\cdot \left( \! 1 \! + \! O \! \left( \! \frac1g \! \right) \! \right)
				\\
				& = \frac{8}{\pi^2 \cdot \sqrt{\pi}} \cdot g^2
				\cdot \frac{\sqrt{\pi}}{4} \cdot \left( \frac{\pi^2}{4} - 1 \right)
				\cdot \left( 1+ O \left( \frac{1}{g^{\sfrac14}} \right) \right)
				\\
				& = \left( \frac12 - \frac{2}{\pi^2} \right)
				\cdot g^2
				\cdot \left( 1 + O \left( \frac{1}{g^{\sfrac14}}\right) \right)
				. \qedhere
			\end{align*}
		\end{proof}
	\end{prop}
	
	From \cite[Lemma 14.2]{eskin_masur_zorich_03}, we also see that loops of multiplicity $1$ in $\mathcal{H}^\hyp(g-1,g-1)$ come from configurations with one two-hole construction and a cylinder whereas loops of multiplicity~$2$ come from configurations with two two-hole constructions and no cylinders. 
	In particular, in both cases, 
	any loop at a fixed zero is homologous to a loop at the other zero. Hence,
	the Siegel--Veech constant for loops at a fixed zero is the same as for loops at any zero.
	
	\bigskip
	
	For the odd and the even or the non-hyperelliptic component as well as for the whole stratum $\mathcal{H}(g-1,g-1)$, 
	we can deduce the Siegel--Veech constants for saddle connections from \cref{thm:distinct_fixed_zeros} as
	\begin{align}
		\label{eq:distinct_zeros_H_gminusone_gminusone}
		& c(\gminusone \saddleconnection \gminusone, p \leq g, \mathcal{H}(g-1, g-1) / \mathcal{H}^\oen(g-1, g-1)) \notag \\
		& = \frac14 \cdot \left( \frac{\pi^2}{6} \right)^{p-1} \cdot \frac{1}{(2g-1)^{2p-4}}
		\cdot \left( 1 + O \left( \frac1g \right) \cdot O(1)^p \right) 
		\intertext{and for loops at a fixed zero from \cref{thm:loops_one_fixed_zero} as}
		\label{eq:loops_H_gminusone_gminusone_fixed_zero}
		& c(\gminusone \scloop, p, \onefixedzero, \mathcal{H}(g-1, g-1) / \mathcal{H}^\oen(g-1, g-1)) \notag \\
		& = \frac18 \cdot \left( \frac{\pi^2}{6} \right)^{p-1} \cdot \frac{1}{(2g-1)^{2p-4}}
		\cdot \left( 1 + O \left( \frac1g \right) \cdot O(1)^p \right) 
		.
		\intertext{
			Recall from the discussion after \cref{cor:loop_any_mult} that for $\mathcal{H}(g-1,g-1)$, the dominant configurations for loops 
			involve only homology classes of loops at one fixed zero,
			hence the Siegel--Veech constant for any zero is obtained by multiplying with the number of zeros as}
		\label{eq:loops_H_gminusone_gminusone_any_zero}
		& c(\gminusone \scloop, p, \anyzero, \mathcal{H}(g-1, g-1) / \mathcal{H}^\oen(g-1, g-1)) \notag \\
		& = \frac14 \cdot \left( \frac{\pi^2}{6} \right)^{p-1} \cdot \frac{1}{(2g-1)^{2p-4}}
		\cdot \left( 1 + O \left( \frac1g \right) \cdot O(1)^p \right) 
		.
	\end{align}

	\section{Table of values} \label{sec:table}
	
	We summarize the saddle connection Siegel--Veech constants that we have recalled and determined in \cref{tab:SV_constants}. If there are two options listed for a property, it means that the asymptotics are the same for both options.
	
	Having available all saddle-connection Siegel--Veech constants now, we do a last calculation here. The Siegel--Veech constant for any stratum or any odd/even/non-hyperelliptic component and any saddle connection (whether between distinct zeros or loops) can be determined as the sum of all Siegel--Veech constants for (non-oriented!) saddle connections between distinct zeros $z_1$ and $z_2$ of orders $m_1$ and $m_2$ and for loops at a zero $z$ of order~$m$.
	Hence we obtain
	\begin{align}
		\label{eq:SV_all}
		& c(\textnormal{any saddle connection}, \mathcal{H}) \notag \\
		& = \frac12 \sum_{z_1 \neq z_2} (m_1+1)(m_2+1) 
		\cdot \left( 1 + O \left( \frac1g \right) \right)
		+ \sum_{z} \frac{(m+1)^2}{2}
		\cdot \left( 1 + O \left( \frac1g \right) \right) \notag \\
		& = \frac12 \sum_{z_1, z_2} (m_1+1)(m_2+1)
		\cdot \left( 1 + O \left( \frac1g \right) \right) \notag \\
		& = \frac{(2g+\ell-2)^2}{2}
		\cdot \left( 1 + O \left( \frac1g \right) \right)
		.
	\end{align}

		\begin{landscape}
			
			\thispagestyle{empty}
			
			\begin{table} \label{tab:SV_constants}
				\vspace*{-1.5cm}
				\hspace*{-4.2cm}
				\small
				\begin{tabular}{r|l|l|l|l|l|l}
					stratum & multiplicity & loops? & order of zeros & fixed/any zeros & large-genus asymptotics & details in \\ \hline &&&&&& \\[-0.5em]
					\ldelim\{{6.8}{*}[\parbox{3.9cm}{any stratum $\mathcal{H}$ / \\ connected component \\ $\mathcal{H}^\odd$ / $\mathcal{H}^\even$ / $\mathcal{H}^\nonhyp$}]
					& $1$ / any & distinct zeros & order $m_1$, $m_2$ & fixed zeros & $(m_1+1)(m_2+1) \cdot (1+O(\frac{1}{g}))$ & \cref{eq:saddle_connections_fixed_all_multiplicities}  \\
					& $p\leq \min \{m_1, m_2\} +1$ & distinct zeros & order $m_1$, $m_2$ & fixed zeros & $(\frac{\pi^2}{6})^{p-1} \cdot \frac{(m_1+1)(m_2+1)}{(2g + \ell -2p)^{2p-2}} \cdot (1+O(\frac1g)\cdot O(1)^p)$ & \cref{thm:distinct_fixed_zeros}  \\
					& $\geq \min \{m_1, m_2\} +2$ & distinct zeros & order $m_1$, $m_2$ & fixed zeros & $0$ & \cref{thm:distinct_fixed_zeros}  \\
					& $1$ / any & loops & order $m$ & fixed zeros & $\frac{(m+1)^2}{2} \cdot (1+O(\frac{1}{g}))$ & Eq.~\eqref{eq:mult_one_loop_labelled} / Cor.~\ref{cor:loop_any_mult} \\
					& $p \geq 1$ & loops & order $m$ & fixed zeros & $\frac{(m+1)(m-2p+1)}{(2g+\ell-3)^{2p-2}} \cdot O(1)^p$ & \cref{thm:loops_one_fixed_zero}  \\
					& $1$ / any & any & & any zeros & $\frac{(2g+\ell-2)^2}{2} \cdot (1+O(\frac1g))$ & \cref{eq:SV_all} \\[1em]
					
					\ldelim\{{4.5}{*}[$\mathcal{H}(1,\ldots,1)$] 
					& $1$ / any & distinct zeros & & any zeros & $8g^2 \cdot (1+O(\frac1g))$ & \cref{eq:mult_any_distinct_unlabelled_zeros_principal} \\
					& $2$ & distinct zeros & & any zeros & $\frac{\pi^2}{12} \cdot (1 +O(\frac1g))$ & \cref{eq:mult_two_distinct_unlabelled_zeros_principal} \\
					& $1$ / any & loops & & any zeros & $2g \cdot (1+ O(\frac1g))$ & \cref{eq:loop_mult_1_simple_zero} \\
					& $p \geq 1$ & loops & & any zeros & $\frac{1}{2} \cdot (\frac{\pi^2}{3})^{p-1} \cdot \frac{1}{(4g-5)^{2p-3}} \cdot ( 1 +  O( \frac{1}{g}) \cdot O(1)^p )$ & \cref{prop:mult_higher_loops_unlabelled_principal} \\[1em]
					
					\ldelim\{{5.5}{*}[\parbox{3.74cm}{$\mathcal{H}(g-1, g-1)$ /  \\ $\mathcal{H}^\odd(g-1,g-1)$ / \\ $\mathcal{H}^\even(g-1, g-1)$ / \\ $\mathcal{H}^\nonhyp(g-1, g-1)$ }] 
					& $1$ / any & distinct zeros & & fixed/any zeros & $g^2 \cdot (1+O(\frac1g))$ & \cref{eq:distinct_zeros_H_gminusone_gminusone} \\
					& $p \leq g$ & distinct zeros & & fixed/any zeros & $\frac14 \cdot ( \frac{\pi^2}{6} )^{p-1} \cdot \frac{1}{(2g-1)^{2p-4}}
					\cdot (1+O(\frac1g) \cdot O(1)^p)$ & \cref{eq:distinct_zeros_H_gminusone_gminusone} \\
					& $1$ / any & loops & & fixed zeros & $\frac12 \cdot g^2 \cdot (1+O(\frac1g))$ & \cref{eq:loops_H_gminusone_gminusone_fixed_zero} \\
					& $p \geq 1$ & loops & & fixed zeros & $\frac18 \cdot ( \frac{\pi^2}{6} )^{p-1} \cdot \frac{1}{(2g-1)^{2p-4}}
					\cdot (1+O(\frac1g) \cdot O(1)^p)$ & \cref{eq:loops_H_gminusone_gminusone_fixed_zero} \\
					& $p \geq 1$ & loops & & any zero & $\frac14 \cdot ( \frac{\pi^2}{6} )^{p-1} \cdot \frac{1}{(2g-1)^{2p-4}}
					\cdot (1+O(\frac1g) \cdot O(1)^p)$ & \cref{eq:loops_H_gminusone_gminusone_any_zero} \\[1em]
					
					\ldelim\{{6.3}{*}[$\mathcal{H}^\hyp (g-1,g-1)$] 
					& $1$ & distinct zeros & & fixed/any zeros & $\frac{2}{\pi} \cdot g^2 \cdot (1 + O(\frac1g) )$ & \cref{prop:distinct_hyperelliptic_g-1} \\
					& $2$ & distinct zeros & & fixed/any zeros & $(1 - \frac{2}{\pi}) \cdot g^2 \cdot ( 1 + O( \frac{1}{g^{\sfrac14}} ))$ & \cref{prop:distinct_hyperelliptic_g-1} \\
					& $\geq 3$ & distinct zeros & & fixed/any zeros & $0$ & \cref{thm:emz_hyperelliptic_distinct} \\
					& $1$ & loops & & fixed/any zeros & $\frac{2}{\pi^2} \cdot g^2 \cdot (1 + O(\frac1g) )$ & \cref{prop:loops_hyperelliptic_g-1} \\
					& $2$ & loops & & fixed/any zeros & $(\frac12 - \frac{2}{\pi^2}) \cdot g^2 \cdot ( 1 + O( \frac{1}{g^{\sfrac14}} ))$ & \cref{prop:loops_hyperelliptic_g-1} \\
					& $\geq 3$ & loops & & fixed/any zeros & $0$ & \cref{thm:emz_hyperelliptic_gminusone_gminusone_loops} \\[1em]
					
					\ldelim\{{2.1}{*}[\parbox{2.97cm}{$\mathcal{H}(2g-2)$ / \\ $\mathcal{H}^\oddeven(2g-2)$}] & $1$ / any & loops &  & fixed/any zeros & $2g^2 \cdot (1 + O(\frac1g))$  & \cref{eq:loops_principal} \\
					& $p \geq 1$ & loops & & fixed/any zeros & $\frac12 \cdot (\frac{\pi^2}{6})^{p-1} \cdot \frac{1}{(2g-2)^{2p-4}} \cdot (1 + O(\frac{1}{g}) \cdot O(1)^p)$ & \cref{eq:loops_principal} \\[1em]
					
					\ldelim\{{3.3}{*}[$\mathcal{H}^\hyp (2g-2)$] 
					& $1$ & loops & & fixed/any zeros & $(\frac{2}{\pi} + \frac{2}{\pi^2}) \cdot g^2 \cdot (1 + O(\frac1g) )$ & \cref{prop:loops_hyperelliptic_minimal} \\
					& $2$ & loops & & fixed/any zeros & $( \frac32 - \frac{2}{\pi} - \frac{2}{\pi^2} ) \cdot g^2 \cdot ( 1 + O( \frac{1}{g^{\sfrac14}} ))$ & \cref{prop:loops_hyperelliptic_minimal} \\
					& $\geq 3$ & loops & & fixed/any zeros & $0$ & \cref{thm:emz_hyperelliptic_twogminustwo_loops} \\
				\end{tabular}
				\caption{Lookup table for Siegel--Veech constants.}
			\end{table}
			
		\end{landscape}
		
		\normalsize	
	
	\bibliographystyle{amsalpha}
	\bibliography{/home/office/Documents/Mathematik/literature/BibTex/Literatur}

\newcommand{\etalchar}[1]{$^{#1}$}
\providecommand{\bysame}{\leavevmode\hbox to3em{\hrulefill}\thinspace}
\providecommand{\MR}{\relax\ifhmode\unskip\space\fi MR }
\providecommand{\MRhref}[2]{%
  \href{http://www.ams.org/mathscinet-getitem?mr=#1}{#2}
}
\providecommand{\href}[2]{#2}
\begin{thebibliography}{ACM{\etalchar{+}}24}

\bibitem[ACM{\etalchar{+}}24]{aulicino_calderon_matheus_salter_schmoll_24}
David Aulicino, Aaron Calderon, Carlos Matheus, Nick Salter, and Martin
  Schmoll, \emph{{Siegel--Veech Constants for Cyclic Covers of Generic
  Translation Surfaces}}, arXiv:2409.06600, 2024.

\bibitem[ADG{\etalchar{+}}20]{aggarwal_delecroix_goujard_zograf_zorich_20}
Amol Aggarwal, Vincent Delecroix, {\'E}lise Goujard, Peter Zograf, and Anton
  Zorich, \emph{Conjectural large genus asymptotics of {Masur}-{Veech} volumes
  and of area {Siegel}-{Veech} constants of strata of quadratic differentials},
  Arnold Mathematical Journal \textbf{6} (2020), no.~2, 149--161.

\bibitem[AEZ16]{athreya_eskin_zorich_16}
Jayadev~S. Athreya, Alex Eskin, and Anton Zorich, \emph{Right-angled billiards
  and volumes of moduli spaces of quadratic differentials on
  {{\(\mathbb{C}\mathrm{P}^1\)}}}, Annales Scientifiques de l'{\'E}cole Normale
  Sup{\'e}rieure. Quatri{\`e}me S{\'e}rie \textbf{49} (2016), no.~6,
  1311--1386.

\bibitem[Agg19]{aggarwal_19}
Amol Aggarwal, \emph{{Large genus asymptotics for Siegel--Veech constants}},
  Geometric and Functional Analysis \textbf{29} (2019), no.~5, 1295--1324.

\bibitem[Agg20]{aggarwal_20}
\bysame, \emph{{Large Genus Asymptotics for Volumes of Strata of Abelian
  Differentials}}, {Journal of the American Mathematical Society} \textbf{33}
  (2020), no.~4, 941--989, {With an appendix by Anton Zorich}.

\bibitem[Agg21]{aggarwal_21}
\bysame, \emph{Large genus asymptotics for intersection numbers and principal
  strata volumes of quadratic differentials}, Inventiones Mathematicae
  \textbf{226} (2021), no.~3, 897--1010.

\bibitem[CMS23]{chen_moeller_sauvaget_23}
Dawei Chen, Martin Möller, and Adrien Sauvaget, \emph{Masur-{Veech} volumes
  and intersection theory: the principal strata of quadratic differentials},
  Duke Mathematical Journal \textbf{172} (2023), no.~9, 1735--1779, {With} an
  appendix by {Ga{\"e}tan} {Borot}, {Alessandro} {Giacchetto} and {Danilo}
  {Lewanski}.

\bibitem[CMSZ20]{chen_moeller_sauvaget_zagier_20}
Dawei Chen, Martin Möller, Adrien Sauvaget, and Don Zagier,
  \emph{Masur-{Veech} volumes and intersection theory on moduli spaces of
  abelian differentials}, Inventiones Mathematicae \textbf{222} (2020), no.~1,
  283--373.

\bibitem[CMZ18]{chen_moeller_zagier_18}
Dawei {Chen}, Martin {Möller}, and Don {Zagier}, \emph{{Quasimodularity and
  large genus limits of Siegel-Veech constants}}, {Journal of the American
  Mathematical Society} \textbf{31} (2018), no.~4, 1059--1163.

\bibitem[CMZ24]{costantini_moeller_zachhuber_24}
Matteo Costantini, Martin M{\"o}ller, and Jonathan Zachhuber, \emph{The area is
  a good enough metric}, Annales de l'Institut Fourier \textbf{74} (2024),
  no.~3, 1017--1059.

\bibitem[DGY25]{duryev_goujard_yakovlev_25}
Eduard Duryev, Elise Goujard, and Ivan Yakovlev, \emph{Volumes of odd strata of
  quadratic differentials}, 2025, arXiv:2502.13121, pp.~1--49.

\bibitem[DGZZ21]{delecroix_goujard_zograf_zorich_21}
Vincent Delecroix, {\'E}lise Goujard, Peter Zograf, and Anton Zorich,
  \emph{Masur-{Veech} volumes, frequencies of simple closed geodesics, and
  intersection numbers of moduli spaces of curves}, Duke Mathematical Journal
  \textbf{170} (2021), no.~12, 2633--2718.

\bibitem[EKZ14]{eskin_kontsevich_zorich_14}
Alex Eskin, Maxim Kontsevich, and Anton Zorich, \emph{Sum of {Lyapunov}
  exponents of the {Hodge} bundle with respect to the {Teichm{\"u}ller}
  geodesic flow}, Publications Math{\'e}matiques \textbf{120} (2014), 207--333.

\bibitem[EM01]{eskin_masur_01}
Alex Eskin and Howard Masur, \emph{Asymptotic formulas on flat surfaces},
  Ergodic Theory and Dynamical Systems \textbf{21} (2001), no.~2, 443--478.

\bibitem[EMZ03]{eskin_masur_zorich_03}
Alex {Eskin}, Howard {Masur}, and Anton {Zorich}, \emph{{Moduli spaces of
  Abelian differentials: the principal boundary, counting problems, and the
  Siegel-Veech constants}}, {Publications Math\'ematiques, Institut des Hautes
  \'Etudes Scientifiques} \textbf{97} (2003), 61--179.

\bibitem[EO01]{eskin_okounkov_01}
Alex Eskin and Andrei Okounkov, \emph{Asymptotics of numbers of branched
  coverings of a torus and volumes of moduli spaces of holomorphic
  differentials}, Inventiones Mathematicae \textbf{145} (2001), no.~1, 59--103.

\bibitem[EZ15]{eskin_zorich_15}
Alex Eskin and Anton Zorich, \emph{Volumes of strata of abelian differentials
  and {Siegel}--{Veech} constants in large genera}, Arnold Mathematical Journal
  \textbf{1} (2015), no.~4, 481--488.

\bibitem[Hub59]{huber_59}
Heinz Huber, \emph{Zur analytischen {Theorie} hyperbolischer {Raumformen} und
  {Bewegungsgruppen}}, Mathematische Annalen \textbf{138} (1959), 1--26.

\bibitem[Hub61]{huber_61}
\bysame, \emph{Zur analytischen {Theorie} hyperbolischer {Raumformen} und
  {Bewegungsgruppen}. {II}}, Mathematische Annalen \textbf{142} (1961),
  385--398.

\bibitem[KZ03]{kontsevich_zorich_03}
Maxim Kontsevich and Anton Zorich, \emph{Connected components of the moduli
  spaces of abelian differentials with prescribed singularities}, Inventiones
  mathematicae \textbf{153} (2003), no.~3, 631--678.

\bibitem[Mas82]{masur_82}
Howard Masur, \emph{Interval exchange transformations and measured foliations},
  Annals of Mathematics. Second Series \textbf{115} (1982), 169--200.

\bibitem[Mas88]{masur_88}
\bysame, \emph{Lower bounds for the number of saddle connections and closed
  trajectories of a quadratic differential}, Holomorphic Functions and Moduli I
  (D.~Drasin, I.~Kra, C.~J. Earle, A.~Marden, and F.~W. Gehring, eds.),
  Springer New York, 1988, pp.~215--228.

\bibitem[Mas90]{masur_90}
\bysame, \emph{The growth rate of trajectories of a quadratic differential},
  Ergodic Theory and Dynamical Systems \textbf{10} (1990), no.~1, 151--176.

\bibitem[Mir08]{mirzakhani_08}
Maryam Mirzakhani, \emph{Growth of the number of simple closed geodesics on
  hyperbolic surfaces}, Annals of Mathematics. Second Series \textbf{168}
  (2008), no.~1, 97--125.

\bibitem[Rei24]{reichert_24}
Maurice Reichert, \emph{Geometric invariants and asymptotics of translation
  surfaces}, Ph.D. thesis, Heidelberg University,
  \url{https://archiv.ub.uni-heidelberg.de/volltextserver/35709/}, 2024.

\bibitem[Sau18]{sauvaget_18}
Adrien Sauvaget, \emph{Volumes and {Siegel}-{Veech} constants of
  {{\({\mathcal{H}}(2G - 2)\)}} and {Hodge} integrals}, Geometric and
  Functional Analysis \textbf{28} (2018), no.~6, 1756--1779.

\bibitem[Sau21]{sauvaget_21}
\bysame, \emph{The large genus asymptotic expansion of {Masur}-{Veech}
  volumes}, International Mathematics Research Notices \textbf{2021} (2021),
  no.~20, 15894--15910.

\bibitem[Sie45]{siegel_45}
Carl~Ludwig Siegel, \emph{A mean value theorem in geometry of numbers}, Annals
  of Mathematics. Second Series \textbf{46} (1945), 340--347.

\bibitem[Val24]{vallejos_24}
Hunter Vallejos, \emph{Random geometric structures on high genus surfaces},
  Ph.D. thesis, University of Texas at Austin,
  \url{https://hunter-vallejos.com/s/Hunter\_Dissertation\_Final.pdf}, 2024.

\bibitem[Vee82]{veech_82}
William~A. Veech, \emph{Gauss measures for transformations on the space of
  interval exchange maps}, Annals of Mathematics. Second Series \textbf{115}
  (1982), 201--242.

\bibitem[Vee98]{veech_98}
\bysame, \emph{Siegel measures}, Annals of Mathematics. Second Series
  \textbf{148} (1998), no.~3, 895--944.

\bibitem[Vor05]{vorobets_05}
Yaroslav Vorobets, \emph{Periodic geodesics on generic translation surfaces},
  Algebraic and topological dynamics. Proceedings of the conference, Bonn,
  Germany, May 1--July 31, 2004, Providence, RI: American Mathematical Society
  (AMS), 2005, pp.~205--258.

\bibitem[YZZ20]{yang_zagier_zhang_20}
Di~Yang, Don Zagier, and Youjin Zhang, \emph{Masur-{Veech} volumes of quadratic
  differentials and their asymptotics}, Journal of Geometry and Physics
  \textbf{158} (2020), 12, Id/No 103870.

\bibitem[Zor02]{zorich_02}
Anton Zorich, \emph{Square tiled surfaces and {Teichm{\"u}ller} volumes of the
  moduli spaces of {Abelian} differentials}, Rigidity in dynamics and geometry.
  Contributions from the programme Ergodic theory, geometric rigidity and
  number theory, Isaac Newton Institute for the Mathematical Sciences,
  Cambridge, UK, January 5--July 7, 2000, Springer, 2002, pp.~459--471.

\end{thebibliography}
	
\end{document}